\documentclass[12pt,a4paper,twoside]{amsart}
\usepackage{geometry}
\usepackage{fancyhdr}
\usepackage{txfonts}
\usepackage{bbm}
\usepackage{threeparttable}
\usepackage{mathrsfs}
\usepackage{amsfonts}
\usepackage{amsmath}
\usepackage{amssymb}
\usepackage{wasysym}
\usepackage{psfrag}
\usepackage{graphics,epsfig,amsmath}  
\usepackage{comment}

\usepackage{xcolor}
\usepackage{enumerate}
\usepackage{hyperref}
\usepackage[bbgreekl]{mathbbol}
\usepackage{bm}

\DeclareSymbolFontAlphabet{\mathbb}{AMSb}
\DeclareSymbolFontAlphabet{\mathbbl}{bbold}

\newtheorem{theorem}{Theorem}[section]

\newtheorem{corollary}{Corollary}

\newtheorem{proposition}{Proposition}

\newtheorem{lemma}{Lemma}[section]

\newtheorem{remark}{Remark}[section]

\newcommand{\ZZ}{\mathbb{Z}}
\newcommand{\NN}{\mathbb{N}}
\newcommand{\RR}{\mathbb{R}}
\newcommand{\CC}{\mathbb{C}}
\newcommand{\TT}{\mathbb{T}}
\newcommand{\mc}{\mathcal}

\newcommand{\red}{\textcolor{red}}
\newcommand{\blue}{\textcolor{blue}}

\newcommand{\orange}{\textcolor{orange}}

\newcommand{\scr}{\mathscr}

\allowdisplaybreaks

\numberwithin{equation}{section}

\title[linear stability of KAM tori]{On linear stability of KAM tori via the Craig-Wayne-Bourgain method}

\author[X. He]{Xiaolong He}
\address{
Department of mathematics\\
Hangzhou Normal University\\
Hangzhou, 311121, P.R. China}

\email{hexiaolong@hnu.edu.cn}

\author[J. Shi]{Jia Shi}
\address{
School of Mathematical Science\\
Fudan University\\
 Shanghai, 200433, P.R. China}
\email{15110180007@fudan.edu.cn}

\author[Y. Shi]{Yunfeng Shi}
\address{School of Mathematical Sciences\\
 Peking University\\
  Beijing, 100871, P.R. China}
\email{yunfengshi18@gmail.com}

\author[X. Yuan]{Xiaoping Yuan}
\address{
School of Mathematical Science\\
Fudan University\\
 Shanghai, 200433, P.R. China}
\email{xpyuan@fudan.edu.cn}

\date{\today}

\begin{document}

\vspace{0.1in}

\begin{abstract}
In this paper, we prove the Melnikov's persistency theorem by combining the
traditional Kolmogorov-Arnold-Moser (KAM) technique and the Craig-Wayne-Bourgain (CWB)  method. The aim of this paper is twofold.
One is to establish the linear stability of the perturbed invariant tori by
using the CWB method without the second Melnikov condition. The other one is
to illustrate the CWB method in detail and make the CWB method more accessible.
\end{abstract}


\maketitle
\tableofcontents

\section{Introduction}

\subsection{Background}
The  celebrated Kolmogorov-Arnold-Moser (KAM) theory concerns with the stable motions
in nearly integrable Hamiltonian system.
For a smooth Hamiltonian of $n$-degree
\begin{equation*}
  H=H_{0}(I)+\epsilon R(\theta,I)
\end{equation*}
with the standard symplectic structure
$d\theta\wedge d I$ on $\TT^{n}\times\RR^{n}$
and the angle-action variable $(\theta,I)$
belongs to some domain in $\TT^{n}\times D\subset\TT^{n}\times \RR^{n}$.
Assume the unperturbed Hamiltonian $H_{0}(I)$
is independent of $\theta$ and satisfies
the Kolmogorov  non-degenerate condition
\begin{equation*}
  \det (\partial^{2}{H_{0}}/{\partial I^{2}})\neq 0,\quad I\in D.
\end{equation*}
Then  the invariant torus
$\mc{T}_{0}=\{\omega(I_{0}) t: t\in\RR\}\times \{I_{0}\}$
with a prescribed Diophantine frequency
$\omega(I_{0})=\partial H (I_{0})/\partial I$
for some $I_{0}\in D$
persists under sufficient small perturbation
$\epsilon R(\theta, I)$.
This is the well-known KAM theorem (\cite{Kol54,Arn61, Mos62})
for the finitely dimensional Hamiltonian system.
It is worthy mentioning that the dimension
of the persisted torus equals to the degree
of the Hamiltonian system.

To explore the existence of those invariant tori
whose dimensions are less than the degree of
the Hamiltonian, consider the following
Hamiltonian
\begin{equation*}
  E= \langle\omega,y\rangle+\sum_{j=1}^{n}
  \Omega_{j} z_{j} \bar{z}_{j}
\end{equation*}
defined on the phase space $(x,y,z,\bar{z})\in \TT^{d}\times \RR^{d}\times \CC^{n}\times\CC^{n}$ with the symplectic structure
$d x\wedge d y+ \sqrt{-1} d z\wedge d \bar{z}$.
Obviously, one finds that $
\mathcal{T}^{d}_{0}= \{\omega t: t\in\RR\}\times \{y=0\}\times\{z=0\}\times \{\bar{z}=0\}$
is an invariant torus of the Hamiltonian
vector field $X_{E}$. Moreover, the dimension
of $\mc{T}^{d}_{0}$ is less than the
degree $d+n$ of the Hamiltonian $E$.
In 1965 , Melnikov \cite{Mel68} announced that, under
suitable non-resonant conditions, the lower
dimensional invariant tori can persist under
sufficiently small Hamiltonian perturbation $\epsilon R(x,y, z,\bar{z})$. In the late 1980's,
Eliasson \cite{Eli88}, P\"oschel \cite{Pos89},
Kuksin \cite{Kuk93} provided a complete proof
of the problem, well-known nowadays as
Melnikov's persistency theorem.

We briefly explain the main idea of the proof
in \cite{Eli88, Pos89, Kuk93}.
Roughly speaking, we expand the
perturbation $R$ into Taylor series in
$(y,z,\bar{z})$
\begin{equation}
\begin{aligned}
  R= & \underbrace{R^{x}(x)+ \langle R^{y}(x), y\rangle}_{(I)}
  + \underbrace{\langle R^{z}(x), z\rangle+\langle R^{\bar{z}}(x)
    ,\bar{z}\rangle}_{(II)}
     +  \underbrace{\langle R^{zz}(x) z,z \rangle+ \langle R^{z\bar{z}}(x) z, \bar{z}\rangle
    +\langle R^{\bar{z}\bar{z}}(x)\bar{z},\bar{z}\rangle}
  _{(III)}\\
  & + O(|y|^{2}+|y|\cdot|z|+|z|^{3}),
  \end{aligned}
\end{equation}
Then we apply the symplectic transformation
to eliminate the items $(I), (II)$ and
$(III)$ of lower order.
Consequently, $\mc{T}_{0}^{d}$ is
an invariant torus of the  Hamiltonian
$H=E+O(|y|^{2}+|y|\cdot |z|+|z|^{3})$.

To eliminate the item $(I)$, we need the usual
Diophantine condition
\begin{equation*}
  \langle k,\omega\rangle\neq 0 \quad
  \textrm{for all}~0\neq k\in\ZZ^{d}.
\end{equation*}
To eliminate the item $(II)$, we need
the first Melnikov condition
\begin{equation*}
  \Omega_{j}\pm\langle k,\omega\rangle\neq 0,\quad \textrm{for all}~k\in\ZZ^{d}, 1\leq j\leq n.
\end{equation*}
To eliminate the expressions $\langle R^{zz}(x)z,z\rangle$ and $\langle R^{\bar{z}\bar{z}}(x)\bar{z},\bar{z}\rangle$ in
$(III)$, we need the second Melnikov condition
\begin{equation}\label{003}
  \Omega_{i}+\Omega_{j}\pm\langle
  k,\omega\rangle\neq 0,\quad
  \textrm{for all}~k\in\ZZ^{d}, 1\leq i,j\leq n.
\end{equation}
To eliminate $\langle R^{z\bar{z}}z, \bar{z}\rangle$ in $(III)$, we still
need the second Melnikov condition but
in the following form
\begin{equation}\label{002}
  \Omega_{i}-\Omega_{j}\pm \langle
  k,\omega\rangle\neq 0,\quad\textrm{for all}~
  |i-j|+|k|\neq 0, k\in\ZZ^{d}, 1\leq i,j\leq n.
\end{equation}
We see from \eqref{002} that when
$k=0$ there is $\Omega_{i}\neq \Omega_{j}$
for any $i\neq j$, i.e.,
the multiplicity of the norm frequency
should be one, which excludes lots of important
applications.

In 1997,
Bourgain \cite{ Bou97_MRL}  improved
 Craig-Wayne's method  \cite{CW93}
to study the Melnikov's problem,
which is completely free from
the second Melnikov condition
and also applies
to infinitely dimensional Hamiltonian system
\cite{Bou97_MRL, Bou98_Ann}. In his famous book
\cite{Bou05_Green} published in 2005,
 Bourgain developed further the method to prove the existence of invariant torus (or
quasi-periodic solution) for NLS
and NLW of arbitrary dimension.
This method now is known as the Craig-Wayne-Bourgain (CWB)
method.
The CWB method is less dependent on the Hamiltonian structure.
It is essentially  based on applying the Newton iteration
to solve directly the differential equation for
quasi-periodic solutions.
However, one has to pay for the price that
the homological equation (or the linearized equation) not only has small divisor problem but also contains variable coefficients.
Moreover, we are not able to obtain a local norm form
around the persisted invariant torus.

Back to the Melnikov's problem in \cite{Bou97_MRL}, Bourgain combined the above
CWB method with the KAM technique.
Taking Taylor expansion of the perturbation and
applying the symplectic transformation as
before, Bourgain
put $(III)$ into unperturbed Hamiltonian
$E$ and eliminated sorely $(I)$ and $(II)$,
which results in a norm form around the invariant
torus
\begin{equation*}
  H_{\infty}=E_{\infty}+O(|y|^{2}+|y|\cdot |z|+|z|^{3})
\end{equation*}
with
\begin{equation*}
  E_{\infty}=\langle \omega, y\rangle
  +\sum_{j=1}^{n}\Omega_{j}z_{j}\bar{z}_{j}
  +(III).
\end{equation*}
Obviously, $\mc{T}_{0}^{d}$ is still an invariant torus of $H_\infty$. The important thing is that since $(III)$ has been putted into  $E_\infty$, it avoids completely the usage of the second Melnikov condition. However, to derive such a normal form $H_\infty$, we have to solve homological equations with variable coefficients. Moreover, the linear stability of the persisted torus is unknown.

Note that we can actually divide the second Melnikov conditions into two parts \eqref{003} and \eqref{002} with \eqref{003} containing terms of the form $\Omega_i+\Omega_j$. Apparently, part \eqref{003} has  essentially the form of the first Melnikov condition. The true difficulty arises from part \eqref{002} which involves the terms $\Omega_i-\Omega_j$. Thus we can eliminate terms associated with part \eqref{003}, and put terms corresponding to part \eqref{002} into the new normal form. As a result, we may obtain a more precise normal form
\begin{align*}
H_\infty'=E_\infty'+O(|y|^2+|y||z|+|z|^3)
\end{align*}
with
\begin{align*}
E_\infty'=\langle\omega,y\rangle+\sum_{j}\Omega_j
z_j\bar z_j+\langle R^{z\bar z}(x)z,\bar z\rangle.
\end{align*}
In particular, $\mc{T}_{0}^{d} =\{\omega t: t\in\RR\}\times \{0\}\times\{0\} \times \{0\}$ is also an invariant torus of $H_\infty'$. Furthermore, the corresponding linearized equation
\begin{align*}
\sqrt{-1}\dot{z}=\Lambda z+R^{z\bar z}(x)z, \quad \Lambda=\mathrm{diag}(\Omega_j)
\end{align*}
admits a $L^2$-conservation law, which implies particularly the linear stability of persisted torus (see Theorem \ref{main th} and Corollary \ref{cor}
 in the following for details).

The aim of this paper is twofold. One is to study the Melnikov's problem by combining the CWB method and
the KAM technique. We show that the linear stability still holds without
the second Melnikov condition \eqref{002}.
The other one is to explain the CWB method in detail and  to make the CWB method more accessible.


\begin{remark}
An alternative method is to put $\langle [R^{z\bar z}(x)]z,\bar z\rangle$ into $E$, where $[R^{z\bar z}(x)]=\int R^{z\bar z}(x){\rm d}x$. This method has the advantage that the homological equations are  of constant coefficients type. The disadvantage  is that the second Melnikov conditions are still  employed, which seems not applicable to higher spatial dimensional  NLS and NLW in infinitely dimensional systems case. This method can be found in an early monograph \cite{BMS76} published in 1969 in Russian. See also
\cite{You99}.

\end{remark}

\subsection{Main result}

Let us recall some basic concepts in the Hamiltonian dynamical systems.
Consider a Hamiltonian function $H=H(x,y,z,\bar{z})$ defined on the phase space
$\mc{P}=\TT^{d}\times \RR^{d}\times \CC^{n}\times\CC^{n}
$ with $\TT^{d}=\RR^{d}/(2\pi \ZZ)^{d}$.
We endow the symplectic form
\begin{equation*}
  \pmb{\omega}=d x\wedge d y+ \sqrt{-1} d z\wedge d \bar{z}=
  \sum_{j=1}^{d} d x_{j}\wedge d y_{j}
  +\sqrt{-1} \sum_{k=1}^{n} d z_{k}\wedge d \bar{z}_{k}.
\end{equation*}
Then the vector field $X_{H}$ given by $X_{H}\lrcorner~ \pmb{\omega}=-d H$ reads
\begin{equation*}
  X_{H}=(\partial_{y} H, -\partial_{x} H,
  \sqrt{-1} \partial_{\bar{z}} H,
  -\sqrt{-1} \partial_{z} H)^{T}.
\end{equation*}
The associated Poisson bracket takes the form of
\begin{equation*}
  \{F,G\}= \langle F_{x}, G_{y}\rangle- \langle F_{y}, G_{x}\rangle
  +\sqrt{-1} \langle F_{z}, G_{\bar{z}}\rangle-\sqrt{-1}
  \langle F_{\bar{z}}, G_{z}\rangle.
\end{equation*}
Given a function $F$, the time-$1$-map of the flow
$X_{F}^{t}$ of the  Hamiltonian vector field $X_{F}$ is symplectic. Moreover,
\begin{equation*}
  \frac{d}{dt} G\circ X_{F}^{t}= \{G,F\}\circ X_{F}^{t}.
\end{equation*}

In this paper, we consider small perturbation of a finite dimensional Hamiltonian
in the parameter dependent normal form
\begin{equation*}
  {E}_{0}= \langle \omega_{0}(\xi),y\rangle+ \langle \Omega z, \bar{z}\rangle,
  \quad (x,y,z,\bar{z})\in \TT^{d} \times \RR^{d}
  \times \CC^{n}\times\CC^{n},
\end{equation*}
where $\Omega=\textrm{diag}(\Omega_{j}: 1\leq j\leq n)$ with $\Omega_{j}>0$.
The tangent frequency $\omega_{0}$
depends on $d$ parameters $\xi\in\Pi_{0}\subset\RR^{d}$,
where $\Pi_{0}$ is a given open set.
The associated Hamiltonian vector field $X_{{E}_{0}}$ of
the normal form ${E}_{0}$
is given by
\begin{equation*}
  X_{{E}_{0}}= (\omega_{0}(\xi), 0, \sqrt{-1} \Omega z, -\sqrt{-1} \Omega \bar{z} )^{T},
\end{equation*}
where $(\cdot)^{T}$ represents the transpose of a matrix (or a vector).
Obviously, for each $\xi\in\Pi_{0}$, there is a $d$-dimensional invariant
torus
\begin{equation*}
  \mathcal{T}^{d}_{0}= \TT^{d}\times \{y=0\}\times\{z=0\}\times \{\bar{z}=0\},
\end{equation*}
carrying a quasi-periodic flow $x= \omega_{0} t+x_{0}$ with fixed torus
frequency $\omega_{0}=\omega_{0}(\xi)$.

The Melnikov's problem is to study the persistence of $\mathcal{T}^{d}_{0}$
under sufficiently small perturbation of the Hamiltonian. We consider
perturbation
\begin{equation*}
  H= {E}_{0}+ P_{0}
\end{equation*}
of ${E}_{0}$ that are \emph{real analytic}\footnote{The real analyticity of $H$ means that $H$ is analytic on the
complex domain $\scr{D}(s,r)$, and takes real value when
$x,y$ are real and $z,\bar{z}$ are complex conjugated. } on some
complex neighborhood
\begin{equation}
  \mathscr{D}(s,r): \quad |\textrm{Im} x|_{\infty}< s,\quad |y|< r^{2}, \quad
  |z|< r,\quad |\bar{z}|<r
\end{equation}
of $\mathcal{T}^{d}_{0}$ in the complex space
$\mathcal{P}_{\CC}=(\CC^{d}/2\pi \ZZ^{d})\times \CC^{d}\times\CC^{n}\times\CC^{n}$, where
$|\cdot|_{\infty}$ denotes the supremum norm and $|\cdot|$ denotes
the Euclidean norm.
It should be pointed out that $z$ and $\bar{z}$ are independent variables.
We also introduce
\begin{equation*}
  \mathscr{D}_{\RR}(s,r)=\{(x,y,z,\bar{z})\in\mathscr{D}(s,r): x,y\in\RR^{d} \},
\end{equation*}
in which $x,y$ are real but $z$ and $\bar{z}$  stay in
the complex space and are complex conjugated.

For $r>0$, we define the weighted phase norm
\begin{equation*}
  _{r}\pmb{|}W\pmb{|}= |X|+\frac{1}{r^{2}}|Y|+\frac{1}{r}|Z|+\frac{1}{r}|\bar{Z}|,\quad
  \textrm{for}~W=(X,Y,Z,\bar{Z})\in \mathcal{P}_{\CC}.
\end{equation*}
For a map $W: \mathscr{D}(s,r)\times \scr{O} \rightarrow \mathcal{P}_{\CC}$, define
\begin{equation*}
  _{r}\pmb{|}W\pmb{|}_{\mathscr{D}(s,r)\times \scr{O}}=
  \sup_{(u,\xi)\in \mathscr{D}(s,r)\times
  \scr{O}} ~
  _{r}\pmb{|}W(u,\xi)\pmb{|}
\end{equation*}
and
\begin{equation*}
  _{r}\pmb{|}W\pmb{|}^{\mathcal{L}}_{\mathscr{D}(s,r)\times \scr{O}}=\sup_{(u,\xi)\in
  \mathscr{D}(s,r)\times
  \scr{O}} ~
  _{r}\pmb{|}\partial_{\xi} W(u,\xi)\pmb{|},
\end{equation*}
where $\partial_{\xi}$ is the derivative with respect to $\xi$ and $\scr{O}\subset\CC^{d}$ is
an open set.

%


Now we state the basic assumptions on the Melnikov's problem.

\bigskip

\noindent\textbf{Assumption A} (Analyticity w.r.t. parameters).
Assume that  $\omega_{0}$ is \emph{real analytic}\footnote{
We say a function is real analytic on some domain in $\CC^{d}$ when it is analytic on that domain and is real for real arguments.
} in $\xi$ on $\scr{O}_{0}\subset\CC^{d}$,
where $\scr{O}_{0}=\scr{O}(\Pi_{0},\rho_{0})=\{z\in\CC^{d}: |z-\xi|<\rho_{0}~\textrm{for some}~\xi\in\Pi_{0}\}$ and $\Pi_{0}\subset\RR^{d}$ is an open interval. When saying an open interval in $\RR^{d}$, we always mean any open set
of the form $\{(\xi_{1},\cdots,\xi_{d}): a_{j}<\xi_{j}<b_{j}, 1\leq j\leq d)\}$.
\bigskip

\noindent\textbf{Assumption B} (Non-degeneracy).
There is some absolute constant $C>0$ such that
\begin{equation*}
  \sup_{\xi\in\mathscr{O}_{0}}|\partial_{\xi}\omega|< C,\quad
  \sup_{\xi\in\mathscr{O}_{0}}|\partial_{\xi} \omega^{-1}|<C.
\end{equation*}

\bigskip

\noindent\textbf{Assumption C} (Regularity).
Let $s_{0}, r_{0}$ be positive constants. Assume the perturbation
$P_{0}(x,y,z,\bar{z};\xi)$  is real analytic in $(x,y,z,\bar{z})$ on the domain
$\mathscr{D}(s_{0},r_{0})$.
For each $\xi\in\mathscr{O}_{0}$, the Hamiltonian
vector field
\begin{equation*}
X_{P_{0}}= (\partial_{y} P_{0}, -\partial_{x} P_{0}, \sqrt{-1} \partial_{\bar{z}} P_{0},
 -\sqrt{-1} \partial_{z}P_{0} )^{T}
\end{equation*}
defines near $\mathcal{T}^{d}_{0}$ an analytic map
\begin{equation*}
  X_{P_{0}}: \scr{D}(s_{0}, r_{0})\subset \mathcal{P}_{\CC}\rightarrow \mathcal{P}_{\CC}.
\end{equation*}
Also assume that $X_{P_{0}}$ is real analytic in
$\xi\in\mathscr{O}_{0}$.
\bigskip

\noindent\textbf{Assumption D} (Reality).
For any $(x,y,z,\bar{z},\xi)\in\mathscr{D}_{\RR}(s_{0},r_{0})
\times\Pi_{0}$, the perturbation $P_{0}$
satisfies the reality condition, i.e.,
\begin{equation*}
  \overline{P_{0}(x,y,z,\bar{z},\xi)}= P_{0}(x,y,z,\bar{z},\xi),
\end{equation*}
where the overline denotes the complex conjugate.
\bigskip

\begin{theorem}\label{main th}
  Suppose $H= {E}_{0}+P_{0}$ satisfies \emph{\textbf{Assumptions A-D}}
  and assume the  smallness condition
  \begin{equation*}
    {}_{r_{0}}
    \pmb{|}X_{P_{0}}\pmb{|}_{\scr{D}(s_{0},r_{0})\times \mathscr{O}_{0}}< \epsilon,\quad
    {}_{r_{0}}\pmb{|}X_{P_{0}}\pmb{|}_{\scr{D}(s_{0},r_{0})\times \mathscr{O}_{0}}^{\mathcal{L}}<\epsilon^{1/3}.
  \end{equation*}
  Then there is a sufficiently small $\epsilon_{*}=\epsilon_{*}(n,d,r_{0},s_{0},\rho_{0},\Pi_{0})>0$
   such that for any $0<\epsilon<\epsilon_{*}$, there
  is a subset $\Pi_{\infty}\subset\Pi_{0}$ with
  \begin{equation*}
    \emph{mes}~(\Pi_{\infty})\geq (1-{O}(\epsilon^{1/2}))~ \emph{mes}~ (\Pi_{0}),
  \end{equation*}
  and there are a family of embedding $\Phi: \TT^{d}\times \Pi_{\infty}\rightarrow
  \mc{P}$,  a map $\omega_{*}: \Pi_{\infty}\rightarrow \RR^{d}$ and
  a matrix function $B^{z\bar{z}}: \TT^{d}\times \Pi_{\infty}\rightarrow \RR^{n\times n}$ such that
  for each $\xi\in\Pi_{\infty}$, the transformation $\Phi$
  and the matrix $B^{z\bar{z}}$ are real analytic on $\TT^{d}_{s_{0}/2}$ giving rise to
  \begin{equation*}
    H\circ\Phi|_{\TT^{d}\times\{\xi\}}= \langle \omega_{*}(\xi),y\rangle
    +\langle \Omega z,\bar{z}\rangle+ \langle B^{z\bar{z}}(x;\xi) z, \bar{z}\rangle
    + {O}(|y|\cdot|z|+|y|^{2}+|z|^{3}).
  \end{equation*}

\end{theorem}

From Theorem
\ref{main th}, one readily see that, for each $\xi\in\Pi_{\infty}$,
 the vector $ \Phi_{\TT^{d}\times \{\xi\}}$ is an analytic embedding
 of rotational torus with frequency $\omega_{*}(\xi)$ for the Hamiltonian
 $H$ at $\xi$. Moreover, following the analysis of \eqref{linear eq}, we further
 obtain the linear stability of the invariant torus.

\begin{corollary}\label{cor}
  Under the assumptions of Theorem \ref{main th}, the perturbed invariant tori
  are linearly stable in the sense
  that the associated Lyapunov exponent is zero.
\end{corollary}

\noindent\textbf{Proof.} For each $\xi\in\Pi_{\infty}$, we consider
the Hamiltonian vector field induced by the Hamiltonian $H\circ \Phi$.
We immediately find  that
\begin{equation*}
  \mathcal{T}^{d}_{0}=\TT^{d}\times \{y=0\}\times \{z=0\}\times\{\bar{z}=0\}
\end{equation*}
is a $d$-dimensional invariant torus of the vector field $X_{H\circ \Phi}$.
Then the linearized equation around $\mathcal{T}^{d}_{0}$ is
\begin{equation}\label{linear eq}
  \left\{
  \begin{aligned}
    \dot x=& \omega,\\
    \dot y=& 0,\\
    \dot z=& \sqrt{-1}(\Omega+B^{z\bar{z}}(x)) ~z,\\
    \dot {\bar{z}}=& -\sqrt{-1} (\Omega+B^{z\bar{z}}(x))~ \bar{z}.
  \end{aligned}
  \right.
\end{equation}
 Along the trajectory $z=z(t)$ of \eqref{linear eq}, we have the $L^{2}$-conservation, i.e.,
\begin{equation*}
  \frac{\textrm{d}}{\textrm{d} t} |z|^{2}= \frac{\textrm{d}}{\textrm{d} t} \langle z,\bar{z}
  \rangle=0,
\end{equation*}
which implies $(z,\bar{z})=0$ is a center equilibrium in \eqref{linear eq}.
This proves the linear stability of the perturbed invariant torus.
\qed

\section{The KAM Iterative Lemma}

In this section, we establish the KAM Iterative Lemma, upon which our main Theorem
\ref{main th} is an immediate result.
To begin with, we summarize the notations
and the iterative constants in subsection
\ref{notation} for reader's quick reference. Next we present
and prove  the KAM Iterative Lemma in
subsection \ref{sect iter lemma} and subsection
\ref{proof iter l}, respectively.
In subsection \ref{proof main}, we prove our main
Theorem \ref{main th}.

\subsection{Notations and the iterative constants}
\label{notation}

We first introduce some  general notations. For two vectors $a, b$ in $\RR^{d}$ or $\CC^{n}$, we denote
$\langle a,b\rangle=\sum_{j} a_{j} b_{j}$. We use the notation
$A\setminus B$ for the set theoretical difference. For
$k\in\ZZ^{d}$ and $U\subset\ZZ^{d}$, $k+U$ denotes the set
$\{k'=k+p: p\in U\}$.
The symbols $\wedge$ and $\vee$ describes
the minimal and maximal operators, respectively.
The measure of a set $\scr{V}\subset\RR^{d}$,
denoted by $\textrm{mes} (\scr{V})$, always
refers to the Lebesgue measure.
By some abuse of notation, we denote by
$|\scr{J}|$ the diameter of a set $\scr{J}\subset\RR^{d}$.

Following the notations in KAM theory, we denote in the
sequel various constants by the same letter $C$.
Of course, these numbers depend only on the universal constants
$d,n,\rho_{0},r_{0},s_{0},\Pi_{0}$ and could be made
explicit  by the  context where they arise, but need not be. For further simplicity,
we write $a\lesssim b$ in estimates to suppress the multiplicative constant in
$Ca< b$. The notation $a\ll b$ indicates $Ca<b$ for sufficiently
large $C>0$ and $a\sim b$ means both $a\ll b$ and $b\ll a$
hold. Furthermore, $\varepsilon^{1-}$ means
$\varepsilon^{1-\delta}$ with some small $\delta>0$ (
the precise meaning of "small" can again be derived
from the context), in which the exponent "$1-$" might be different from line
to line.

 If not specified,
the norm  for vectors in real or complex space
 refers to the Euclidean norm. The norm of a matrix is
 the induced operator norm on the vectors.
 For a Fourier series $q(x)=\sum_{k\in\ZZ^{d}} \widehat{q}(k)
 e^{\sqrt{-1} \langle k,x\rangle}$, we define the
 truncation operator $\Gamma_{N}$ by
 \begin{equation}\label{trunc}
	(\Gamma_{N}q)(x)=\sum_{k\in\ZZ^{d}, |k|\leq N}
	\widehat{q}(k)
	e^{\sqrt{-1}~ \langle k, x\rangle }.
 \end{equation}

Next,
we define
the following iterative
constants and domains:
\begin{itemize}
  \item $s_{0}>0$ and  $r_{0}>0$ are fixed and  given in \textbf{Assumption C};

  \item
  $A=A(n,d)>0$ is sufficiently large;

  \item $l\in\mathbb{N}$ is the number of the  KAM iterative steps;

  \item $\epsilon_{l}=A^{-(\frac{4}{3})^{l}}$ measures the size of the perturbation at the
  $l^{th}$ step;

  \item $e_{l}=\frac{1^{-2}+2^{-2}+\cdots+l^{-2}}{2(1^{-2}+2^{-2}+\cdots)}$ (so
  $0<e_{l}<\frac{1}{2}$ for all $l$);

  \item $s_{l}=s_{0}(1-e_{l})$ (so $s_{l}>\frac{1}{2}s_{0}$ for all $l$), which measures
  the width of the analyticity strip for the angle variable $x$ at the $l^{th}$ step;

  \item $r_{l}=r_{0}(1-e_{l})$ (so $r_{l}>\frac{1}{2} r_{0}$ for all $l$), which
  measure the analyticity radius for the action variable $y$, as well as the normal
  coordinates $z,\bar{z}$, at the $l^{th}$ step;

  \item $s_{l}^{{(j)}}=(1-\frac{j}{100}) s_{l}+ \frac{j}{100} s_{l+1}$ $(j=0,\cdots,100)$
  are the intermediate points between $s_{l}$ and $s_{l+1}$ dividing
  $[s_{l},s_{l+1}]$ into $100$ subintervals with the same length;

  \item $r_{l}^{{(j)}}=(1-\frac{j}{100}) r_{l}+ \frac{j}{100} r_{l+1}$ $(j=0,\cdots,100)$
  are the intermediate points between $r_{l}$ and $r_{l+1}$ dividing
  $[r_{l},r_{l+1}]$ into $100$ subintervals with the same length;

  \item $\scr{D}(s_{l},r_{l})=\{(x,y,z,\bar{z})\in\mathcal{P}_{\CC}: |\textrm{Im} x|_{\infty}<s_{l},
  |y|<r_{l}^{2}, |z|<r_{l}, |\bar{z}|<r_{l}\}$ denotes a neighborhood of the torus
  \begin{equation*}
    \mc{T}^{d}_{0}=\TT^{d}\times \{y=0\}\times\{z=0\}\times\{\bar{z}=0\},
  \end{equation*}
  where $|\cdot|_{\infty}$ denotes the supremum norm. Obviously,
  \begin{equation*}
    \scr{D}(s_{l},r_{l})\supset \scr{D}(s_{l+1},r_{l+1})\supset\cdots\supset
    \scr{D}(\frac{s_{0}}{2}, \frac{r_{0}}{2});
  \end{equation*}

  \item $\TT^{d}_{s_{l}}=\{x\in\CC^{d}/(2\pi \ZZ)^{d}: |\textrm{Im} x|_{\infty}< s_{l} \}$
  denotes a neighborhood of $\TT^{d}$ with strip width $s_{l}$ and obviously
  \begin{equation*}
    \TT^{d}_{s_{l}}\supset\TT^{d}_{s_{l+1}}\supset\cdots\supset\TT^{d}_{s_{0}/2};
  \end{equation*}

  \item Given a sequence of open sets  $\Pi_{l}$ in $\RR^{d}$, we denote
  \begin{equation*}
    \scr{O}_{l}=\scr{O}(\Pi_{l}, A^{-l^{C_{3}}})=\{\xi\in\CC^{d}: |\xi-\xi'|<A^{-l^{C_{3}}}
    ~\textrm{for some}~\xi'\in\Pi_{l}\}\subset\CC^{d}.
  \end{equation*}

\end{itemize}

Finally, we define some matrices depending on the iteration, which are used to solve
the homological equations. At the $l$-th KAM iterative step, we define the matrix
\begin{equation}\label{T l 0}
  T_{l}=D_{l}+S_{l}
\end{equation}
defined on $\{1,\cdots,n\}\times \ZZ^{d}$,
where $D_{l}$ is a diagonal matrix
\begin{equation}\label{T l 1}
  D_{l}(j,k)= \Omega_{j}+ \langle k,\omega_{l}\rangle
\end{equation}
and $S_{l}$ is a non-diagonal matrix
\begin{equation}\label{T l 2}
  S_{l}((j,k),(j',k'))= (B_{l;jj'}+R^{z\bar{z}}_{l;jj'})^{\wedge}(k-k'),
\end{equation}
with $1\leq j,j'\leq n$ and $k,k'\in\ZZ^{d}$.
The matrix-valued functions $B_{l}$ and $R^{z\bar{z}}_{l}$ are defined in
the following Iterative Lemma and the
hat $\widehat{(\cdot)}(k)$ (or $(\cdot)^{\wedge}(k)$) denotes
the $k$-th Fourier coefficient of the associated function.

For $U\subset\ZZ^{d}$, we denote by $T_{l;U}$ the
restriction of the matrix $T_{l}$ on $\{1,\cdots,n\}\times U$, i.e.,
\begin{equation*}
T_{l;U}((j,k),(j',k'))	=T_{l}((j,k),(j',k')),
\end{equation*}
when $(j,k),(j',k')\in\{1,\cdots,n\}\times U$.
By some abuse of notation, we sometimes denote by
$T_{l;M}$ the restriction of $T_{l}$ on
$\{1,\cdots,n\}\times ([-M,M]^{d}\cap \ZZ^{d})$ for
 any integer $M>0$. The inverse matrix of
$T_{l; M}$ (or $T_{l;\Lambda}$), if exists, is always denoted by $G_{l;M}$ (or $G_{l;\Lambda}$).

In addition, we define another matrix $\textbf{T}_{l}$ on $\{\textbf{j}=(j_{1},j_{2}):1\leq j_{1},j_{2}\leq n \}\times \ZZ^{d}$
by $\textbf{T}_{l}=\textbf{D}_{l}+\textbf{S}_{l}$.
The diagonal matrix $\textbf{D}_{l}$ is defined by
$\textbf{D}_{l}(\textbf{j},k)= \mathbf{\Omega}_{\textbf{j}}
  +\langle k,\omega_{l}\rangle$ and $\mathbf{\Omega}_{\textbf{j}}
 =\Omega_{j_{1}}+\Omega_{j_{2}}.$ The nondiagonal matrix $\textbf{S}_{l}$ is
 defined  in \eqref{homo 4 1-3}. As we shall see later,
 $T_{l}$ and $\textbf{T}_{l}$ have essentially the same structure except
 the difference between the finite indices $j$ and $\textbf{j}$.
 Similarly, $\textbf{T}_{l;M}$ denotes the restriction of
 $\textbf{T}_{l}$ on $\{\textbf{j}=(j_{1},j_{2}): 1\leq j_{1},j_{2}\leq n\}\times ([-M,M]^{d}\cap\ZZ^{d})$
 and $\textbf{G}_{l;M}$ is the inverse of $\textbf{T}_{l;M}$ whenever the matrix
 $\textbf{T}_{l;M}$ is invertible.

\bigskip

\subsection{The Iterative Lemma}\label{sect iter lemma}

Choose and fix
the various constants $C_{0},C_{1},\cdots, C_{7}$  such that
  \begin{equation}\label{C}
  \begin{array}{llll}
& C_{1}>C_{0}\gg 1,\quad & C_{2}>2C_{1}+10,\quad & C_{4}>C_{3}>C_{1},\\
& C_{5}>C_{6}+2,
    \quad & C_{6}>2 C_{4},\quad & C_{7}> (C_{4}+10)\vee C_{5}.
    \end{array}
  \end{equation}

Unlike the usual KAM theorems, the following Iterative Lemma
starts from $l_{*}$ with $l_{*}=l_{*}(\epsilon)$ large enough.
To keep the consistency of the notations, we set
\begin{equation}\label{initial func}
H_{l_{*}}=H_{0}, ~ P_{l_{*}}=P_{0},~ B_{l_{*}}=B_{l_{*}-1}=0,~\Pi_{l_{*}-1}=\Pi_{0},~s_{l_{*}}=s_{0},~r_{l_{*}}=r_{0},~
\omega_{l_{*}}=\omega_{l_{*}-1}=\omega_{0}.	
\end{equation}

\begin{lemma}\label{iter lemma}

Consider a family of Hamiltonian functions $H_{l}~ (l_{*} \leq l\leq m)$,
\begin{equation}\label{decom. H 1}
  H_{l}= {E}_{l}+ P_{l},
\end{equation}
defined on $\scr{D}(s_{l},r_{l})\times \scr{O}_{l}$ with $\scr{O}_{l}=\scr{O}(\Pi_{l}, A^{-l^{C_{3}}})$,
where $${E}_{l}= \langle \omega_{l}(\xi),y \rangle+\langle \Omega z,\bar{z}\rangle
 + \langle B_{l}(x) z,\bar{z}\rangle$$ is a normal form and
 the perturbation
 $$P_{l}=
 P_{l}^{\emph{low}}+P_{l}^{\emph{high}},\quad
 P_{l}^{\emph{high}}=O(|y|\cdot|z|+
|y|^{2}+|z|^{3}).$$
Assume the Hamiltonian $H_{l}$ and the parameter set $\Pi_{l}$ satisfy the following
properties.
\begin{enumerate}[$(l.1)$]
\item The frequency $\omega_{l}$ is real analytic on
$\scr{O}_{l}$ and
\begin{equation*}
  \sup_{\xi\in\scr{O}_{l}}\left\{
  |\partial_{\xi} \omega_{l}|,
   |\partial_{\xi} \omega^{-1}_{l}|
   \right\}
  \lesssim 1.
\end{equation*}
Furthermore, we have
\begin{equation*}
  \sup_{\xi\in\scr{O}_{l}}|\omega_{l}-\omega_{l-1}|< \epsilon_{l}^{1/10}.
\end{equation*}

\item The matrix $B_{l}$ is analytic in $x\in \TT_{s_{l}}^{d}$ and
real analytic in $\xi\in\scr{O}_{l}$. For any fixed $(x,\xi)\in\TT^{d}\times \Pi_{l}$,
 the matrix $B_{l}$ is real symmetry, i.e.,
\begin{equation*}
B_{l; jk}(x,\xi)= B_{l;kj}(x,\xi),\quad \overline{B_{l;jk}(x,\xi)}= B_{l;jk}(x,\xi),
\end{equation*}
in which the indices $j$ and $k$ indicate the row or column.
  Furthermore, we have
  \begin{equation*}
    \sup_{\TT_{s_{l}}^{d}\times \scr{O}_{l} }
 \left\{ \epsilon^{-1}\|B_{l}(x,\xi)\|, \epsilon^{-1/3}
 \|\partial_{\xi} B_{l}(x,\xi)\|\right\}\lesssim 1,
  \end{equation*}
  and
  \begin{equation*}
    \sup_{\TT_{s_{l}}^{d}\times \scr{O}_{l} } |B_{l}-B_{l-1}|< {\epsilon_{l}^{1/10}}.
  \end{equation*}

\item The perturbation $P_{l}$ is analytic
 on $\mathscr{D}(s_{l},r_{l})\times \scr{O}_{l}$ and satisfies the reality condition
\begin{equation*}
  \begin{aligned}
   \overline{ P_{l}(x,y,z,\bar{z};\xi)}= P_{l}(x,y,z,\bar{z};\xi),\quad
   \emph{for}~ (x,y,z,\bar{z},\xi)\in\mathscr{D}_{\RR}(s_{l},r_{l})\times \Pi_{l}.
  \end{aligned}
\end{equation*}
 The Hamiltonian vector field
$$
X_{P_{l}}=(\partial_{y} P_{l}, -\partial_{x} P_{l}, \sqrt{-1} \partial_{\bar{z}} P_{l}, -\sqrt{-1} \partial_{z} P_{l})^{T},
$$
defines an analytic map
\begin{equation*}
  X_{P_{l}}: \scr{D}(s_{l},r_{l})\subset \mathcal{P}_{\CC}\rightarrow \mathcal{P}_{\CC}.
\end{equation*}
and satisfies
\begin{equation*}
\begin{aligned}
& {~}_{r_{l}} \pmb{|} X_{P_{l}^{\emph{low}}}\pmb{|}
 _{\mathscr{D}(s_{l},r_{l})\times \scr{O}_{l}} < \epsilon_{l},\quad
  {~}_{r_{l}} \pmb{|} X_{P_{l}^{\emph{low}}}\pmb{|}
  _{\mathscr{D}(s_{l},r_{l})\times \scr{O}_{l}}^{\mathcal{L}}
< \epsilon_{l}^{1/3},\\
& {~}_{r_{l}} \pmb{|} X_{P_{l}^{\emph{high}}}\pmb{|}
 _{\mathscr{D}(s_{l},r_{l})\times \scr{O}_{l}}
 \lesssim \epsilon,\quad
  {~}_{r_{l}} \pmb{|} X_{P_{l}^{\emph{high}}}\pmb{|}
  _{\mathscr{D}(s_{l},r_{l})\times \scr{O}_{l}}^{\mathcal{L}} \lesssim \epsilon^{1/3}.
  \end{aligned}
\end{equation*}

  \item The parameter set $\Pi_{l}$ is the union of a collection $\Lambda_{l}$ of disjoint
  intervals $\scr{J}\subset \RR^{d}$ of size
  $A^{-l^{C_{3}}}$, i.e., $\Pi_{l}= \cup_{\scr{J}\in\Lambda_{l}}\scr{J}$ with
  $|\scr{J}|= A^{-l^{C_{3}}}$. Moreover, the following properties hold.
  \begin{enumerate}
    \item [$(l.4.1)$] For any interval $\scr{J}\in\Lambda_{l}$, there is
    a unique $\scr{J}'\in\Lambda_{l-1}$ such that $\scr{J}\subset\scr{J}'$.

    \item[$(l.4.2)$]  The parameter set $\Pi_{l}$ is contained in
   \begin{equation*}
   \begin{aligned}
	&\{\xi\in\RR^{d}: |\langle k,\omega_{l-1}(\xi)\rangle|>\sqrt{\epsilon}~ (1+2^{-(l-1)})~ |k|^{-\tau},
      0\neq |k|\leq A^{l}\}\\
    \cap~ & \Big\{\xi\in\RR^{d}: \|G_{l-1; A^{l}}\|<A^{(l\log A)^{C_{1}}}\\
&\qquad\emph{and}\quad |G_{l-1; A^{l}}(k,k')|<
    \exp\{-(s_{l-1}-(l\log A)^{-8}) |k-k'|\}~\emph{for}~
    |k-k'|>(l \log A)^{C_{2}}\Big\}\\
    \cap ~ & \big\{\xi\in\RR^{d}: \|\textbf{G}_{l-1;A^{l}}\|<A^{(l\log A)^{C_{1}}}\\
&\qquad\emph{and}\quad |\textbf{G}_{l-1; A^{l}}(k,k')|<
    \exp\{-(s_{l-1}-(l\log A)^{-8}) |k-k'|\}~\emph{for}~
    |k-k'|>(l \log A)^{C_{2}}\Big\},
   \end{aligned}
   \end{equation*}
   where $G_{l-1;A^{l}}$ and $\textbf{G}_{l-1; A^{l}}$ are defined
   at the beginning of this section.

%
    \item [$(l.4.3)$] There is the measure estimate
    \begin{equation*}
  \emph{mes}~\left(\Pi_{l-1}\setminus
  \Pi_{l}\right)< A^{-(\log l)^{C_{4}}}.
\end{equation*}
  \end{enumerate}

\end{enumerate}

Then there is an absolute positive constant $\epsilon_{*}>0$ such that
if $0<\epsilon<\epsilon_{*}$, there is a parameter set
$\Pi_{m+1}$ and a change of variables
$\Phi_{m+1}: \mathscr{D}(s_{m+1},r_{m+1})\times \scr{O}_{m+1}\rightarrow
\mathscr{D}(s_{m},r_{m})$
being real analytic in
space coordinates and also real analytic in
$\xi$ on the complex domain $\scr{O}_{m+1}$.
The transformation is close to the identity in the sense that
\begin{equation*}
  {}_{r_{m}} \pmb{|} \Phi_{m+1}-id~ \pmb{|}_{\mathscr{D}(s_{m+1},r_{m+1})\times \scr{O}_{m+1}}
< \epsilon_{m}^{1/3},\quad
  {}_{r_{m}} \pmb{|} \Phi_{m+1}-id~ \pmb{|}_{\mathscr{D}(s_{m+1},r_{m+1})\times
  \scr{O}_{m+1}}^{\mathcal{L}}< \epsilon_{m}^{1/4},
\end{equation*}
where $\scr{O}_{m+1}=\scr{O}(\Pi_{m+1}, A^{-(m+1)^{C_{3}}})$.
Furthermore, the new Hamiltonian $H_{m+1}= H_{m}\circ \Phi_{m+1}$
has the form of \eqref{decom. H 1},
and
 the  properties
$(l.1)-(l.4)$ hold with $l$ being replaced by $m+1$.

\end{lemma}

\begin{remark}\label{rmk iter}
  Let us briefly explain the property $(l.4.2)$.
  The parameter set $\Pi_{l}$ is contained in the intersection of three sets.
  The first one refers to the Diophantine conditions.

  For the second set, we look at the definition of $T_{l-1; A^{l}}$,
  \begin{equation*}
    T_{l-1; A^{l}}= T_{l-1}|_{[-A^{l},A^{l}]^{d}}=(D_{l-1}+S_{l-1})|_{[-A^{l},A^{l}]^{d}},
  \end{equation*}
  which originates from solving the homological equation of the following form
  \begin{equation*}
    \partial_{\omega} F^{z}+ \sqrt{-1}~ (\Omega+
  B(x)+ R^{z\bar{z}}(x)) F^{z}= \scr{E}
  \end{equation*}
  by the Fourier expansion and the truncation of the Fourier modes.
  With sufficiently small perturbation,
  for those initial KAM steps ($l$ is close to $l_{*}$), the  matrices $T_{l-1; A^{l}}$
  is diagonally dominated if
  \begin{equation}\label{G dom}
    |D_{l-1}(j,k)|^{-1}=|\Omega_{j}+\langle k,\omega_{l-1}\rangle|^{-1}\lesssim (1/\epsilon)^{1-},\quad \emph{for all}~1\leq j\leq n, |k|\leq A^{l}.
  \end{equation}
  The condition \eqref{G dom} corresponds to the first Melnikov condition
  \begin{equation*}
    |\Omega_{j}+\langle k,\omega_{l-1}\rangle|\neq 0.
  \end{equation*}

  For the third set in $(l.4.2)$, the construction of $\textbf{G}_{l-1;A^{l}}$
  originates from solving
  \begin{equation*}
    \partial_{\omega} F^{zz}+ \sqrt{-1}~ (\Omega+  B(x)
    + R^{z\bar{z}}(x)) F^{zz}+F^{zz}
    (\Omega+  B(x)+ R^{z\bar{z}}(x))=\mathscr{S}.
  \end{equation*}
  Similarly, for those initial KAM iterations, $\textbf{G}_{l-1;A^{l}}$ is  also derived
  from the dominance of the diagonal matrix $\textbf{D}_{l-1}$, which requires
  \begin{equation}\label{G dom 1}
    |\textbf{D}_{l-1}(\textbf{j},k)|^{-1}=|\Omega_{j_{1}}+\Omega_{j_{2}}+
    \langle k,\omega_{l-1}\rangle|^{-1}\lesssim (1/\epsilon)^{1-},\quad
    \emph{for all}~ 1\leq j_{1},j_{2}\leq n, |k|\leq A^{l}.
  \end{equation}
  The condition \eqref{G dom 1} corresponds to the second Melnikov condition
  \begin{equation}\label{2nd mel}
    |\Omega_{j_{1}}\pm\Omega_{j_{2}}+\langle k,\omega_{l-1}\rangle|\neq 0.
  \end{equation}
  However, only the plus sign in \eqref{2nd mel} occurs in our case,
  which can be essentially regarded as the first Melnikov condition since
  $\Omega_{j_{1}}+\Omega_{j_{2}}$ never vanishes.

\end{remark}

\subsection{Proof of  the Iterative Lemma \ref{iter lemma}}\label{proof iter l}

In what follows, we drop the subscript $m$ for simplicity and let
$$
\varepsilon=\epsilon_{m}=A^{-(\frac{4}{3})^{m}},\quad N= A^{m+1}.
$$
The intermediate points $s_{m}^{{(j)}}$ between $s_{m}$ and
$s_{m+1}$ are also written by $s^{{(j)}}$ for $0\leq j\leq 100$.

\subsubsection{Derivation of homological equations.}
Recall that $P^{\textrm{low}}$ is  a polynomial in $y, z, \bar{z}$ of low order and we write
$P^{\textrm{low}}=P^{\textrm{low}}_{\clubsuit}+\langle R^{z\bar{z}}(x)z,\bar{z}\rangle$ with
\begin{equation*}
  P^{\textrm{low}}_{\clubsuit}= R^{x}(x)+ \langle R^{y}(x), y\rangle+ \langle R^{z}(x),z\rangle+
  \langle R^{\bar{z}}(x),\bar{z}\rangle+ \langle R^{zz}(x)z,z\rangle
  + \langle R^{\bar{z}\bar{z}}(x)\bar{z},\bar{z}\rangle.
\end{equation*}
We are looking for a symplectic transformation $\Phi=X_{F}^{t}|_{t=1}$ to
eliminate  $P^{\textrm{low}}_{\clubsuit}$ in the Hamiltonian $H$. As a result,  we take
$F$  in the form of
\begin{equation*}
 F(x,y,z,\bar{z})=  F^{x}(x)+ \langle F^{y}(x), y\rangle+ \langle F^{z}(x), z\rangle+ \langle F^{\bar{z}}(x), \bar{z}\rangle+\langle F^{zz}(x)z, z\rangle+
  \langle F^{\bar{z}\bar{z}}(x)\bar{z}, \bar{z}\rangle.
\end{equation*}
Putting the unsolved term $R^{z\bar{z}}$ into the normal form ${E}$, we get a corrected
normal form
\begin{equation*}
  \underline{E}= {E}+ \langle R^{z\bar{z}}z,\bar{z}\rangle.
\end{equation*}
Then we have
\begin{equation*}
\begin{aligned}
  H\circ \Phi= \underline{E} &+\langle \{\underline{E},F\}^{z\bar{z}} z,\bar{z}\rangle+
  \langle \{P^{\textrm{high}},F\}^{z\bar{z}} z,\bar{z}\rangle
  +\grave{P}+P^{\textrm{high}}
    + \{P^{\textrm{high}},F\}^{\textrm{high}}\\
&+ P^{\textrm{low}}_{\clubsuit}+ \{\underline{E},F\}_{\clubsuit}+
   \{P^{\textrm{high}},F\}^{\textrm{low}}_{\clubsuit},
\end{aligned}
\end{equation*}
where $ \{\underline{E},F\}_{\clubsuit}= \{\underline{E},F\}-\{\underline{E},F\}^{z\bar{z}}$,
$\{P^{\textrm{high}},F\}^{\textrm{low}}_{\clubsuit}=
   \{P^{\textrm{high}},F\}^{\textrm{low}}-\{P^{\textrm{high}},F\}^{z\bar{z}}$
   and
\begin{align}
\grave{P}= \int_{0}^{1}  \{(1-t)\{\underline{E},F\}+P^{\textrm{low}}_{\clubsuit},F\}\circ X_{F}^{t} dt
  +\int_{0}^{1} (1-t) \{\{P^{\textrm{high}},F\},F\}\circ X_{F}^{t} dt,\label{grave P}
    \end{align}

We aim at  solving
\begin{equation}\label{coh-0}
   \{F,\underline{E}\}_{\clubsuit}+ \{F,P^{\textrm{high}}\}^{\textrm{low}}_{\clubsuit}=
   P_{\clubsuit}^{\textrm{low}}.
\end{equation}
As usual, we shall employ the truncation technique.
Recalling the truncation operator $\Gamma_{N}$ defined in \eqref{trunc},
 we solve \eqref{coh-0} up to a admissible error and get the following
homological equations:
\begin{align}
&\partial_{\omega} F^{x}= \Gamma_{N} R^{x}, \label{homo 1}\\
&\partial_{\omega} F^{z}+ \sqrt{-1}~ \Gamma_{N}[(\Omega+
  B(x)+ R^{z\bar{z}}(x)) F^{z}] = \Gamma_{N}\scr{E}, \label{homo 2}\\
&\partial_{\omega} F^{\bar{z}}-\sqrt{-1}~ \Gamma_{N}[ (\Omega+ B(x)
  + R^{z\bar{z}}(x)) F^{\bar{z}}]= \Gamma_{N} \scr{E}',\label{homo 2'}\\
& \partial_{\omega} F^{y}
  = \Gamma_{N}\mathscr{R}-\widehat{\mathscr{R}}(0),   \label{homo 3}\\
& \partial_{\omega} F^{zz}+ \sqrt{-1}~ \Gamma_{N}[(\Omega+  B
    + R^{z\bar{z}}) F^{zz}+F^{zz}
    (\Omega+  B+ R^{z\bar{z}})]=\Gamma_{N}\mathscr{S},\label{homo 4}\\
& \partial_{\omega} F^{\bar{z}\bar{z}}+ \sqrt{-1}~ \Gamma_{N}
   [(\Omega+  B+   R^{z\bar{z}})
    F^{\bar{z}\bar{z}}+ F^{\bar{z}\bar{z}}(\Omega+   B
    +  R^{z\bar{z}})]=
  \Gamma_{N}\mathscr{S}', \label{homo 4'}
\end{align}
where
\begin{align}
& \scr{E}=R^{z}-
   P^{yz}~ \partial_{x} F^{x},\quad \scr{E}'
   = R^{\bar{z}}-P^{y\bar{z}}~
  \partial_{x} F^{x}, \label{E}\\
& \mathscr{R}= R^{y}+\sqrt{-1} (P^{yz} F^{\bar{z}}-P^{y\bar{z}} F^{z})- P^{yy}
  \partial_{x} F^{x}, \label{R}\\
& \mathscr{S}=R^{zz}+[\sqrt{-1}(P^{zzz} F^{\bar{z}}-P^{zz\bar{z}} F^{z})-P^{yz}
  \partial_{x} F^{z}-P^{yzz}
   \partial_{x} F^{x}],\label{S}\\
& \mathscr{S}'= R^{\bar{z}\bar{z}}+
  [\sqrt{-1}(P^{z\bar{z}\bar{z}} F^{\bar{z}}-P^{\bar{z}\bar{z}\bar{z}} F^{z})-P^{y\bar{z}}
\partial_{x} F^{\bar{z}}-P^{y\bar{z}\bar{z}}
   \partial_{x} F^{x}].\label{S'}
\end{align}

Without loss of generality, we assume
\begin{equation*}
  \widehat{R^{x}}(0)=\int_{\TT^{d}} R^{x}(x) dx=0,
\end{equation*}
since the dynamics of the Hamiltonian vector field are unaffected.
The homological equations to be solved are divided into four classes.
The first one is \eqref{homo 1},
which is well known in the KAM theory.
Upon solving \eqref{homo 1}, we turn to the second homological equations
\eqref{homo 2}-\eqref{homo 2'}
Observe that the two equations are complex conjugated if some symmetry is preserved
(to be specified later).
Moreover, \eqref{homo 2}-\eqref{homo 2'} contain the variable  coefficients $B(x)$ and
$R^{z\bar{z}}(x)$. We will
employ the techniques developed in the Anderson localization theory  to solve
 \eqref{homo 2}-\eqref{homo 2'}.
Then we come up with the third  homological equation
\eqref{homo 3}
in which $\scr{R}$
is known from \eqref{homo 1}-\eqref{homo 2'}.
The unsolved constant $\widehat{\mathscr{R}}(0)$ corresponds to the
shift of the tangent frequency during the iterations.
The last homological equations are
\eqref{homo 4}-\eqref{homo 4'},
which are  essentially the same to \eqref{homo 2}-\eqref{homo 2'}.

Once \eqref{homo 1}-\eqref{homo 4'} are solved, we get
\begin{equation*}
  H_{+}= H\circ \Phi= {E}_{+}+ P_{+},
\end{equation*}
where
\begin{equation}\label{E_+}
\begin{aligned}
  E_{+}
  =& \langle\omega+ \widehat{\mathscr{R}}(0),y\rangle+\langle\Omega z,\bar{z}\rangle+\langle
  (B+R^{z\bar{z}}+\{\underline{E},F\}^{z\bar{z}}
  +\{R^{\textrm{high}},F\}^{z\bar{z}})z,\bar{z}\rangle\\
\equiv& \langle\omega_{+}, y\rangle+ \langle\Omega z,\bar{z}\rangle+\langle
B_{+} z , \bar{z}\rangle
  \end{aligned}
\end{equation}
and
\begin{align}
 P_{+}=&P^{\textrm{high}}+\{P^{\textrm{high}},F\}^{\textrm{high}}
 +\grave{P}\label{new 0}\\
 + & (1-\Gamma_{N} ) R^{x}\label{new 1}\\
 +& \langle (1-\Gamma_{N})\scr{E}, z\rangle
   +\langle (1-\Gamma_{N}) \scr{E}',\bar{z}\rangle
   \label{new 2}\\
   +& \sqrt{-1} \langle(1-\Gamma_{N})[(B+R^{z\bar{z}})F^{z}],z\rangle
   - \sqrt{-1} \langle(1-\Gamma_{N})[(B+R^{z\bar{z}})F^{\bar{z}}],\bar{z}\rangle
   \label{new 3}\\
   +& \langle(1-\Gamma_{N})\scr{R},y\rangle\label{new 4}\\
   +& \langle(1-\Gamma_{N})\mathscr{S} z, z\rangle+ \langle(1-\Gamma_{N})
   \mathscr{S}' \bar{z},\bar{z}\rangle\label{new 5}\\
   -& \sqrt{-1}\langle ((1-\Gamma_{N}) [(B^{z\bar{z}}+R^{z\bar{z}})
   F^{zz}+ F^{zz}(B^{z\bar{z}}+R^{z\bar{z}})]) z, z\rangle\label{new 6}\\
   -& \sqrt{-1}\langle(1-\Gamma_{N})[(B^{z\bar{z}}+R^{z\bar{z}})
   F^{\bar{z}\bar{z}}+F^{\bar{z}\bar{z}}(B^{z\bar{z}}+R^{z\bar{z}})]
   \bar{z}, \bar{z}\rangle\label{new 7}.
    \end{align}

\subsubsection{Reality conditions.}
In this part, we establish the reality property of the transformation function $F$ and
the new perturbation $P_{+}$.

\begin{lemma}\label{sym 1}
  Assume $P$ satisfies the following reality condition
  \begin{equation}\label{P real}
  \begin{aligned}
   \overline{ P(x,y,z,\bar{z};\xi)}= P(x,y,z,\bar{z};\xi),\quad
   \forall~ (x,y,z,\bar{z},\xi)\in\mathscr{D}_{\RR}(s,r)\times \Pi.
  \end{aligned}
\end{equation}
Let $R_{ij}^{zz}, R_{ij}^{z\bar{z}}, R_{ij}^{\bar{z}\bar{z}}$ be the elements
of the matrices $R^{zz}, R^{z\bar{z}}, R^{\bar{z}\bar{z}}$ respectively. Then
we have that, for real $x$,
\begin{align*}
& \overline{R^{x}(x)}= R^{x}(x),\quad \overline{R^{y}(x)}= R^{y}(x),\quad
  \overline{R^{z}(x)}= R^{\bar{z}}(x),\\
& R_{ij}^{zz}(x)=R_{ji}^{zz}(x),\quad \overline{R_{ij}^{zz}(x)}
  =R_{ij}^{\bar{z}\bar{z}}(x),\quad
   R_{ij}^{z\bar{z}}(x)=R_{ji}^{z\bar{z}}(x),\quad \overline{R_{ij}^{z\bar{z}}(x)}
  =R^{{z}\bar{z}}_{ji}(x),\\
&\overline{P^{yz}(x)}= P^{y\bar{z}}(x),\quad \overline{P^{yzz}(x)}=P^{y\bar{z}\bar{z}}(x),
  \quad\overline{P^{zzz}(x)}=P^{\bar{z}\bar{z}\bar{z}}(x),\quad
  \overline{P^{zz\bar{z}}(x)}=P^{\bar{z}\bar{z}z}(x)
\end{align*}
\end{lemma}
\noindent\textbf{Proof.} Taking $n=1$ for example, we write
\begin{equation*}
  P(x,y,z,\bar{z};\xi)=\sum_{s,t\geq 0} p_{s,t} z^{s}\bar{z}^{t},
\end{equation*}
where $p_{s,t}= p(x,y;\xi)$. It follows from the reality condition of $P$ that
\begin{equation*}
  \overline{p_{s,t}(x,y;\xi)}=p_{t,s}(x,y;\xi),\quad \forall x\in\TT^{d}, y\in\RR^{d}, \xi\in\Pi.
\end{equation*}
Then
\begin{equation*}
  \overline{R^{x}(x)}=\overline{P(x,y,0,0)}=\overline{p_{0,0}(x,0;\xi)}=p_{0,0}(x,0;\xi)=R^{x}(x).
\end{equation*}
For $\overline{R^{z\bar{z}}(x)}$, we see that
\begin{equation*}
  R^{z\bar{z}}=\partial_{z}\partial_{\bar{z}} P|_{y=z=\bar{z}=0}=
  \left(\sum_{s,t} st p_{st} z^{s-1} \bar{z}^{t-1}\right)|_{y=z=\bar{z}=0}=p_{11}(x,0;\xi)
\end{equation*}
and hence
\begin{equation*}
  \overline{R^{z\bar{z}}(x)}=\overline{p_{11}(x,0;\xi)}=p_{11}(x,0;\xi)=R^{z\bar{z}}(x).
\end{equation*}
The remaining relationships can be verified similarly and are omitted. \qed
\bigskip

By  $(m.2)$ in the Iterative Lemma and the above lemma \ref{sym 1}, we know that ${E}$ is real for $x\in\TT^{d}$,
$y\in\RR^{d}$ and $z, \bar{z}\in\CC^{n}$.

\begin{lemma}\label{sym 2}
Suppose $P$ satisfies the reality condition \eqref{P real}.
  If $F$ in some sub-domain $\scr{D}(s',r')\times \Pi'$ of $\scr{D}(s,r)\times\Pi$
  is the unique solution of the homological equations \eqref{homo 1}-\eqref{homo 4'},
  then
  $F$ satisfies the reality condition:
  \begin{equation}\label{F real}
  \begin{aligned}
   \overline{ F(x,y,z,\bar{z};\xi)}= F(x,y,z,\bar{z};\xi),\quad
   \forall~ (x,y,z,\bar{z},\xi)\in\mathscr{D}_{\RR}(s',r')\times \Pi'.
  \end{aligned}
\end{equation}
\end{lemma}

\noindent\textbf{Proof.} Taking complex conjugation on both sides of \eqref{homo 1}, we see
from $\overline{R^{x}}=R^{x}$ that $\overline{F^{x}}$ is also a solution to
\eqref{homo 1}, which, by the uniqueness of solution, implies that $\overline{F^{x}}=F^{x}$.
For \eqref{homo 2} and \eqref{homo 2'}, it then follows that $\scr{E}'=\overline{\scr{E}}$.
As a result, if $(F^{z},F^{\bar{z}})$ solves \eqref{homo 2}-\eqref{homo 2'}, so does
$(\overline{F^{\bar{z}}},\overline{F^{z}})$. Also by using the uniqueness assumption, we have
$\overline{F^{{z}}}=F^{\bar{z}}$. Similarly, we can show
$\overline{F^{y}}=F^{y}, \overline{F^{zz}}=F^{\bar{z}\bar{z}}$.
Consequently, the reality condition \eqref{F real} of $F$ holds true.
\qed

\begin{proposition}
  Under the assumption of Lemma \ref{sym 2}, we have that the matrix $B_{+}$ is self-adjoint,
  i.e.,
  \begin{equation*}
       \overline{B_{+}^{T}(x;\xi)}= B_{+}(x;\xi),\quad \forall~x\in\TT^{d}, \xi\in\Pi'.
    \end{equation*}
    Moreover, the new perturbation $P_{+}$ satisfies the reality condition
    \eqref{P real} on $\scr{D}_{\RR}(s',r')\times\Pi'$.
\end{proposition}

\noindent\textbf{Proof.} This is an immediate result of the reality property
of $P$ and $F$.\qed

\subsubsection{Solutions of the homological equations.}
In this part, we establish several propositions to solve the homological equations
\eqref{homo 1}-\eqref{homo 4'} in order.

 \emph{Firstly}, we solve the homological equation
\eqref{homo 1}, which is very standard in the classical KAM theory.

\begin{proposition}\label{solu 1}
  \textbf{\emph{(Solution of \eqref{homo 1})}} Under the assumptions of Lemma \ref{iter lemma}, there is
  a parameter set $V^{(1)}$ with $\emph{mes}~(V^{(1)})< A^{-m\tau}$ such that for all $\xi\in\Pi\setminus V^{(1)}$, the Diophantine
  condition holds
 \begin{equation}\label{dioph 1}
|\langle k,\omega(\xi)\rangle|>\frac{\sqrt{\epsilon} (1+2^{-m})}{ |k|^{\tau}},
\quad 0\neq |k|\leq N, k\in\ZZ^{d}.
\end{equation}
Then
equation \eqref{homo 1} has an analytic solution $F^{x}$ defined
on $\TT^{d}_{s^{{(1)}}}\times \scr{O}(\Pi\setminus V^{(1)}, 10 A^{-(m+1)^{C_{3}}})$.
Moreover, we have
\begin{equation*}
  \sup_{(x,\xi)\in \TT^{d}_{s^{{(2)}}}\times \scr{O}(\Pi\setminus V^{(1)}, 8 A^{-(m+1)^{C_{3}}})}\left\{ |F^{x}|, |\partial_{x} F^{x}, |\partial_{\xi} F^{x}|, |\partial_{\xi}\partial_{x} F^{x}|\right\}< \varepsilon^{1-}.
\end{equation*}

\end{proposition}
\noindent\textbf{Proof.}
Recall that  $\widehat{R^{x}}(0)=0$. Moreover,
\begin{equation*}
|R^{x}|\leq 2 \left|\frac{\partial P(x,0,0,0)}{\partial x}\right|\lesssim \pmb{|} X_{P^{\textrm{low}}} \pmb{|} \lesssim \varepsilon
\end{equation*}
and thus $\sup_{x\in\TT^{d}_{s}} |\Gamma_{K} R^{x}| \lesssim \varepsilon.$
 From $(m.4.2)$ in the iterative lemma, there is
  \begin{equation*}
    |\langle k,\omega_{m-1}\rangle|>\frac{\sqrt{\epsilon} (1+2^{-(m-1)})}{ |k|^{\tau}},\quad 0\neq |k|\leq A^{m}.
  \end{equation*}
  Then by $|\omega-\omega_{m-1}|< \epsilon_{m-1}^{1/10}$ in $(m.1)$, we have
  \begin{equation*}
    |\langle k,\omega_{m}\rangle|\geq \frac{\sqrt{\epsilon} (1+2^{-(m-1)})}{ |k|^{\tau}}
    -A^{m} \epsilon_{m-1}^{1/10}>\frac{\sqrt{\epsilon} (1+2^{-m})}{ |k|^{\tau}}
  \end{equation*}
  for all $0\neq |k|\leq A^{m}$.
  It remains to exclude the parameter $\xi$ in
  \begin{equation*}
    {V}^{(1)}=\bigcup_{A^{m}\leq |k|\leq A^{m+1}}
    \left\{\xi\in\Pi_{m}: |\langle k,\omega(\xi)\rangle |< \frac{\sqrt{\epsilon} (1+2^{-m})}
    {|k|^{\tau}}
    \right\},
  \end{equation*}
   where $|\partial_{\xi}\omega|$
  and $|\partial_{\xi}\omega^{-1}|$ are uniformly
  bounded  along the iterations.
  As a result, the total excluded measure
  \begin{equation*}
   \textrm{mes}~{V}^{(1)} \lesssim \sqrt{\epsilon} \sum_{A^{m}\leq |k|\leq A^{m+1}}
    \frac{1}{|k|^{\tau}}< A^{-m\tau}
    \ll A^{-(\log (m+1))^{C_{4}}},
  \end{equation*}
  which is allowed in the measure estimate in the iterative lemma.
  Moreover, the Diophantine condition \eqref{dioph 1} remains valid on an $10 A^{-(m+1)^{C_{3}}}$-neighborhood
  of $\Pi_{m}\setminus {V}^{(1)}$ due to the fact that
  $A^{-(m+1)^{C_{3}}}\ll A^{-100 \tau (m+1)}$.
  Using the Diophantine condition, the existence of the solution $F^{x}$ of
  \eqref{homo 1}, as well as its estimates on derivatives, is well known in KAM theory.
  We omit the proof here.
  \qed

%
%

%
%
%
%
%
%
%
%

\bigskip

\emph{Secondly}, we solve the homological equations \eqref{homo 2}-\eqref{homo 2'},
by using some techniques from the Craig-Wayne-Bourgain method.
From Lemma \ref{sym 2}, it suffices to solve \eqref{homo 2} since
$F^{\bar{z}}=\overline{F^{z}}$.

\begin{proposition}\label{solu 2}
 \textbf{\emph{(Solution of \eqref{homo 2}-\eqref{homo 2'})}}
  Under the assumptions of Lemma \ref{iter lemma}, there exists a parameter set
  $\Pi_{+}^{(1)}\subset\Pi\setminus V^{(1)}$ such that
  for all $\xi\in\Pi_{+}^{(1)}$, equations
  \eqref{homo 2}-\eqref{homo 2'} have a unique  solution $(F^{z}, F^{\bar{z}})$
  with $F^{\bar{z}}=\overline{F^{z}}$, which admits analytic extension to
  $\TT^{d}_{s^{{(3)}}}\times \scr{O}(\Pi_{+}^{(1)}, 8 A^{-(m+1)^{C_{3}}})$
 and satisfies
  \begin{equation*}
    \sup_{\TT^{d}_{s^{{(4)}}}\times \scr{O}(\Pi_{+}^{(1)}, 6 A^{-(m+1)^{C_{3}}})} \left\{
    |F^{z}|, |\partial_{x} F^{z}|, |\partial_{\xi} F^{z} |, |\partial_{\xi}\partial_{x} F^{z}|
    \right\}< \varepsilon^{1-}.
  \end{equation*}
  Moreover,
  $\Pi_{+}^{(1)}$ is the union of  a family of  disjoint intervals $\scr{J}'$ with  $|\scr{J}'|=A^{-(m+1)^{C_{3}}}$. For each
$\scr{J}'$, there is a unique $\scr{J}\in\Lambda_{m}$ such that $\scr{J'}\subset \scr{J}$.
The total removed set satisfies
\begin{equation*}
  \emph{mes}~(\Pi\setminus \Pi_{+}^{(1)})<\frac{1}{3} A^{-(\log (m+1))^{C_{4}}}.
\end{equation*}

\end{proposition}

\noindent\textbf{Proof.}~
Due to Lemma \ref{sym 2}, we only
solve the homological equation \eqref{homo 2} with $\scr{E}$ given by \eqref{E}
\begin{equation*}
 \partial_{\omega} F^{z}+ \sqrt{-1}~ \Gamma_{N}(\Omega+
 \Gamma_{N} B(x)+\Gamma_{N} R^{z\bar{z}}(x)) F^{z} = \Gamma_{N}\mathscr{E},
\end{equation*}
which, by matching the components of the vector-valued functions, turns out to be
\begin{equation*}
-\sqrt{-1}\partial_{\omega} F^{z}_{j}+  \Gamma_{N} \Omega_{j} F^{z}_{j}
+\Gamma_{N}\left(\sum_{1\leq r\leq n} (\Gamma_{N} B_{jr}(x)) F^{z}_{r}(x))
+\sum_{1\leq r\leq n} (\Gamma_{N} R^{z\bar{z}}_{jr}(x)) F^{z}_{r}(x)\right)=
-\sqrt{-1} \Gamma_{N}\mathscr{E}_{j}
\end{equation*}
for  $1\leq j\leq n$.
We are looking for solution $F^{z}(x)$ with compact support in the Fourier modes
\begin{equation*}
F^{z}_{j}(x)=\sum_{|k|\leq N, k\in\ZZ^{d}} \widehat{F^{z}_{j}}(k)
e^{\sqrt{-1} \langle k, x\rangle}.
\end{equation*}
Passing to the Fourier transformation, we then get
\begin{equation}\label{homo 2-1'}
 (\langle k,\omega\rangle+ \Omega_{j}) \widehat{F^{z}_{j}}(k)+\sum_{1\leq r\leq n}
 \sum_{|p|\leq N,  |k-p|\leq N}
 (B_{jr}+R^{z\bar{z}}_{jr})^{\wedge}(k-p) \widehat{F^{z}_{r}}(p)
 =-\sqrt{-1} \widehat{\mathscr{E}_{j}}(k),
\end{equation}
where $1\leq j\leq n$, $|k|\leq N$, $k\in\ZZ^{d}$ and $(\cdot)^{\wedge}(k)$ stands
for the $k$-th Fourier coefficient of the indicated function.

Thinking of  $\widehat{F}$ as a vector defined on $\{1,\cdots,n\}\times \ZZ^{d}$,
we write \eqref{homo 2-1'} in a matrix form. To this end, we
let
\begin{equation}\label{T}
  T=D+S
\end{equation}
where  $D$  is a diagonal matrix
\begin{equation*}
D(j, k)= \Omega_{j}+ \langle k,\omega\rangle
\end{equation*}
and  $S$ is a non-diagonal matrix
\begin{equation*}
S((j,k),(r, p))= (B_{jr}+R^{z\bar{z}}_{jr})^{\wedge}(k-p).
\end{equation*}
We denote by $T_{N}$ and $\widehat{\scr{E}}_{N}$ the restriction of the matrix  $T$
and the restriction of the vector $\widehat{\scr{E}}$ on
$\{1,\cdots,n\}\times\{k\in\ZZ^{d}: |k|\leq N\}$
respectively.
With the notations introduced above, equation
\eqref{homo 2-1'} is equivalent to
\begin{equation}\label{homo 2-1''}
  T_{N} \widehat{F}^{z}=-\sqrt{-1} \widehat{\mathscr{E}}_{N}.
\end{equation}
Then our goal is to establish the existence and the decay property of
 the inverse matrix of $T_{N}$ (see \eqref{G N decay}),
which is also called the Green's function estimate.
Indeed, we have the following results, whose proof is delayed
to the next section.

\begin{lemma}\label{Green}
  Under the assumptions of \emph{Lemma \ref{iter lemma}}, there exists a parameter set
  $\Pi_{+}^{(1)}\subset\Pi\setminus V^{(1)}$ such that for any $\xi\in\scr{O}(\Pi_{+}^{(1)}, 10 A^{-(m+1)^{C_{3}}})$, the inverse matrix $G_{N}=T_{N}^{-1}$ exists and satisfies
  \begin{equation}\label{G N decay}
  \begin{aligned}
&\|G_{N}\|< A^{(\log N)^{C_{1}}},\\
& |G_{N}(x,y)|< e^{-(s- (\log N)^{-8})~ |x-y| }\quad \emph{for}~|x-y|> (\log N)^{C_{2}}.
  \end{aligned}
\end{equation}
Moreover, $\Pi_{+}^{(1)}$ is the union of  a family of  disjoint intervals $\scr{J}'$ with  $|\scr{J}'|=A^{-(m+1)^{C_{3}}}$. For each
$\scr{J}'$, there is a unique $\scr{J}\in\Lambda_{m}$ such that $\scr{J'}\subset \scr{J}$.
The total removed set satisfies
\begin{equation*}
  \emph{mes}~(\Pi\setminus \Pi_{+}^{(1)})<\frac{1}{3} A^{-(\log (m+1))^{C_{4}}}.
\end{equation*}
\end{lemma}

\bigskip

Now we apply Lemma \ref{Green} to solve \eqref{homo 2-1''} and then to prove Proposition
\ref{solu 2}.
Recall the definition of $\scr{E}$ in \eqref{E}. Observe that
$
  R^{z}=\partial_{z} P(x,0,0,0),
$
and it follows from Cauchy's estimate that
\begin{equation*}
  \sup_{x\in\TT^{d}_{s}} |R^{z}(x)|< {~}_{r} \pmb|X_{P^{\textrm{low}}}\pmb|\cdot r\lesssim\varepsilon.
\end{equation*}
Moreover, $P^{yz}=\partial_{y}\partial_{z} P(x,0,0,0)$ and then
\begin{equation}\label{R yz}
  \sup_{x\in\TT^{d}_{s}}|P^{yz}|\leq \frac{1}{r}
  \sup |\partial_{y} P(x,0,z,0)|< {~}_{r} \pmb|X_{P^{\textrm{high}}}\pmb|\cdot r\lesssim \epsilon.
\end{equation}
Combining with Proposition \ref{solu 1}, we obtain
\begin{equation}\label{E decay}
  |\widehat{\scr{E}}(k)|\lesssim \varepsilon^{1-} e^{-s^{2} |k|}
\end{equation}
for any $k\in\ZZ^{d}$.

Back to \eqref{homo 2-1''}, we have
\begin{equation*}
\begin{aligned}
&\widehat{F}^{z}= - \sqrt{-1} G_{N} \widehat{\mathscr{E}}_{N},\\
&\widehat{F}^{z}(k)=-\sqrt{-1} \sum_{|p|\leq N, p\in\ZZ^{d}} G_{N}(k,p) \widehat{\scr{E}}_{N}(p).
  \end{aligned}
\end{equation*}
It then follows that
\begin{equation}\label{F z decay}
\begin{aligned}
  |\widehat{F}^{z}(k)|\leq & \varepsilon^{1-} \sum_{|k-p|\leq (\log N)^{C_{2}}} A^{(\log N)^{C_{1}}}
  e^{-s^{{(2)}} |p|}
  + \varepsilon^{1-} \sum_{|k-p|>(\log N)^{C_{2}}}
  e^{-(s-(\log N)^{-8}) |k-p|} e^{-s^{{(2)}} |p|}\\
<& \varepsilon^{1-} (\log N)^{C} A^{(\log N)^{C_{1}}} e^{s^{{(2)}} (\log N)^{C_{2}}} e^{-s^{{(2)}}
  |k|}
  +C \varepsilon^{1-} e^{-s^{{(2)}} |k|}\\
<& \varepsilon^{1-} e^{-s^{{(2)}} |k|}.
\end{aligned}
\end{equation}

Passing to the function $F^{z}$, we then have, for any
$ x\in\TT^{d}_{s^{{(3)}}}$,
\begin{equation*}
   |F^{z}(x)|\leq \sum_{k\in\ZZ^{d}} |\widehat{F}^{z}(k)|\cdot
   |e^{\sqrt{-1} \langle k, x\rangle}|< \varepsilon^{1-}.
\end{equation*}
The remaining estimates for $\partial_{x} F^{z}$, $\partial_{\xi} F^{z}$ and
$\partial_{\xi}\partial_{\xi} F^{z}$ follow from the Cauchy's estimate by
using the fact that $(m+1)^{2}\ll A^{(m+1)^{C_{3}}} \ll (1/\varepsilon)^{0+}$.
This completes the proof of Proposition \ref{solu 2}.
\qed
\bigskip

\emph{Thirdly}, we solve the homological equation \eqref{homo 3}, which is similar to
that of \eqref{homo 1}. To begin with, we need to estimate $\scr{R}$ defined
by \eqref{R}.
Observe that
$ R^{y}=\partial_{y}P(x,0,0,0)$
and it follows from Cauchy estimate that
\begin{equation*}
  \sup_{x\in\TT_{s}^{d} } |R^{y}(x)|< {~}_{r}\pmb|X_{P^{\textrm{low}}}\pmb|\cdot r<\varepsilon
\end{equation*}
Moreover, $P^{yy}=\partial_{yy}^{2} P(x,0,0,0)$ and
\begin{equation*}
  \sup_{x\in\TT^{d}_{s}}|P^{yy}(x)|\leq \frac{1}{r^{2}}
  \sup |\partial_{y} P(x,0,0,0)|\leq {~}_{r}\pmb{|} X_{P^{\textrm{high}}}\pmb{|}\leq \epsilon.
\end{equation*}
Then we see from \eqref{R yz}, Proposition \ref{solu 1} and Proposition \ref{solu 2} that
\begin{equation}\label{scr R}
  \sup_{x\in\TT^{d}_{s^{{(4)}}}} |\scr{R}|< \varepsilon^{1-}.
\end{equation}
Recalling the frequency shift induced by the unsolved term $\widehat{\scr{R}}(0)$, we have
\begin{equation}\label{Del omega}
  |\omega_{+}-\omega|= |\widehat{\scr{R}}(0)|< \varepsilon^{1-}.
\end{equation}

\begin{proposition}\label{solu 3}
\textbf{\emph{(Solution of \eqref{homo 3})}}
Under the assumptions of Lemma \ref{iter lemma}, equation \eqref{homo 3}
has an analytic solution $F^{y}$ defined on
$\TT^{d}_{s^{{(5)}}}\times \scr{O}(\Pi_{+}^{(1)}, 6 A^{-(m+1)^{C_{3}}})$
satisfying
\begin{equation*}
  \sup_{(x,\xi)\in \TT^{d}_{s^{{(6)}}}\times \scr{O}(\Pi_{+}^{(1)}, 4 A^{-(m+1)^{C_{3}}})}\left\{ |F^{y}|, |\partial_{x} F^{y}, |\partial_{\xi} F^{y}|, |\partial_{\xi}\partial_{x} F^{y}|\right\}< \varepsilon^{1-}.
\end{equation*}
\end{proposition}

\emph{Finally}, we solve the homological equations \eqref{homo 4}-\eqref{homo 4'},
which are essentially the same to \eqref{homo 2}-\eqref{homo 2'}.
To begin with, we need to estimate $\scr{S}$ and $\scr{S}'$ defined in \eqref{S}
and \eqref{S'} respectively.
From the estimates of $F^{x}, F^{z}, F^{\bar{z}}$, it is easy to see
\begin{equation}\label{scr S decay}
  \sup_{x\in\TT^{d}_{s^{{(6)}}}} \left\{|\scr{S}|,|\scr{S}'|\right\}
  \lesssim \varepsilon^{1-}.
\end{equation}

\begin{proposition}\label{solu 4}
\textbf{\emph{(Solution of \eqref{homo 4}-\eqref{homo 4'})}}
  Under the assumptions of Lemma \ref{iter lemma}, there exists a parameter set
 $\Pi_{+}^{(2)}\subset\RR^{d}$ such that
  for all $\xi\in \Pi_{+}=\Pi_{+}^{(1)}\cap \Pi_{+}^{(2)}$, equations
  \eqref{homo 4}-\eqref{homo 4'} have a unique  solution $(F^{zz}, F^{\bar{z}\bar{z}})$
  with $F^{\bar{z}\bar{z}}=\overline{F^{zz}}$, which admits analytic extension to
  $\TT^{d}_{s^{{(7)}}}\times \scr{O}(\Pi_{+}, 6 A^{-(m+1)^{C_{3}}})$
  and satisfies
  \begin{equation*}
    \sup_{\TT^{d}_{s^{{(8)}}}\times \scr{O}(\Pi_{+}, 4 A^{-(m+1)^{C_{3}}})} \left\{
    |F^{zz}|, |\partial_{x} F^{zz}|, |\partial_{\xi} F^{zz} |, |\partial_{\xi}\partial_{x} F^{zz}|
    \right\}< \varepsilon^{1-}.
  \end{equation*}
  Moreover, $\Pi_{+}^{(2)}$ is the union of  a family of  disjoint intervals $\scr{J}'$ with  $|\scr{J}'|=A^{-(m+1)^{C_{3}}}$. For each
$\scr{J}'$, there is a unique $\scr{J}\in\Lambda_{m}$ such that $\scr{J'}\subset \scr{J}$.
The total removed set satisfies
\begin{equation*}
  \emph{mes}~(\Pi\setminus \Pi_{+}^{(2)})<\frac{1}{3} A^{-(\log (m+1))^{C_{4}}}.
\end{equation*}
\end{proposition}

\noindent\textbf{Proof.}
Note that the unknown function $F^{zz}$ is of matrix value, i.e.,
$F^{zz}(x)=(F^{zz}_{ij}(x))_{1\leq i,j\leq n}$.
Writing equation \eqref{homo 4} into components yields
\begin{equation}\label{homo 4 1-1}
  \begin{aligned}
    \partial_{\omega}F^{zz}_{ij}+\sqrt{-1}\left[\Gamma_{N} (\Omega_{i}+\Omega_{j}) F^{zz}_{ij}
    +\Gamma_{N}\sum_{p=1}^{n} \Gamma_{N}( B^{z\bar{z}}_{ip}+R^{z\bar{z}}_{ip}) F^{zz}_{pj}
    + F^{zz}_{i p} \Gamma_{N}( B^{z\bar{z}}_{pj}+R^{z\bar{z}}_{pj})
    \right]= \Gamma_{N} \scr{S}_{ij}.
  \end{aligned}
\end{equation}
Let
\begin{equation*}
  \textbf{j}=(i,j),\quad \Omega_{\textbf{j}}= \Omega_{i}+\Omega_{j},\quad 1\leq i, j\leq n,
\end{equation*}
and hence  $\textbf{j}$ is an index taking $n^{2}$ many values.
Then \eqref{homo 4 1-1} is equivalent to
\begin{equation}\label{homo 4 1-2}
  \partial_{\omega} F^{zz}_{\textbf{j}}+\sqrt{-1} \Gamma_{N} \Omega_{\textbf{j}}
  F^{zz}_{\textbf{j}}+\sqrt{-1}\sum_{\textbf{j}'}\Gamma_{N}
  (\textbf{B}_{\textbf{j}\textbf{j}'}^{z\bar{z}}+\textbf{R}_{\textbf{j}\textbf{j}'}^{z\bar{z}})
   F^{zz}_{\textbf{j}'}=\Gamma_{N}\scr{S}_{\textbf{j}},
\end{equation}
where $\textbf{j}'=(i',j')$,
\begin{equation*}
   \textbf{B}^{z\bar{z}}_{\textbf{j}\textbf{j}'}=
  \left\{
  \begin{array}{lcl}
& B^{z\bar{z}}_{ii'},\quad &\textrm{for}~j'=j, i'\neq i,\\
& B_{j'j}^{z\bar{z}},\quad &\textrm{for}~j'\neq j, i'=i,\\
& B_{ii}^{z\bar{z}}+B^{z\bar{z}}_{jj},\quad &\textrm{for}~j'=j, i'=i,\\
& 0,\quad &\textrm{otherwise},
  \end{array}
  \right.~
  \textrm{and}\quad \textbf{R}^{z\bar{z}}_{\textbf{j}\textbf{j}'}=
  \left\{
  \begin{array}{lcl}
& R^{z\bar{z}}_{ii'},\quad &\textrm{for}~j'=j, i'\neq i,\\
& R_{j'j}^{z\bar{z}},\quad &\textrm{for}~j'\neq j, i'=i,\\
& R_{ii}^{z\bar{z}}+R^{z\bar{z}}_{jj},\quad &\textrm{for}~j'=j, i'=i,\\
& 0,\quad &\textrm{otherwise}.
  \end{array}
  \right.\\
\end{equation*}
We look for solution $F^{zz}_{\textbf{j}}(x)$ of \eqref{homo 4 1-2} in the following form
\begin{equation*}
  F^{zz}_{\textbf{j}}(x)= \sum_{|k|\leq N, k\in\ZZ^{d}} \widehat{F^{zz}_{\textbf{j}}}(k)
  e^{\sqrt{-1} \langle k, x\rangle}.
\end{equation*}
Expanding \eqref{homo 4 1-2} into Fourier series and matching the coefficients yield
\begin{equation}\label{homo 4 1-3}
  (\langle k,\omega\rangle+\Omega_{\textbf{j}}) \widehat{F^{zz}_{\textbf{j}}}(k)
  +\sum_{\textbf{j}'} \sum_{\substack{k'\in\ZZ^{d}, |k'|\leq N,\\ |k-k'|\leq N }}
  (\textbf{B}_{\textbf{j}\textbf{j}'}+\textbf{R}^{z\bar{z}}_{\textbf{j}\textbf{j}'})^{\wedge}
  (k-k')  \widehat{F^{zz}_{\textbf{j}}}(k')= -\sqrt{-1}~ \widehat{\scr{S}_{\textbf{j}}}(k).
\end{equation}
where $k\in\ZZ^{d}$ and $|k|\leq N$.

Writing further \eqref{homo 4 1-3} into a matrix equation like \eqref{homo 2-1''},
we obtain an essentially same matrix $\textbf{T}_{N}$ except the difference between the finite
index $\textbf{j}$ and $j$. Note also the fact that $\Omega_{\textbf{j}}=\Omega_{i}+\Omega_{j}$
with $\textbf{j}=(i,j)$ and $\Omega_{\textbf{j}}$ never vanishes.
Therefore, for $\textbf{G}_{N}$ defined on $\{\textbf{j}=(i,j): 1\leq i,j\leq n\}\times \{k\in\ZZ^{d}: |k|\leq N\}$, we are also able to establish the Green's function estimate like Proposition \ref{solu 2} and obtain the desired parameter set
$\Pi_{+}^{(2)}$.

The remaining estimate of $F^{zz}$ is the same to that of $F^{z}$ and we omit it here.
\qed
\bigskip

Let $\Pi_{m+1}=\Pi_{+}=\Pi_{+}^{(1)}\cap \Pi_{+}^{(2)}$. One easily finds that
$\Pi_{m+1}$ satisfies $((m+1).4)$. See  Remark \ref{Pi +}
for more details.

\subsubsection{The estimate of new error.}
Now we are at the stage of estimating the new terms after the symplectic transformation,
which are given in \eqref{E_+} and \eqref{new 1}-\eqref{new 6}.
The majority of them arise from the remaining terms after the truncation.

For $\omega_{+}$, it follows from \eqref{Del omega} that
\begin{equation*}
  |\omega_{+}-\omega|=|\Delta\omega|<\varepsilon^{1-}.
\end{equation*}
For $B_{+}$, we have
\begin{equation*}
  B_{+}-B= R^{z\bar{z}}-\{\underline{E}, F\}^{z\bar{z}}- \{P^{\textrm{high}}, F\}^{z\bar{z}}
\end{equation*}
in which, for $|\textrm{Im} x|_{\infty}\leq {s}$,
\begin{align*}
& |R^{z\bar{z}}(x)|=|\partial_{z\bar{z}} P(x,0,0,0)|\leq \frac{1}{r}
  \sup_{\scr{D}(s,r)} |\partial_{\bar{z}}P^{\textrm{low}}(x,0,z,0)|\lesssim \varepsilon.
\end{align*}
It then follows from Proposition \ref{solu 3} that
\begin{align*}
 {\sup_{|\textrm{Im} x|_{\infty}<s^{(6)}}}|\{\underline{E}, F\}^{z\bar{z}}|\leq {\sup_{|\textrm{Im} x|_{\infty}<s^{(6)}}}|\partial_{x} B+\partial_{x}R^{z\bar{z}}|\cdot |F^{y}|
\lesssim \varepsilon^{1-} (m+1)^{2}< \varepsilon^{1-}.
\end{align*}
%
Similarly, we obtain $\sup_{x\in\TT^{d}_{s^{(4)}}}|\{P^{\textrm{high}}, F\}^{z\bar{z}}|< \varepsilon^{1-}$
and then we have
\begin{equation*}
  \sup_{|\textrm{Im} x|_{\infty}\leq s^{(6)}} |B_{+}(x)-B(x)| < \varepsilon^{1-},
\end{equation*}
and
\begin{equation*}
  \sup_{|\textrm{Im} x|_{\infty}\leq s^{(6)}} |B_{+}(x)|\leq \sup |B_{l_{*}}|+ \sum_{l=l_{*}}^{m+1} \sup |B_{l}-B_{l-1}|\lesssim \epsilon.
\end{equation*}

Next we estimate \eqref{new 1}-\eqref{new 6}.
Obviously, using \cite[Lemma A.2]{Pos01}, we have
\begin{equation*}
  \sup_{x\in\TT^{d}_{s^{(2)}},\alpha\in\{0,1\} }
  \left
  |\partial_{x}^{\alpha}[(1-\Gamma_{N}) R^{x}]\right|
  \lesssim \frac{N^{d}}{(s-s^{(2)})^{C}} e^{-N (s-s^{(2)})}
  \sup_{\TT^{d}_{s}} |R^{x}|<\frac{1}{100} A^{-(\frac{4}{3})^{m+1}}=\frac{1}{100}\varepsilon^{4/3}
\end{equation*}
provided $A\gg 1$. For $\grave{P}$ defined in \eqref{grave P}, we have
\begin{equation*}
  \sup_{\scr{D}({s^{(10)},r^{(10)}})} \left\{|\grave{P}|, |\nabla \grave{P}|\right\}\lesssim \varepsilon^{2-}
  (m+1)^{C}<\frac{1}{100} \varepsilon^{4/3}.
\end{equation*}
where $\nabla\grave{P}=(\partial_{x}\grave{P},\partial_{y}\grave{P},\partial_{z}\grave{P},
\partial_{\bar{z}}\grave{P})$.
For \eqref{new 2}, we see from \eqref{E decay} that
\begin{equation}\label{hua E}
  \sup_{x\in\TT^{d}_{ s^{(4)}}, \alpha\in\{0,1\}} \left|\partial_{x}^{\alpha} [(1-\Gamma_{N})\scr{E}]\right|\lesssim N^{d} e^{-N (s^{(2)}-s^{(4)})}
  \frac{\varepsilon^{1-}}{(s^{(2)}-s^{(4)})^{C}}<\frac{1}{100} \varepsilon^{4/3}.
\end{equation}
Estimate \eqref{hua E} still holds when
replacing $\scr{E}$ by $\scr{E}'$.
By \eqref{scr R}, \eqref{scr S decay} we get
\begin{equation*}
  \sup_{x\in\TT^{d}_{s^{(8)}}, \alpha\in\{0,1\}}
  \left\{|\partial_{x}^{\alpha}[(1-\Gamma_{N}) \scr{R}]|,
  |\partial_{x}^{\alpha}[(1-\Gamma_{N}) \scr{S}]|,
  |\partial_{x}^{\alpha}[(1-\Gamma_{N}) \scr{S}']|\right\}
< \frac{1}{100} \varepsilon^{4/3},
\end{equation*}
which controls \eqref{new 4} and \eqref{new 5}.

For \eqref{new 3}, we see that
\begin{equation*}
  (1-\Gamma_{N}) ((B+R^{z\bar{z}}) F^{z})=\sum_{|k|\geq N} ((B+R^{z\bar{z}}) F^{z})^{\wedge}(k)
  e^{\sqrt{-1} \langle k, x\rangle}.
\end{equation*}
Since $B_{l}(x)$ and $R_{l}^{z\bar{z} }(x)=\partial_{\bar{z}}\partial_{z} P_{l}^{\textrm{low}}(x,0,0,0)$ are analytic
in $x\in\TT^{d}_{s_{l}}$, we have
\begin{equation}\label{T regularity}
|\widehat{B}_{l}(k)|\lesssim \epsilon e^{-s_{l} |k|}, \quad |\widehat{R^{z\bar{z}}}(k)|<  \frac{1}{r_{l}}
\sup_{\mathscr{D}(s_{l},r_{l})} |\partial_{z} P_{l}| e^{- s_{l} |k|}< \varepsilon_{l}
 e^{-s_{l} |k|}.
\end{equation}
It follows from \eqref{T regularity} and \eqref{F z decay} that
\begin{equation*}
\begin{aligned}
  |((B+R^{z\bar{z}}) F^{z})^{\wedge}(k)|=&|\sum_{|p|\leq N} (B+R^{z\bar{z}})^{\wedge}(k-p)
  \widehat{F^{z}}(p)|\\
  \lesssim & \sum_{|p|\leq N, p\in\ZZ^{d}}  (\epsilon+\varepsilon) e^{-s|k-p|}
  \varepsilon^{1-} e^{-{s^{(4)}} |p|}
< \varepsilon^{1-} e^{-{s^{(4)}} |k|}.
\end{aligned}
\end{equation*}
Then
\begin{equation*}
  \sup_{|\textrm{Im} x|_{\infty}\leq {s^{(6)}}} |(1-\Gamma_{N}) [(B+R^{z\bar{z}}) F^{z}] |<
  \varepsilon^{1-} \sum_{|k|\geq N} e^{-{s^{(4)}} |k|} e^{{s^{(6)}} |k|}
  \lesssim \varepsilon^{1-} N^{d} e^{-({s^{(4)}-s^{(6)}}) N}<\frac{1}{100} \varepsilon^{4/3}.
\end{equation*}
With the margins in our estimate, we further have
\begin{equation*}
   \sup_{|\textrm{Im} x|\leq {s^{(6)}}, \alpha\in\{0,1\}} \left\{|\partial_{x}^{\alpha} [(1-\Gamma_{N}) ((B+R^{z\bar{z}}) F^{z})] |,
   |\partial_{x}^{\alpha}[(1-\Gamma_{N}) ((B+R^{z\bar{z}}) F^{\bar{z}})] |
   \right\}<
   \frac{1}{100} \varepsilon^{4/3}.
\end{equation*}
The estimate of \eqref{new 6}-\eqref{new 7} is the same to that of \eqref{new 3} and
reads
\begin{equation*}
  \sup_{|\textrm{Im} x|\leq {s^{(10)}}, \alpha\in\{0,1\}} \left\{
  \begin{aligned}
    |\partial_{x}^{\alpha}[(1-\Gamma_{N}) ((B^{z\bar{z}}+R^{z\bar{z}})
   F^{zz}+ F^{zz}(B^{z\bar{z}}+R^{z\bar{z}}))]|,\\
   |\partial_{x}^{\alpha}[(1-\Gamma_{N}) ((B^{z\bar{z}}+R^{z\bar{z}})
   F^{\bar{z}\bar{z}}+ F^{\bar{z}\bar{z}}(B^{z\bar{z}}+R^{z\bar{z}}))]|
  \end{aligned}
  \right\}< \frac{1}{100} \varepsilon^{4/3}.
\end{equation*}

Note that
\begin{equation*}
  P_{+}^{\textrm{low}}=\grave{P}^{\textrm{low}}+ \eqref{new 1}+\eqref{new 2}
  +\cdots+\eqref{new 7},
\end{equation*}
and
\begin{equation*}
  P_{+}^{\textrm{high}}= P^{\textrm{high}}+ \grave{P}^{\textrm{high}}
  +\{P^{\textrm{high}}, F\}^{\textrm{high}}.
\end{equation*}
Take
\begin{equation*}
\varepsilon_{+}=\varepsilon^{4/3}=A^{-(\frac{4}{3})^{m+1}},\quad s_{+}=s_{m+1},\quad r_{+}=r_{m+1},\quad
\scr{O}_{+}=\scr{O}(\Pi_{+}, A^{-(m+1)^{C_{3}}}).	
\end{equation*}

We obtain from the above analysis that
\begin{equation*}
  ~_{r_{+}}\pmb{|}X_{P^{\textrm{low}}_{+}}\pmb{|}_{\scr{D}(s_{+},r_{+})\times \scr{O}_{+}}
< \varepsilon_{+}
  \end{equation*}
and
\begin{equation*}
   ~_{r_{+}}\pmb{|}X_{P^{\textrm{high}}_{+}}\pmb{|}_{\scr{D}(s_{+},r_{+})\times \scr{O}_{+}}
   \lesssim \epsilon.
\end{equation*}
The transformation $\Phi=X_{F}^{t}|_{t=1}$ is also close to the identity in the sense that
\begin{equation*}
  ~_{r_{+}}\pmb{|} \Phi-id|_{\scr{D}(s_{+},r_{+})}< \varepsilon^{1/3}
\end{equation*}
since $\sup_{\scr{D}(s_{+},r_{+})}|\nabla F|< \varepsilon^{1-}$.

\subsection{Proof of  the main Theorem \ref{main th}}\label{proof main}
Let the constant $A$  be sufficiently large.
In the proof of the iterative lemma, we take extensively advantage of
the largeness of the iteration step $l$. As a result, it suffices to
start the iteration from  {$l=l_{*},l_{*}=l_{*}(\epsilon)\gg 1$} (independent of the
iterations) instead of $l=0$.

We then need to verify the induction statements at $l={l_{*}}$.
Recall our imposition \eqref{initial func} in the first step
\begin{equation*}
	H_{l_{*}}=H_{0}, ~ P_{l_{*}-1}=P_{l_{*}}=P_{0},~ B_{l_{*}}=B_{l_{*}-1}=0,~\Pi_{l_{*}-1}=\Pi_{0},~s_{l_{*}}=s_{0},~r_{l_{*}}=r_{0},~
\omega_{l_{*}}=\omega_{l_{*}-1}=\omega_{0}.
\end{equation*}
Obviously,  the statements $(l_{*}.1), (l_{*}.2)$
and $ (l_{*}.3)$  hold.
It suffices to find the set $\Pi_{l_{*}}$ such that the
statement
$(l_{*}.4)$ holds, which can essentially be
described by the Diophantine condition and  the first Melnikov condition.
The construction of  the set $\Pi_{l_{*}}$ is given below.

We first pave the set $\Pi_{l_{*}-1}=\Pi_{0}$ into
a $\tilde{\Lambda} $ family of disjoint intervals of diameter
$A^{-l_{*}^{C_{3}}}$,
i.e., $\Pi_{0}=\cup_{\scr{J}\in\tilde{\Lambda}} \scr{J}$
with $|\scr{J}|=A^{-l_{*}^{C_{3}}}$ for each
$\scr{J}\in\tilde\Lambda$. If there exists
some $\xi_{0}\in\scr{J}$ such that
\begin{equation}\label{0 dioph}
|\langle k, \omega_{0}(\xi_{0})\rangle|> \sqrt{\epsilon}
(1+2^{-(l_{*}-1)}) |k|^{-\tau}
\end{equation}
violates for some $|k|\leq A^{l_{*}}$,
then for any $\xi\in\scr{J}$, there is
\begin{equation*}
	|\langle k, \omega_{0}(\xi)\rangle|\leq
	\sqrt{\epsilon}
(1+2^{-(l_{*}-1)}) |k|^{-\tau}+ C A^{l_{*}} A^{-l_{*}^{C_{3}}}
<2 \sqrt{\epsilon}
(1+2^{-(l_{*}-1)}) |k|^{-\tau}
\end{equation*}
Let
\begin{equation*}
\Lambda_{l_{*}}^{(1)}=\{\scr{J}\in\tilde\Lambda:~
\eqref{0 dioph}~\textrm{holds for all}~\xi\in\scr{J}
~\textrm{and all}~0\neq |k|\leq A^{l_{*}}
\},	
\end{equation*}
and one easily sees from the twist condition in
\textbf{Assumption B} that
\begin{equation}\label{0 mes}
	\textrm{mes}~\left(\bigcup_{\scr{J}\in\tilde{\Lambda}\setminus
	\Lambda_{l_{*}}^{(1)}} \scr{J}\right)\lesssim \sqrt{\epsilon}.
\end{equation}

Next we consider the matrix $T_{l_{*-1}}=D_{l_{*}-1}+S_{l_{*}-1}$ with
\begin{equation*}
	D_{l_{*}-1}(j,k)=\Omega_{j}+\langle k,\omega_{0}\rangle
\end{equation*}
and $S_{l_{*}-1}((j,k),(j',k'))=\widehat{P_{0;jj'}^{z\bar{z}}}(k-k')$
since $B_{l_{*}-1}=0$ (see \eqref{T l 0}-\eqref{T l 2}).
To describe the first Melnikov's condition, we take
\begin{equation*}
\begin{aligned}
\Lambda_{l_{*}}^{(2)}=\{\scr{J}\in\tilde\Lambda: &
|\langle k,\omega_{0}(\xi)\rangle+\Omega_{j}|>\sqrt{\epsilon} |k|^{-\tau}~\\
&\textrm{holds for all}~1\leq j\leq n, |k|\leq A^{l_{*}}
~\textrm{and all }~\xi\in\scr{J}\}.
\end{aligned}
\end{equation*}
For $\xi\in\cup_{\scr{J}\in\Lambda_{l_{*}}^{(2)}} \scr{J}$, there is
\begin{equation}\label{D -1}
  |D_{l_{*}-1}(j,k)|^{-1}\leq A^{\tau l_{*}}\epsilon^{-1/2}
\end{equation}
and hence the diagonal matrix
$
  \|D_{l_{*}-1}^{-1}\|\leq A^{\tau l_{*}}\epsilon^{-1/2}.
$
Observing that  $\|S_{l_{*}-1}\|\lesssim \epsilon$,
we obtain from the Neumann series that the inverse
matrix $G_{l_{*}-1; A^{l_{*}}}$ of $T_{l_{*-1}; A^{l_{*}}}$
satisfies
\begin{equation}\label{G l_0}
  \|G_{l_{*}-1;A^{l_{*}}}\|\leq 2 A^{\tau l_{*}}\epsilon^{-1/2}< A^{l_{*}^{C_{1}}}
\end{equation}
if we take
$l_{*}=l_{*}(\epsilon)\sim \log_{A}\frac{1}{\epsilon}
$ (more precisely, 	$A^{\tau l_{*}}= \epsilon^{-1/3}$).
Moreover, there is
\begin{equation}\label{G l_0 decay}
  |G_{l_{*}-1;A^{l_{*}}}(k,k')|< e^{-s_{0}|k-k'|},\quad \textrm{for} ~|k|,|k'|\leq A^{l_{*}}.
\end{equation}
Similarly, letting
\begin{equation*}
\begin{aligned}
	\Lambda_{l_{*}}^{(3)}=\{\scr{J}\in\tilde\Lambda: &
|\langle k,\omega_{0}(\xi)\rangle+\Omega_{j_{1}}+\Omega_{j_{2}}|>\sqrt{\epsilon} |k|^{-\tau}~\\
&\textrm{holds for all}~1\leq j_{1}, j_{2}\leq n, |k|\leq A^{l_{*}}
~\textrm{and for all }~\xi\in\scr{J}\},
\end{aligned}
\end{equation*}
we have \eqref{G l_0} and \eqref{G l_0 decay} hold
on $\cup_{\scr{J}\in\Lambda_{l_{*}}^{(3)}}\scr{J}$
when replacing $G_{l_{*}-1;A^{l_{*}}}$
by $\textbf{G}_{l_{*}-1;A^{l_{*}}}$.
There is also the measure estimate as that in
\eqref{0 mes} for $\Lambda_{l_{*}}^{(2)}$
and $\Lambda_{l_{*}}^{(3)}$, which implies
\begin{equation*}
	\textrm{mes}~\left(\bigcup_{\scr{J}\in\tilde{\Lambda}\setminus (
	\Lambda_{l_{*}}^{(1)}\cap
	\Lambda_{l_{*}}^{(2)}\cap
	\Lambda_{l_{*}}^{(3)} )} \scr{J}\right)\lesssim \sqrt{\epsilon}< A^{-(\log l_{*})^{C_{4}}}.
\end{equation*}
Then, taking
\begin{equation*}
	\Pi_{l_{*}}=\bigcup_{\scr{J}\in \Lambda_{l_{*}}^{(1)}\cap
	\Lambda_{l_{*}}^{(2)}\cap
	\Lambda_{l_{*}}^{(3)}} \scr{J},
\end{equation*}
we obtain the desired parameter set in the statement $(l_{*}.4)$.
This verifies the first step of the iteration in the Iterative Lemma.

Letting $\Pi_{\infty}=\cap_{l\geq l_{*}} \Pi_{l}$, the convergence of the iteration on the uniform domain $\scr{D}(\frac{s_{0}}{2},\frac{r_{0}}{2})\times \Pi_{\infty}$
is standard and we omit the details.
\qed

\section{Green's function estimate}

This section is devoted to the proof of Lemma \ref{Green}.
in which, for simplicity, we have dropped
the iterative subscript $m$ for some expressions.
For reader's convenience, we recall and explain
some notations at the beginning.

Recall that
\begin{equation*}
\varepsilon=A^{-(\frac{4}{3})^{m}},\quad N=A^{m+1},\quad
\Pi=\Pi_{m}=\bigcup_{\scr{J}\in\Lambda_{m}}\scr{J},\quad
(r,s)=(r_{m},s_{m}), \quad\omega=\omega_{m}.	
\end{equation*}
The matrix $T=T_{m}$ in Lemma \ref{Green} (depending
on the $m$-th iteration) is defined by
\begin{equation*}
	T=D+S, \quad D=D_{m},\quad S=S_{m},
\end{equation*}
where the diagonal matrix
\begin{equation}\label{D m}
	D_{m}(j,k)=\Omega_{j}+\langle k,\omega_{m}\rangle
\end{equation}
and non-diagonal matrix
\begin{equation}\label{S m}
	S_{m}((j,k),(j',k'))=(B_{m;jj'}+R^{z\bar{z}}_{m;jj'})^{\wedge}(k-k'),
\end{equation}
with $(\cdot)^{\wedge}(k)$ being the $k$-th Fourier
coefficient of the associated function.
In what follows, we shall also consider   those matrices  depending on
the $l$-th iteration. To make a distinction, we recall
\eqref{T l 0}-\eqref{T l 2} that
\begin{equation}\label{T l}
T_{l}= D_{l}+S_{l},\quad l_{*}\leq l\leq m,
\end{equation}
where $D_{l}$ and $S_{l}$ are defined by
\eqref{D m} and \eqref{S m} upon replacing $m$ by $l$,
respectively.

For any set $U\subset\ZZ^{d}$, we denote by
$T_{l;U}$ the restriction of $T_{l}$ on
$\{1,\cdots,n\}\times U$. For any integer $M>0$, we
write $T_{l;M}=T_{l; [-M,M]^{d}\cap \ZZ^{d}}$
by some abuse of notation.
As a result,  $T_{N}=T_{m;N}$ in Lemma \ref{Green} denotes
the restriction of $T_{m}$ on
$\{1,\cdots,n\}\times ([-A^{m+1},A^{m+1}]^{d}\cap\ZZ^{d})$.
Our goal in this section is to construct and control the inverse of the matrix $T_{N}$, i.e.,
to establish the Green's function estimate for $T_{N}$.

By the definition of $S_{l}$ in \eqref{S m} (replacing
$m$ by any $l_{*}\leq l\leq m$), one readily sees that $S_{l}$ is a Toeplitz matrix with respect to the indices $k,k'$ in $\ZZ^{d}$, i.e.,
\begin{equation}\label{toeplitz}
  S_{l}((j,k+p),(j',k'+p))= S_{l}((j,k),(j',k'))
\end{equation}
for any $p\in\ZZ^{d}$. Moreover,
since $B_{l}(x)$ and $R_{l}^{z\bar{z} }(x)=\partial_{\bar{z}}\partial_{z} P_{l}^{\textrm{low}}(x,0,0,0)$ are analytic
in $x\in\TT^{d}_{s_{l}}$, we have
\begin{equation*}
|\widehat{B}_{l}(k)|\lesssim \epsilon e^{-s_{l} |k|}, \quad |\widehat{R^{z\bar{z}}}(k)|<  \frac{1}{r_{l}}
\sup_{\mathscr{D}(s_{l},r_{l})} |\partial_{z} P_{l}^{\textrm{low}}| e^{- s_{l} |k|}< \varepsilon_{l}
 e^{-s_{l} |k|}.
\end{equation*}
Consequently,
the matrix $S_{l}$ enjoys the off-diagonal exponential decay
\begin{equation}\label{off decay}
  |S_{l}((j,k),(j',k'))|\lesssim \epsilon e^{-s_{l}|k-k'|}.
\end{equation}

Throughout the proof of Lemma \ref{Green}, one easily finds that
the spatial indices $j, j'$ play seldom role in establishing
the Green's function estimate, except those estimates involving absolute constants
 depending only on $n$. For that reason, we omit the finite indices and write
$S(k,k')= S((j,k), (j',k'))$ for simplicity.

Now we give an outline of the construction and estimate of the Green's function
$G_{N}=T_{N}^{-1}$.
By the Iterative Lemma, we are able to obtain the Green's function estimates for
$G_{K}$ with $K\sim (\log N)^{C}$. Then we shall apply the large deviation estimate
to establish the Green's function estimate for all $G_{k_{0}+[-M_{0},M_{0}]^{d}}$
with $K/2\leq k_{0}\leq N $ and $M_{0}\sim (\log N)^{C}$, in which parameter
exclusion should be taken care of by the semialgebraic sets arguments.
Finally, we employ a coupling lemma with two scales ($K$ and $M_{0}$) to
prove \eqref{G N decay} for $G_{N}$, in which one should be careful on the loss of the decay rate.

Due to the  rapid convergence of the Newton iteration, we can  study those
\emph{suitable} matrices $T_{l}$ with $l<m$ and work out the
Green's function estimate for them. Then $G_{K}$ and $G_{k_{0}+[-M_{0},M_{0}]^{d}}$
can be derived directly from $G_{l;K}$ and $G_{l'; k_{0}+[-M_{0},M_{0}]^{d}}$
by employing the Neumann series.

We organize this section as follows. In subsection \ref{preli},
we give some auxiliary lemmas, which are frequently used in this section.
In subsection \ref{LDT main}, we employ the large deviation theorem
and the multiscale analysis method to establish the estimate of the Green's function
 $G_{M_{0}}^{\sigma}$.
In subsection \ref{meas}, we employ the semialgebraic set method
to give the measure estimate and obtain the desired parameter $\Pi_{m+1}$
in the Iterative Lemma.
Finally, we apply the coupling lemma to prove the estimate of the Green's function
 $G_{N}$, which completes the proof of Lemma \ref{Green}.

\subsection{Preliminary}\label{preli}
We first give a quantitative lemma here based on the Neumann series, which is
frequently used throughout this section.
It is worthy mentioning that
the matrix $\mc{T}$ , the integer
$\mc{N}$ and $\pmb{\epsilon}$ in Lemma \ref{Neumann} are arbitrary and
independent of the KAM iterations.

\begin{lemma}\label{Neumann}
Let $U\subset\mathbb{Z}^d$ satisfy the diameter $|U|= \mathcal{N}>0$ and let $\mathcal{T}, \mc{T}'$  be two linear operators  on $\ell^2(\mathbb{Z}^d)$. Denote $\mc{T}_{U}=R_U \mc{T} R_U$ with $R_{U}$ being the restriction operator on $U$. Let further $\alpha>0$, $\rho>0$ and
 $0<b< \theta<1$.

 Assume the following properties hold.
 \begin{enumerate}[(i)]

  \item $\mc{G}_{U}=\mc{T}_{U}^{-1}$ admits the Green's function estimate
 \begin{align*}
& \|\mc{G}_{U}\|\leq e^{\mc{N}^{b}},\\
& |\mc{G}_{U}(x,y)|\leq e^{-\alpha |x-y|}
\quad\emph{\textrm{for}}~ |x-y|>\mc{N}^{\theta}.
 \end{align*}

 \item
For all $x,y\in U$, $$|(\mc{T}_{U}'-\mc{T}_{U})(x,y)|\leq
\pmb{\epsilon } e^{-\rho |x-y|}.$$
\end{enumerate}

Then, if
$
  \pmb{\epsilon}< e^{-4 \rho \mc{N}^{\theta} },
$ we have
\begin{align*}
&\|\mc{G}_{U}'\|\leq 2\|\mc{G}_{U}\|,\\
& |\mc{T}'_{U}(x,y)|\leq 2 e^{-(\alpha\wedge \rho) |x-y|},\quad
\emph{\textrm{for}}~|x-y|>\mc{N}^{\theta},
 \end{align*}
 where $\alpha\wedge\rho=\min\{\alpha,\rho\}$.
\end{lemma}

\noindent\textbf{Proof.}
It is easy to see
$\mc{T}_{U}'=\mc{T}_{U}(Id+\mc{G}_{U}(\mc{T'}_{U}-\mc{T}_{U}))$
and we write $$\Delta=\mc{G}_{U} (\mc{T}'_{U}-\mc{T}_{U}).$$ Then by assumptions, $\|\Delta\|\leq {1}/{2}$,
which together with Neumann series argument implies
$\|\mc{G}'_{U}\|\leq 2 \|\mc{G}_{U}\|.$ For any integer
$s\geq 1$, we compute
\begin{equation*}
  \begin{aligned}
    \Delta^{s}(x,y)=&\sum_{k_{1},\cdots, k_{s-1}\in U}\Delta(k_{0},k_{1})\Delta(k_{1},k_{2})\cdots
    \Delta(k_{s-1},k_{s})\\
    =& \sum_{k_{1},\cdots,k_{s-1}, l_{0},\cdots, l_{s-1}\in U}
    \prod_{j=0}^{s-1} \mc{G}_{U}(k_{j},l_{j}) (\mc{T}'-\mc{T})(l_{j},k_{j+1}),
  \end{aligned}
\end{equation*}
where $k_{0}=x$ and $k_{s}=y$.
If $|k_{j}-l_{j}|> \mc{N}^{\theta}$, there is
\begin{equation*}
|  \mc{G}_{U}(k_{j},l_{j})|\cdot |(\mc{T}'-\mc{T})(l_{j},k_{j+1})| < \pmb{\epsilon}
e^{-(\alpha\wedge\rho) |k_{j}-k_{j+1}|},
\end{equation*}
and if $|k_{j}-l_{j}|\leq \mc{N}^{\theta}$, there is
\begin{equation*}
  |  \mc{G}_{U}(k_{j},l_{j})|\cdot |(\mc{T}'-\mc{T})
  (l_{j},k_{j+1})|
< \pmb{\epsilon} e^{\mc{N}^{b}+\rho \mc{N}^{\theta} -\rho |k_{j}-k_{j+1}|}.
\end{equation*}
It follows from $\pmb{\epsilon}< e^{-4\rho \mc{N}^{\theta}}$ that
\begin{equation*}
  |\Delta^{s}(x,y)|<(C \mc{N})^{2 d s} {\pmb{\epsilon}}^{s}
   e^{s(\mc{N}^{b}+\rho \mc{N}^{\theta} )} e^{-\rho |x-y|}
< e^{-2\rho \mc{N}^{\theta} s} e^{-\rho |x-y|}
\end{equation*}
and thus
\begin{equation*}
 \left | \sum_{s=1}^{\infty} \Delta^{s}(x,y)\right|< 2 e^{-2\rho \mc{N}^{\theta}}
 e^{-\rho |x-y|}.
\end{equation*}

Finally, for any $x,y\in U$ we have
\begin{equation*}
  \begin{aligned}
    |\mc{G}_{U}'(x,y)|<&|\mc{G}_{U}(x,y)|+
    \sum_{l\in U}|\sum_{s=1}^{\infty} \Delta^{s}(x,l)|
    \cdot |\mc{G}_{U}(l,y)|\\
<& |\mc{G}_{U}(x,y)|+
    \sum_{l\in U, |l-y|>\mc{N}^{\theta}} |\sum_{s=1}^{\infty} \Delta^{s}(x,l)|
    \cdot |\mc{G}_{U}(l,y)|\\
&\quad +\sum_{l\in U, |l-y|\leq \mc{N}^{\theta}} |\sum_{s=1}^{\infty} \Delta^{s}(x,l)|
    \cdot |\mc{G}_{U}(l,y)|\\
<& |\mc{G}_{U}(x,y)|+ (C\mc{N})^{d} e^{-2\rho \mc{N}^{\theta}}
    e^{-\rho|x-y|}+ (C\mc{N})^{d} e^{-2\rho \mc{N}^{\theta}} e^{\mc{N}^{b}+\rho \mc{N}^{\theta}}
    e^{-\rho |x-y|}\\
<& |\mc{G}_{U}(x,y)|+\frac{1}{10} e^{-\rho |x-y|}.
  \end{aligned}
\end{equation*}
As a result, whenever $|x-y|> \mc{N}^{\theta}$,
\begin{equation*}
  |\mc{G}_{U}'(x,y)|< 2 e^{-(\alpha\wedge \rho) |x-y|}.
\end{equation*}
This completes the proof. \qed

\bigskip

Next we describe quantitatively the variation of $T_{l}$, which enables us to
apply Lemma \ref{Neumann} in what follows.

\begin{lemma}\label{variation}

  Let $l_{*}\leq l<l'\leq m$ and consider the linear operator $T_{l}, T_{l'}$ defined in \eqref{T l}.
  Let further $T=T_{l;A^{l'}}$ and $T'=T_{l'; A^{l'}}$ be the restriction of $T_{l}, T_{l'}$ on
  $[-A^{l'},A^{l'}]^{d}$.
  Then,
  we have
  \begin{equation*}
    |(T'-T)(k,k')|\lesssim A^{l'}\cdot \epsilon_{l}^{1/10}\exp(-s_{l'}~ |k-k'|).
  \end{equation*}
\end{lemma}
\noindent\textbf{Proof.}~
By definition we have
\begin{equation*}
(T_{l'}-T_{l})(k,k')=\langle\omega_{l'}-\omega_{l}, k\rangle\cdot \delta_{k k'}+ (B_{l'}-B_{l})^{\wedge}(k-k')
+ (R^{z\bar{z}}_{l'}-R^{z\bar{z}}_{l})^{\wedge}(k-k'),
\end{equation*}
where $\delta_{k k'}$ equals to one if $k=k'$ and vanishes otherwise.
By  the Iterative Lemma, there is
\begin{equation*}
|\omega_{l'}-\omega_{l}|\lesssim \epsilon_{l}^{1/10}.
\end{equation*}
Moreover,
\begin{equation*}
B_{l'}-B_{l}=\sum_{r=l}^{l'-1}R_{r}^{z\bar{z}}+\{N_{r},F_{r}\}^{z\bar{z}}+\{P^{\textrm{high}}_{r}, F_{r}\}^{z\bar{z}}
\end{equation*}
and
\begin{equation*}
R_{l'}^{z\bar{z}}-R_{l}^{z\bar{z}}= \partial_{z\bar{z}} P^{\textrm{low}}_{l'} (x,0,0,0)-
\partial_{z\bar{z}} P^{\textrm{low}}_{l} (x,0,0,0).
\end{equation*}

The property ${~}_{r_{l}} \pmb{|} X_{P_{l}^{\textrm{low}}}\pmb{|}_{\mathscr{D}(s_{l},r_{l})}
<\epsilon_{l}$
ensures
\begin{equation*}
\sup_{x\in \TT^{d}_{s_{l'}}}|R_{l'}^{z\bar{z}}(x)-R_{l}^{z\bar{z}}(x)|\lesssim \epsilon_{l}.
\end{equation*}
Since $R^{z\bar{z}}_{l}= \partial_{z\bar{z} } P^{\textrm{low}}_{l}(x,0,0,0)$, there is also
\begin{equation*}
\sup_{x\in\TT^{d}_{s_{l}}} |R^{z\bar{z}}_{l}(x)|\lesssim \epsilon_{l}.
\end{equation*}
and
\begin{equation*}
\sup_{x\in\TT^{d}_{s_{l}^{(1)}}}\left|\{E_{l},F_{l}\}^{z\bar{z}}(x)+\{R^{\textrm{high}}_{l}, F_{l}\}^{z\bar{z}}(x) \right|
\lesssim \epsilon_{l}^{1/3} (l+1)^{C}< \epsilon_{l}^{1/4}.
\end{equation*}
Hence
\begin{equation*}
\sup_{x\in\TT^{d}_{l'}}|B_{l'}(x)-B_{l}(x)|< \epsilon_{l}^{1/10}
\end{equation*}
and the conclusion follows.
\qed

\bigskip

We finally  cite here a decomposition lemma in \cite[Lemma 9.9]{Bou05_Green}.
\begin{lemma}\label{decomposition}
  Let $\mathcal{S}\subset[0,1]^{2n}$ be a semi-algebraic set of degree
  $B$ and $\emph{mes}_{2n} (\mc{S})< \eta$, $\log B\ll \log \frac{1}{\eta}$.
  We denote $(\omega,x)\in[0,1]^{n}\times[0,1]^{n}$ the product variable.
  Fix $\pmb{\epsilon}> \eta^{\frac{1}{2n}}$. Then there is a decomposition
  \begin{equation*}
    \mathcal{S}=\mathcal{S}_{1}\cup \mathcal{S}_{2},
  \end{equation*}
  $\mathcal{S}_{1}$ satisfying
  \begin{equation*}
    |\emph{Proj}_{\omega} \mathcal{S}_{1}|< B^{C} \pmb{\epsilon}
  \end{equation*}
  and $\mathcal{S}_{2}$ satisfying the transversality property
  \begin{equation*}
    \emph{mes}_{n} (\mathcal{S}_{2}\cap L)< B^{C} \pmb{\epsilon}^{-1} \eta^{\frac{1}{2n}}
  \end{equation*}
  for any $n$-dimensional hyperplane $L$ such that
  $\max_{0\leq j\leq n-1} |\emph{Proj}_{L}(e_{j})|<\frac{1}{100}\pmb{\epsilon}$
  (we denote $(e_{0}, \cdots, e_{n-1})$ the $\omega$- coordinate vectors.)
\end{lemma}

\subsection{Estimate of $G^{\sigma}_{M_{0}}$}\label{LDT main}

In this part, our goal is to establish the following type of  Green's function estimate
\begin{equation}\label{G M 0}
  \begin{aligned}
& \|G_{ U(k_{0})}\|< e^{M_{0}^{b}},\\
& |G_{ U(k_{0})} (k,k') |< e^{-\alpha''|k-k'|}\quad
\textrm{for}~  |k-k'|> M_{0}^{\theta},
\end{aligned}
\end{equation}
for all
$T_{U(k_{0})}$ with $K/2\leq |k_{0}|\leq N$, where
$0<b<\theta<1$, $\alpha''>0$ is to be specified, $U(k_{0})=k_{0}+[-M_{0},M_{0}]^{d}$ and
\begin{equation*}
M_{0}= (\log N)^{C_{0}},\quad \log K=(\log M_{0})^{C_{7}}.
\end{equation*}

As mentioned before, we shall work on some $T_{l_{0}; U(k_{0})}$
 with $l_{0}<m$ rather than on
$T_{m;U(k_{0})}$ directly, due to the rapid
convergence of the Newton iteration. This can be resolved by a simple
application of  Neumann series (see Lemma \ref{Neumann}).
Indeed, choosing
\begin{equation}\label{l 0 M 0}
 l_{0}=C_{8} \log M_{0},\quad \textrm{with}~C_{8}> \frac{1+\log 10}{\log \frac{4}{3}},
\end{equation}
we see from Lemma \ref{variation} that
\begin{equation*}
  |T_{U(k_{0})}(k,k')-T_{l_{0}; U(k_{0})}(k,k')|\lesssim \epsilon_{l_0}^{1/10}\cdot N~ \exp (-s |k-k'|),\quad
  s=s_{m}.
\end{equation*}
Suppose \eqref{G M 0} is valid for $G_{l_{0};U(k_{0})}$, then, by
verifying
\begin{equation*}
\epsilon_{l_0}^{1/10}\cdot N < \epsilon_{l_{0}}^{1/20}<\frac{1}{100}
 e^{-M_{0}^{\theta}},
\end{equation*}
it follows from
Lemma \ref{Neumann} that Green's function estimate \eqref{G M_0} also
holds for $G_{U(k_{0})}$   up to a constant multiplication.
To this end, we shall establish \eqref{G M 0} for $G_{l_{0};U(k_{0})}$
in what follows.

Recalling the Toeplitz property \eqref{toeplitz} for $T_{l}$, $l_{*}\leq l\leq m$,
 we denote
\begin{equation}\label{T l sigma}
  T^{\sigma}_{l}= D_{l}^{\sigma}+S_{l},
\end{equation}
where  $S_{l}$ is defined in \eqref{S m} (replacing
$m$ by $l$) and $D_{l}^{\sigma}$ takes the form of
\begin{equation*}
  D^{\sigma}_{l}(j,k)= \sigma+ \langle k,\omega_{l}\rangle+\Omega_{j},\quad 1\leq j\leq n, ~k\in\ZZ^{d}.
\end{equation*}
Observe  by the Toeplitz property that
\begin{equation*}
  T_{l_{0}; U(k_{0})}=
  T_{l_{0}; U(k_{0})}^{\sigma=0}=
  T_{l_{0}; M_{0}}^{\sigma=\langle k_{0},\omega_{l_{0}}\rangle}.
\end{equation*}
Then it suffices to establish the Green's function estimate of
$T^{\sigma}_{l_{0}; M_{0}}$ for
\begin{equation*}
  \sigma\in\{\langle k,\omega_{l_{0}} \rangle: K/2\leq k\leq N \}.
\end{equation*}

The lemma below is the core of our analysis in this part, which is
independent of the Iterative Lemma and whose proof
is delayed to the appendix. To formulate it,
we need to introduce the elementary regions.
An \emph{elementary region} is defined to be a set $U$ of the form
\begin{equation*}
  U=R\setminus (R+z)
\end{equation*}
where $z\in\ZZ^{d}$  is arbitrary and $R$ is a block in $\ZZ^{d}$, i.e.,
\begin{equation*}
  R= \{y=(y_{1},\cdots, y_{d})\in\ZZ^{d}: y_{i}\in [x_{i}-M_{i}, x_{i}+M_{i}], i=1\cdots, d\}.
\end{equation*}
The diameter of an elementary region $U$ is  denoted by $|U|$. The set of all
elementary regions of diameter $M$ is denoted by $\mathcal{ER}(M)$.
The class of elementary regions consists of $d$-dimensional rectangles, L-shaped regions and
$(d-1)$ -dimensional rectangles with normal vector parallel to the
axis.

\begin{lemma}\label{LDT}
Consider the matrix $\mc{T}^{\sigma}=\mc{D}^{\sigma}+ \mc{S}$,
where $\sigma\in\RR$ and $\mc{D}^{\sigma}$ is a diagonal matrix with
\begin{equation*}\label{D sigma}
  \mc{D}^{\sigma}(j,k)= \langle k,  \pmb{\omega} \rangle+\sigma+\Omega_{j},\quad
  1\leq j\leq n, ~k\in\ZZ^{d},
\end{equation*}
and we omit the finite index $j$ for simplicity.
 Let $\underline{\mc{N}}_{0}, \overline{\mc{N}}_{0}=\underline{\mc{N}}_{0}^{C}$ be sufficiently large and let the various constants below
satisfy
\begin{equation*}
  0<\beta\ll 1,\quad 1-\frac{\beta}{10}<b<\theta<1, \quad \alpha_{0}>0,\quad
  \rho>0.
\end{equation*}
Assume the following properties hold.
\begin{enumerate}[(i)]
  \item The matrix $\mc{S}$ satisfies  the Toeplitz property with respect to the $k$-index and
  \begin{equation*}
  |\mc{S}(x,y)|< \epsilon e^{-\rho |x-y|},\quad x\neq y.
\end{equation*}

\item The frequency $\pmb{\omega}$ satisfies
Diophantine condition
\begin{equation*}
  |\langle k,\pmb{\omega}\rangle|> \frac{\nu}{|k|^{\tau}},\quad 0\neq k\in\ZZ^{d},
  0<\nu<1, \tau>d+1.
\end{equation*}

\item For any $\underline{\mc{N}}_{0}<\mc{N}_{0}< \overline{\mc{N}}_{0}$ and any elementary region $U_{0}\in\mc{ER}(\mc{N}_{0})$, the Green's function estimate\begin{equation}
\begin{aligned}
&\|\mc{G}^{\sigma}_{U_{0}}\|< e^{\mc{N}_{0}^{b}},\\
& |\mc{G}_{U_{0}}^{\sigma}(x,y)|< e^{-\alpha_{0} |x-y|},\quad
  \emph{for}~|x-y|> \mc{N}_{0}^{\theta},
  \end{aligned}
\end{equation}
holds for all $\sigma$ except in a set $\mathscr{E}_{0}(U_{0})$ of measure
at most $ e^{-\mc{N}_{0}^{\beta^{3}}}$.
\end{enumerate}

 Then for any large $\mc{N}>\overline{\mc{N}}_{0}$ and any elementary region $U\in\mc{ER}(\mc{N})$,
the Green's function estimate
\begin{equation}\label{LDT-1}
  \begin{aligned}
&\|\mc{G}^{\sigma}_{U}\|< e^{\mc{N}^{b}},\\
& |\mc{G}_{U}^{\sigma}(x,y)|< e^{-\alpha |x-y|},\quad
  \emph{for}~|x-y|> \mc{N}^{\theta}
  \end{aligned}
\end{equation}
holds for all $\sigma\in\RR$ outside of a set $\mathscr{E}=\mathscr{E}(U)$ with
\begin{equation*}
  \emph{\textrm{mes}}~(\mathscr{E})< e^{-\mc{N}^{\beta^{3}}},
\end{equation*}
where $\alpha>(\alpha_{0}\wedge\rho)-(\log \mc{N}_{0})^{-8}$.
\end{lemma}

Now we apply Lemma \ref{LDT} to prove the following induction statements.

\begin{proposition}\label{indu state}
 Under the assumptions of Lemma \ref{iter lemma}, we consider a family of
  matrices $T_{l}^{\sigma}$ defined by \eqref{T l sigma}.
  Let $q(l)=  \frac{\log \frac{4}{3}}{2\log A} l$.
  Then for all $l_{*}\leq l\leq m$ and any elementary
  region $U\in\mc{ER}(A^{q(l)})$, there is
  \begin{equation}\label{q l}
  \begin{aligned}
& \|G_{l; U}^{\sigma}\|< e^{A^{q(l) b}},\\
& |G_{l; U}^{\sigma}(k,k')|< e^{-\alpha'(l) |k-k'|},\quad\emph{for}~ |k-k'|> A^{q(l)\theta},
  \end{aligned}
  \end{equation}
  for all $\sigma$ except in a set $\scr{E}_{l}=\scr{E}_{l}(U)$ with $\emph{mes}~ (\scr{E}_{l})<
  e^{-A^{q(l) \beta^{3}}}$, where $\alpha'(l)>s_{l+1}$.
\end{proposition}

\noindent\textbf{Proof.} We prove the proposition by the method of inductions on $l$. The initial steps $(l_{*}\leq l\leq l_{*}^{C})$ are essentially a direct application of the Neumann series provided
the perturbation is small enough and we omit it. See also the similar arguments in subsection \ref{proof main}.

Assume by induction that the property \eqref{q l} holds with
$l<m$. We need to
establish \eqref{q l} for $l+1$ and any $U\in\mc{ER}(A^{q(l+1)})$.
Observe first by similar computations in Lemma \ref{variation} that, for any
$V\in\mc{ER}(A^{q(l)})$, there is
\begin{equation*}
  |(T_{l+1;V}^{\sigma}-T_{l;V}^{\sigma})(k,k')|< \epsilon_{l}^{1/20}
  \exp (-s_{l+1} |k-k'|).
\end{equation*}
Since $\epsilon_{l}^{1/20}=A^{-(\frac{4}{3})^{l}\frac{1}{20}}< e^{-A^{b q(l)}}$ holds with
large $A$, it follows from Lemma \ref{Neumann} that
\begin{equation*}
  \begin{aligned}
& \|G_{l+1; V}^{\sigma}\|<  e^{-A^{q(l) b}},\\
& |G_{l+1; V}^{\sigma}(k,k')|< e^{-s_{l+1} |k-k'|}\quad\textrm{for}~ |k-k'|> A^{q(l)\theta},
  \end{aligned}
  \end{equation*}
  essentially holds for all $\sigma$ except in a set $\scr{E}_{l}$ with $\textrm{mes}~ (\scr{E}_{l})<
  e^{-A^{q(l)} \beta^{3}}$.
  Now we apply Lemma \ref{LDT} by taking
  \begin{equation*}
    \mc{T}^{\sigma}= T^{\sigma}_{l+1},\quad \mc{N}= A^{q(l+1)},\quad \overline{\mc{N}_{0}}= A^{q(l)},\quad \rho= s_{l+1},\quad \pmb{\omega}= \omega_{l+1}.
  \end{equation*}
Then we obtain \eqref{q l} for $l+1$ with $\alpha'(l+1)> s_{l+1}-(\log A^{q(l)})^{-8}
>s_{l+2}$.
This completes the proof of the induction statements.\qed

\bigskip

Recall that $l_{0}$ and $M_{0}$ are fixed
in \eqref{l 0 M 0} (depending on $N$). Back to $T_{l_{0}; M_{0}}^{\sigma}$, we have

\begin{proposition}\label{g l 0 m}
Under the assumption of the {Iterative Lemma \ref{iter lemma}}, we have
  \begin{equation}\label{G M_0}
  \begin{aligned}
& \|G^{\sigma}_{l_{0}; M_{0}}\|< e^{M_{0}^{b}},\\
& |G^{\sigma}_{l_{0};M_{0}}(k,k')|< \exp(-\alpha'' |k-k'|)\quad\emph{for}~|k-k'|>M_{0}^{\theta}.
  \end{aligned}
\end{equation}
except for all $\sigma\in\RR$ outside of a set $\mathscr{E}_{M_{0}}$ with
$\emph{mes}~(\mathscr{E}_{M_{0}})< e^{-M_{0}^{\beta^{3}}}$,
where $\alpha''>s$ and $0<\beta\ll 1$.
\end{proposition}

\noindent\textbf{Proof.}~
For fixed $l_{0}=l_{0}(N)\sim \log (m+1)$, we define $\mc{N}_{0}, \underline{l_{0}}'$ and  $\underline{l_{0}}$ in order
as follows
\begin{equation}\label{N_0}
\mc{N}_{0}=2 \exp (l_{0}^{1/4}),\quad  \mc{N}_{0}=A^{\underline{l_{0}}'},\quad
\underline{l_{0}}=\frac{2\log A}{\log \frac{4}{3}}\underline{l_{0}'}\sim \log \mc{N}_{0} .
\end{equation}
By Proposition \ref{indu state}, for $\underline{l_{0}}\sim (\log (m+1))^{1/4}$ and $\mc{N}_{0}$ satisfying
\eqref{N_0}, there is, for any $U_{0}\in\mc{ER}(\mc{N}_{0})$,
the estimate
\begin{equation}\label{G l 0 U 0}
  \begin{aligned}
& |G_{\underline{l_{0}}; U_{0}}^{\sigma}|< e^{N_{0}^{b}},\\
& |G^{\sigma}_{\underline{l_{0}};U_{0}}(k,k')|< e^{-\alpha'(\underline{l_{0}})|x-y|},\quad \textrm{for}~
    |k-k'|> \mc{N}_{0}^{\theta}
  \end{aligned}
\end{equation}
holds for all $\sigma$ except in a set $\mathscr{E}_{\mc{N}_{0}}$ with
$\textrm{mes}~(\mathscr{E}_{\mc{N}_{0}})< e^{-\mc{N}_{0}^{\beta^{3}}}$,
where $\alpha'(\underline{l_{0}})> s_{\underline{l_{0}}+1}$.
It then follows from Lemma \ref{Neumann} that
the Green's function estimate
\eqref{G l 0 U 0} essentially holds  when replacing
$G_{\underline{l_{0}}; U_{0}}^{\sigma}$ by $G_{l_{0}; U_{0}}^{\sigma}$,
since
$
  \epsilon_{\underline{l_{0}}} e^{\mc{N}_{0}^{\theta}}= A^{-(\frac{4}{3})^{\underline{l_{0}}}}\cdot
  e^{A^{\theta \underline{l_{0}'}}}<1.
$
By $(l_{0}.4.2)$ in the Iterative Lemma, we have
\begin{equation}\label{dioph}
  |\langle k,\omega_{l_{0}}(\xi)\rangle|>\frac{\nu}{|k|^{\tau}},\quad 0\neq |k|\leq M_{0}, \quad k\in\ZZ^{d},
\end{equation}
where $\nu\sim\sqrt{\epsilon}$ and $\tau>d+1$.
Taking
\begin{equation*}
	\mc{T}^{\sigma}= T^{\sigma}_{l_{0}},\quad \mc{N}= M_{0},\quad \rho= s_{l_{0}},\quad
  \pmb{\omega}=\omega_{l_{0}}(\xi),
\end{equation*}
and applying Lemma \ref{LDT}, we have
\begin{equation}\label{G M_0}
  \begin{aligned}
& \|G^{\sigma}_{l_{0}; M_{0}}\|< e^{M_{0}^{b}},\\
& |G^{\sigma}_{l_{0};M_{0}}(k,k')|< \exp(-\alpha'' |k-k'|)\quad\textrm{for}~|k-k'|>M_{0}^{\theta},
  \end{aligned}
\end{equation}
except for all $\sigma\in\RR$ outside of a set $\mathscr{E}_{M_{0}}$ with
$\textrm{mes}~(\mathscr{E}_{M_{0}})< e^{-M_{0}^{\beta^{3}}}$,
where
\begin{equation*}
\alpha''=(\alpha'\wedge s_{l_{0}})-
    (\log \mc{N}_{0})^{-8}>s_{l_{0}}-(\log \mc{N}_{0})^{-8}> s_{m^{1/10}}>s=s_{m}.
\end{equation*}
This completes the proof. \qed
%
%

\subsection{Elimination of $\sigma$ and measure estimate}\label{meas}
In this part, we shall eliminate the additional parameter $\sigma$ and establish
the Green's function estimates for all $G_{l_{0}; U(k_{0})}$ with
$K/2\leq |k_{0}|\leq N$ and
$U(k_{0})=(k_{0}+[-M_{0},M_{0}]^{d})\cap\ZZ^{d}$. This requires a further parameter exclusion,
whose measure is estimated by the decomposition theorem for semialgebraic
sets. For that reason, we need to give
a semialgebraic description for the breakdown of the Green's function estimate.
The main result in this part is presented below.

\begin{lemma}\label{M_0 decay}
Under the assumption of Lemma \ref{iter lemma},
there exists a measurable set $\Pi_{+}^{(1)}\subset\RR^{d}$ such that for all
$\xi\in\scr{O}(\Pi_{+}^{(1)},A^{-(m+1)^{C_{3}}})$ there is
\begin{equation}\label{G M_0 k}
\begin{aligned}
& \|G_{l_{0}; U(k_{0})}\|< e^{M_{0}^{b}},\quad 0<b<1,\\
& |G_{l_{0}; U(k_{0})} (k,k') |< e^{-s|k-k'|},\quad
\emph{for}~ |k-k'|> M_{0}^{\theta},
\end{aligned}
\end{equation}
for all $ K/2 \leq |k_{0}|\leq N$. Moreover, $\Pi_{+}^{(1)}$ is the union of  a family of  disjoint intervals $\scr{J}'$ with  $|\scr{J}'|=A^{-(m+1)^{C_{3}}}$. For each
$\scr{J}'$, there is a unique $\scr{J}\in\Lambda_{m}$ such that $\scr{J'}\subset \scr{J}$.
The total removed set satisfies
\begin{equation*}
  \emph{mes}~(\Pi\setminus \Pi_{+}^{(1)})<\frac{1}{3} A^{-(\log (m+1))^{C_{4}}}.
\end{equation*}
\end{lemma}

\begin{remark}\label{Pi +}
In the proof of Lemma \ref{M_0 decay}, we shall
pave $\Pi$ into a collection of intervals of diameter
$A^{-(m+1)^{C_{3}}}$ , with
the shrunken parameter set   $\Pi_{+}^{(1)}\subset\Pi$ being a sub-collection.
  When solving \eqref{homo 4}-\eqref{homo 4'} by the same method, we would
  obtain another set $\Pi_{+}^{(2)}$ which is also a sub-collection of the
  $A^{-(m+1)^{C_{3}}}$-intervals paving $\Pi$.
  Then we obtain the desired set $\Pi_{+}=\Pi_{+}^{(1)}\cap \Pi_{+}^{(2)}$, which
  satisfies $\emph{mes} ~(\Pi\setminus \Pi_{+})< A^{-(\log (m+1))^{C_{4}}}$.
  Moreover, $\Pi_{+}$ is a $\Lambda_{+}$ collection of disjoint $A^{-(m+1)^{C_{3}}}$-intervals.
  For each $\scr{J}'\in\Lambda_{+}$, there is a  unique $\scr{J}\in\Lambda$ such
  that $\scr{J}'\subset\scr{J}$.
\end{remark}

\noindent\textbf{Proof.}~We divide the proof into three steps. Step one is devoted
to the truncation of parameters in the Green's function estimate,
which enables us to make a semialgebraic description.
Step two is devoted to the elimination of the additional parameter $\sigma$.
Step three is devoted to the construction of the desired parameter
set and establishing the associated measure estimate.

\emph{Step one.}
From the Iterative  Lemma, we know that $T_{l_{0}}^{\sigma}$ is analytic in
$\xi\in\scr{O}_{l_{0}}=\scr{O}(\Pi_{l_{0}}, A^{-l_{0}^{C_{3}}})$ with
$\Pi_{l_{0}}=\cup_{\scr{J}\in \Lambda_{l_{0}}} \scr{J}$ and $|\scr{J}|=A^{-l_{0}^{C_{3}}}$.
Fix any $\scr{J}\in\Lambda_{l_{0}}$ and denote the center of $\scr{J}$ by
$\xi_{0}$.
Recall that $T^{\sigma}_{l_{0};M_{0}}= D^{\sigma}_{l_{0}}+ S_{l_{0}}$ with
\begin{equation*}
  D^{\sigma}_{l_{0}}(j,k;\xi)=\sigma+\langle k,\omega_{l_{0}}(\xi)\rangle+\Omega_{j},
\end{equation*}
and
\begin{equation*}
  S_{l_{0}}(k,k')= (B_{l_{0}}(\xi)+R^{z\bar{z}}_{l_{0}}(\xi))^{\wedge} (k-k').
\end{equation*}
Let
$
  p= A^{l_{0}^{C_{5}}}
$
with $C_{5}$ given by \eqref{C}.
By Taylor's formula, we denote
\begin{equation*}
  \begin{aligned}
& \omega_{l_{0}}^{\leq }(\xi)=\sum_{i\leq p} \frac{\omega_{l_{0}}^{(i)}(\xi_{0})}{i!}
     (\xi-\xi_{0})^{i},\\
& B_{l_{0}}^{\leq}(\xi)=\sum_{i\leq p} \frac{B_{l_{0}}^{(i)}(\xi_{0})}{i!}
     (\xi-\xi_{0})^{i},\quad
     R_{l_{0}}^{z\bar{z};\leq}(\xi)=\sum_{i\leq p} \frac{(R^{z\bar{z}}_{l_{0}})^{(i)}(\xi_{0})}{i!}
     (\xi-\xi_{0})^{i},\\
& D^{\sigma;\leq}_{l_{0}}(j,k;\xi)= \sigma+\langle k,\omega_{l_{0}}^{\leq}(\xi)\rangle+\Omega_{j},\quad
     S^{\leq}_{l_{0}}(k,k')= (B_{l_{0}}^{\leq}(\xi)
     +R^{z\bar{z};\leq}_{l_{0}}(\xi))^{\wedge} (k-k'),\\
& T^{\sigma;\leq}_{l_{0}}(k,k')=D^{\sigma;\leq}_{l_{0}}(x)+ S^{\leq}_{l_{0}}(k,k').
  \end{aligned}
\end{equation*}
As a result, $T^{\sigma;\leq}_{l_{0}}$ is a polynomial function in $\xi$, whose degree
\begin{equation*}
  \deg_{\xi} T^{\sigma;\leq}_{l_{0}}(k,k')\leq p.
\end{equation*}
Obviously, the truncation error satisfies
\begin{equation*}
  |(T_{l_{0};M_{0}}^{\sigma}-T^{\sigma;\leq}_{l_{0};M_{0}})(k,k')|\leq
  [M_{0} |\omega_{l_{0}}-\omega_{l_{0}}^{\leq}|+(B_{l_{0}}-B_{l_{0}}^{\leq})+
  (R^{z\bar{z}}_{l_{0}}-R_{l_{0}}^{z\bar{z};\leq})]\cdot \exp(-s_{l_{0}} |k-k'|).
\end{equation*}
For $|\xi-\xi_{0}|\leq \kappa |\mathscr{J}|=\kappa A^{-l_{0}^{C_{3}}}$ with
$\kappa\approx 1$ to be specified, we see from
Cauchy's estimate that
\begin{equation}
 |\omega_{l_{0}}-\omega_{l_{0}}^{\leq}|\leq \sup_{\xi\in\mathscr{J}}|\omega_{l_{0}}|
    \cdot \frac{|\xi-\xi_{0}|^{p+1}}{|\mathscr{J}|^{p}\cdot (|\mathscr{J}|-|\xi-\xi_{0}|)}
    \lesssim \frac{\kappa^{p+1}}{1-\kappa}.
\end{equation}
Since $B_{l_{0}}$ and $R_{l_{0}}^{z\bar{z}}$ stay uniformly bounded in their analytical
domain, there is also
\begin{equation*}
  \sup_{|\xi-\xi_{0}|\leq \kappa |\mathscr{J}|}
  \left\{|B_{l_{0}}-B_{l_{0}}^{\leq}|, |R_{l_{0}}^{z\bar{z}}-R_{l_{0}}^{z\bar{z};\leq}|
  \right\} \lesssim
  \frac{\kappa^{p+1}}{1-\kappa}.
\end{equation*}

On the one hand,
letting
$$\mathscr{V}_{1}=\bigcup_{\scr{J}\in\Lambda_{l_{0}}}
\{\xi\in\mathscr{J}:
 \kappa A^{-l_{0}^{C_{3}}}<|\xi-\xi_{0}|\leq A^{-l_{0}^{C_{3}}} \}
 \subset \RR^{d},$$
 we have
\begin{equation*}
  \textrm{mes}~ \mathscr{V}_{1}\lesssim_{d} \frac{1}{|\scr{J}|^{d}} \cdot (1-\kappa) |\scr{J}|
  \lesssim (1-\kappa) A^{(d-1) l_{0}^{C_{3}}}.
\end{equation*}
Taking
\begin{equation*}
  \kappa=1- A^{-(\log (m+1))^{C_{6}}}
\end{equation*}
with $C_{6}$ given by \eqref{C},
then
\begin{equation*}
  \textrm{mes}~ (\mathscr{V}_{1})\lesssim A^{-(\log (m+1))^{C_{6}}} A^{(d-1) (\log (m+1))^{C_{3}}}
< \frac{1}{100} A ^{-  (\log (m+1))^{C_{4}}}
\end{equation*}
provided $m$ is large.

On the other hand, for $|\xi-\xi_{0}|<\kappa |\scr{J}|$, there is
\begin{equation*}
  |(T^{\sigma}_{l_{0};M_{0}}-T^{\sigma;\leq}_{l_{0};M_{0}})(k,k')|\lesssim M_{0}
  \frac{\kappa^{p+1}}
  {1-\kappa} \exp(-s_{l_{0}}(k-k')).
\end{equation*}
By noticing $p=A^{l_{0}^{C_{5}}}$ and $C_{5}=C_{6}+2$, we
get
\begin{equation*}
 M_{0} \frac{\kappa^{p+1}}{1-\kappa}= M_{0}(1-A^{-(\log (m+1))^{C_{6}}})^{p+1}
   A^{(\log (m+1))^{C_{6}}}< e^{-M_{0}}= e^{- C (m+1)^{C_{0}}}.
\end{equation*}
In conclusion, we have
\begin{equation}\label{trunc error}
  |(T^{\sigma}_{l_{0};M_{0}}-T^{\sigma;\leq}_{l_{0};M_{0}})(k,k')|< e^{-M_{0}}
  e^{-s_{l_{0}}|k-k'|}\quad
  \textrm{for}~ |k|\leq M_{0}, |k'|\leq M_{0}.
\end{equation}

\emph{Step two.}
By Lemma \ref{Neumann} and Proposition
\ref{g l 0 m},  we also essentially have\footnote{We omit the constant multiplier
induced by Lemma \ref{Neumann}, which finally can be absorbed by the margins in our
estimates. See Lemma \ref{2 scale couple} for example.}
\begin{equation}\label{G M_0 <}
  \begin{aligned}
& \|G^{\sigma;\leq}_{l_{0}; M_{0}}\|< e^{M_{0}^{b}},\\
& |G^{\sigma;\leq}_{l_{0};M_{0}}(k,k')|< \exp(-\alpha'' |k-k'|)\quad \textrm{for}~
    |k-k'|>M_{0}^{\theta},
  \end{aligned}
\end{equation}
except
for all $\sigma\in\RR$ outside of a set $\mathscr{E}_{M_{0}}$ with
$\textrm{mes}~\mathscr{E}_{M_{0}}< e^{-M_{0}^{\beta^{3}}}$.
Using the formula
$$G^{\sigma;\leq}_{l_{0};M_{0}}(k,k')=
 (T^{\sigma;\leq}_{l_{0};M_{0}})^{*}(k,k')/\det T^{\sigma;\leq}_{l_{0};M_{0}}$$
 with $(\cdot)^{*}$ being the adjoint matrix,
we consider the set  $\mathfrak{S}$ of the triplets  $(\xi,\omega_{l_{0}}^{\leq}, \sigma)$
such that
\begin{equation*}
  |\langle k,\omega^{\leq}_{l_{0}}\rangle|>\frac{\nu}{|k|^{\tau}},\quad 0\neq |k|\leq M_{0},\quad
  k\in\ZZ^{d}
\end{equation*}
and
\eqref{G M_0 <} fails. Obviously, $\mathfrak{S}\subset (\mathscr{J}\cap\RR^{d})\times \RR^{d}\times\RR$ is a semi-algebraic set of
degree at most $M_{0}^{C} p= M_{0}^{C} A^{(\log (m+1))^{C_{5}}}$.
Since $T$ is restricted to $[-M_{0},M_{0}]^{d}$, we may restrict $\sigma$ to be
in $[-CM_{0}, CM_{0}]$. Otherwise, $T$ is diagonal dominated and it suffices
to apply Neumnan series to $T^{\sigma}$ to get the desired estimate.
We decompose $[-C M_{0}, C M_{0}]$ into intervals of length $1$ and identify
each of them with $[0,1]$. Then $\mathfrak{S}$ is divided into $CM_{0}$ sub-intervals
$\mathfrak{S}'$.

Let $\pmb{\epsilon}$ (in Lemma \ref{decomposition}) be $\pmb{\epsilon}=2/K$ and
\begin{equation*}
  \log K= (\log M_{0})^{C_{7}}.
\end{equation*}
We apply the decomposition Lemma \ref{decomposition} to the semialgebraic set $\mathfrak{S}'$ by
 identifying the algebraic curve $(\xi, \omega_{l_{0}}^{\leq})$ with
 an interval.
 Then we obtain
 \begin{equation*}
   \mathfrak{S}'=\mathfrak{S}_{1}'\cup \mathfrak{S}_{2}'
 \end{equation*}
with
\begin{equation*}
  \textrm{Proj}_{(\xi,\omega_{l_{0}^{\leq}}(\xi))}\mathfrak{S}_{2}'<
  M_{0}^{C} A^{C (\log (m+1))^{C_{5}}} \pmb{\epsilon}
< A^{C (\log (m+1))^{C_{5}}} e^{-(\log (m+1))^{C_{7}}}
< A^{-(\log (m+1))^{C_{4}+6}},
\end{equation*}
since $C_{7}> (C_{4}+10)\vee C_{5}$ in our choice \eqref{C} of constants.

Moreover, for any $|k|>K/2$ and
 the hyperplane $L_{k}= \{(\xi,\omega_{l_{0}}^{\leq}, \langle k,\omega_{l_{0}}^{\leq}\rangle)\}$,
 there is
 \begin{equation*}
 \begin{aligned}
   \textrm{mes}~(\mathfrak{S}'\cap L_{k})<& M_{0}^{C} A^{ C (\log (m+1))^{C_{5}}} \pmb{\epsilon}^{-1}
   e^{-M_{0}^{\beta^{3}}/(2d)}\\
<& A^{C (\log (m+1))^{C_{7}}}  e^{-(m+1)^{C_{0} \beta^{3}}}\\
<& e^{-(m+1)^{\beta^{4} C_{0}} }.
   \end{aligned}
 \end{equation*}
 Therefore, taking into consideration of each $\mathfrak{S}'$,
 there is a set $\mathscr{V}_{2}\subset \mathscr{J}$ satisfying
 \begin{equation*}
   \textrm{mes}~(\mathscr{V}_{2})< C M_{0} \left(A^{-(\log (m+1))^{C_{4}+6}}
   + A^{C(m+1)} e^{-(m+1)^{\beta^{4} C_{0}} }\right)< A^{-(\log (m+1))^{C_{4}+5}},
 \end{equation*}
 such that
\eqref{G M_0 <} holds for all $\sigma= \langle k,\omega_{l_{0}}^{\leq}\rangle$ with
$K/2\leq |k|\leq N$.

\emph{Step three.} We divide $\scr{J}$ into a sequence of disjoint sub-intervals with each
interval of diameter $A^{-(m+1)^{C_{3}}}$, i.e., $\scr{J}=\cup\scr{J'}$ with
$|\scr{J}'|=A^{-(m+1)^{C_{3}}}$. Then for any $\xi\in\scr{J}$ but lying outside  the
boundaries of the subintervals, there is a unique $\scr{J'}$ such that
$\xi\in\scr{J}'$. Suppose $\xi\in \scr{V}_{2}$, i.e., \eqref{G M_0 <} fails for
 $\xi$ and for some $\sigma=\langle k,\omega_{l_{0}}^{\leq}\rangle, K/2\leq |k|\leq N $. We have that, for all $\xi'\in\scr{J}'$, \eqref{G M_0 <} with the above $\sigma$ also
fails but with a smaller constant due to Neumann series (actually, $M_{0}^{C} A^{-(m+1)^{C_{3}}}\ll e^{-M_{0}^{b}}$), i.e., for some $\sigma=\langle k,\omega_{l_{0}}^{\leq}\rangle $, there is
\begin{equation}\label{G M 0 fail}
  \begin{aligned}
& \|G^{\sigma;\leq}_{l_{0}; M_{0}}(\xi')\|>\frac{1}{2} e^{M_{0}^{b}},\\
    \textrm{or}\quad & |G^{\sigma;\leq}_{l_{0};M_{0}}(k,k';\xi')|>\frac{1}{2} \exp(-\alpha'' |k-k'|)\quad \textrm{for some}~
    |k-k'|>M_{0}^{\theta},
  \end{aligned}
\end{equation}
As before, \eqref{G M 0 fail} also has a semialgebraic description. Denoting
by $\scr{V}_{2}'$ the set of all $\xi\in\scr{J}$ such that
\eqref{G M 0 fail} holds for some $\sigma=\langle k,\omega_{l_{0}}^{\leq}\rangle$,
we have $\textrm{mes} (\scr{V}_{2}')< A^{-(\log (m+1))^{C_{4}+5}}$.
Moreover, we have
\begin{equation*}
  \Delta\scr{J}\equiv\bigcup \left\{\xi'\in\scr{J}': \exists~ \xi\in\scr{J}'\subset\scr{J}~\textrm{s.t.}~
  \eqref{G M_0 <}~\textrm{fails for some}~ \sigma=\langle k,\omega_{l_{0}}^{\leq}\rangle\right\}\subset\scr{V}_{2}'.
\end{equation*}
As a result, $\scr{J}\setminus\Delta\scr{J}$ is the union of a sequence of
intervals  with each
interval of diameter $A^{-(m+1)^{C_{3}}}$ and
$\textrm{mes}~(\Delta\scr{J})<A^{-(\log (m+1))^{C_{4}+5}}$.

 Letting $\mathscr{J}$ range over $\Lambda_{l_{0}}$, the total measure of
 the set $\Delta\scr{J}$ removed
 from $\Pi_{l_{0}} $
fulfills
\begin{equation*}
  \textrm{mes}~(\cup_{\scr{J}\in\Lambda_{l_{0}}} \Delta\scr{J})<  A^{l_{0}^{C_{3}}} A^{-(\log (m+1))^{C_{4}+5}}
< \frac{1}{100} A^{-(\log (m+1))^{C_{4}}}
\end{equation*}
in view of $C_{4}>C_{3}$ in \eqref{C}.

Let $\Pi_{+}= \cup_{\scr{J}\in\Lambda_{0}} \Pi\cap (\scr{J}\setminus\Delta\scr{J})$.
Then $\Pi_{+}$ is a collection $\Lambda_{+}^{(1)}$ of disjoint intervals with diameter
$A^{-(m+1)^{C_{3}}}$. Since $\Pi\subset\Pi_{l_{0}}$,  for each interval
$\scr{J}'\in\Lambda_{+}^{(1)}$, there is  a unique $\scr{I}\in\Lambda$ such that
$\scr{J}'\subset\scr{I}$.
On $\Pi_{+}$, \eqref{G M_0 <} essentially holds (up to a constant multiplication by applying the Neumann series).
From \eqref{trunc error}, we see that \eqref{G M_0 k} essentially holds on $\Pi_{+}$.
Since $C_{3}>C_{1}>C_{0}$ in \eqref{C}, it then follows from Lemma \ref{Neumann} that
\eqref{G M_0 k} remains valid on $\scr{O}(\Pi_{+}, A^{-(m+1)^{C_{3}}})$
 by verifying $A^{-(m+1)^{C_{3}}}< A^{-M_{0}^{\theta}}$.

This completes the proof of Lemma \ref{M_0 decay}. \qed

%
%
%
%
%
%
%
%
%
%
%
%
%
%

\subsection{Estimate of $G_{N}$}\label{Green G N}

In this part, we shall establish the Green's function estimate for
$G_{N}$. As we mentioned before, we shall apply a coupling lemma
involving two scales, which is independent of
the KAM iteration.

\begin{lemma}\label{2 scale couple}
 Let the matrix $\mc{T}=\mc{D}+\epsilon \mc{S}$ defined on $[-\mc{N},\mc{N}]^{d}\cap\ZZ^{d}$ satisfy
  \begin{equation*}
    |\mc{S}(x,y)|< e^{-\rho |x-y|},\quad x\neq y.
  \end{equation*}
  Let the integers $0<2 \mc{M}_{0}<\mc{K}<\mc{N}$ be sufficient large and the various constants
  below satisfy
  \begin{align*}
&  C_{1}>C_{0}>10,\quad C_{2}>2 C_{1}+10,\\
& \mc{M}_{0}\sim (\log \mc{N})^{C_{0}},\quad \log\log \mc{K}\sim \log\log \mc{M}_{0},\\
& 0<\frac{\rho_{0}}{2}<\rho<\tilde{\alpha},\quad  0<\frac{\rho_{0}}{2}<\alpha<\tilde{\alpha},\\
& 0<b<\theta<1-\frac{8}{C_{0}}.
  \end{align*}

  Assume
  \begin{itemize}
    \item there is Green's function estimate on $\mc{G}_{\mc{K}}$
    \begin{align*}
& \|\mc{G_{K}}\|\lesssim A^{(\log \mc{K})^{C_{1}}},\\
& |\mc{G_{K}}(x,y)|\lesssim e^{-\alpha |x-y|},\quad \emph{\textrm{for}}~|x-y|> (\log \mc{K})^{C_{2}}.
    \end{align*}
    \item for each $|k_{0}|> \mc{K}/2$, there is
    \begin{align*}
&\|\mc{G}_{k_{0}+[-\mc{M}_{0},\mc{M}_{0}]^{d}}\|\lesssim e^{\mc{M}_{0}^{b}},\\
& |G_{k_{0}+[-\mc{M}_{0},\mc{M}_{0}]^{d}}(x,y)|\lesssim e^{-\tilde{\alpha} |x-y|},
      \quad \emph{\textrm{for}}~ |x-y|> \mc{M}_{0}^{\theta}.
    \end{align*}

  \end{itemize}
  The we have the Green's function estimate
  \begin{align*}
& \|\mc{G_{N}}\|< A^{(\log \mc{N})^{C_{1}}},\\
& |\mc{G_{N}}(x,y)|< e^{-\gamma |x-y|},\quad \emph{\textrm{for}}~|x-y|> (\log \mc{N})^{C_{2}},
  \end{align*}
  where $\gamma> (\alpha\wedge\rho)-(\log \mc{N})^{-8}$.
\end{lemma}

\begin{remark}
  The above lemma also appeared in \cite[Chapter 18]{Bou05_Green}
  and \cite[Lemma 5.1]{BW08}. One should be very careful to establish Lemma
   \ref{2 scale couple} when taking into account the loss of regularity in the
   KAM iteration. In \cite{Bou05_Green} and \cite{BW08}, the off-diagonal
   exponential decay for $G_{M_{0}}(k,k')$ is valid when $|k-k'|>\frac{1}{100} M_{0}$,
   rather than $|k-k'|>M_{0}^{\theta}$ in our imposition.
  We remark that in \cite{Bou05_Green}, this might lead to a great loss of
  regularity at each KAM step, which possibly impedes us to get a uniform analyticity
  domain for the angle variable. In \cite{BW08}, there is no such trouble
  since the matrix therein is of short range.
  This is the main reason why we establish the Green's function estimate
  for those $|k-k'|>M_{0}^{\theta}$.
\end{remark}

\noindent\textbf{Proof.} The proof is based on the application of the
resolvent identity. We divide the proof into two parts which are on
the norm control and the exponential decay estimate, respectively.

\noindent\emph{1. Estimate of the norm.} For any fixed $x\in[-\mc{N},\mc{N}]^{d}$, we define
\begin{equation*}
U(x)=\left\{
\begin{aligned}
& [-\mc{K},\mc{K}]^{d},\quad \textrm{if}~|x|\leq \frac{\mc{K}}{2},\\
& (x+[-\mc{M}_{0}, \mc{M}_{0}]^{d})\cap [-\mc{N},\mc{N}]^{d},\quad \textrm{if}~|x|>\frac{\mc{K}}{2}.
\end{aligned}
\right.
\end{equation*}
For $|x|\leq \mc{K}/2$, we have
\begin{equation*}
\textrm{dist}(x,  [-\mc{N},\mc{N}]^{d}\setminus U(x))\geq \frac{\mc{K}}{2},
\end{equation*}
and  for $|x|> \mc{K}/2$, we have
\begin{equation*}
\textrm{dist}(x, [-\mc{N},\mc{N}]^{d}\setminus U(x))> \mc{M}_{0}.
\end{equation*}
Compute by the resolvent identity
\begin{equation}\label{0}
|\mc{G_{N}}(x,y)|\leq |\mc{G}_{U(x)}(x,y)|~ \chi_{U(x)}(y)+ \sum_{w\in U(x), v\not\in U(x)} |\mc{G}_{U(x)}(x,w)| e^{-\rho |w-v|} |\mc{G_{N}}(v,y)|.
\end{equation}
When $|x|\leq K/2$,
we have
\begin{equation}\label{5}
\begin{aligned}
\mc{G_{N}}(x,y)\leq &
|\mc{G}_{U(x)}(x,y)|~ \chi_{U(x)}(y)
    + \varphi_{\mc{K}}
    e^{-(\rho\wedge \alpha) |x-v| } |\mc{G_{N}}(v,y)|\\
\end{aligned}
\end{equation}
for some $|v|>\mc{K}$,
where
\begin{equation*}
  \varphi_{\mc{K}}=2 \mc{K}^{d} \mc{N}^{d} A^{(\log \mc{K})^{C_{1}}} e^{\rho (\log \mc{K})^{C_{2}}}
<\mc{N}^{2d}.
\end{equation*}
Since $|x-v|>\mc{K}/2$, there is $$\varphi_{\mc{K}} e^{-(\rho\wedge\alpha) |x-y|}
< \mc{N}^{2d} e^{-\rho_{0} \mc{K}/4}<\frac{1}{10}.$$

When $|x|>\mc{K}/2$, we have
\begin{equation}\label{6}
  |\mc{G_{N}}(x,y)|\leq |\mc{G}_{U(x)}(x,y)|~ \chi_{U(x)}(y)+
  \varphi_{\mc{M}_{0}} e^{-\rho |x-v|} |\mc{G_{N}}(v,y)|
\end{equation}
for some $v$ satisfying $|v-x|>\mc{M}_{0}$, where
\begin{equation*}
  \varphi_{\mc{M}_{0}}=2 \mc{M}_{0}^{d} \mc{N}^{d} e^{\mc{M}_{0}^{b}}
    e^{\rho \mc{M}_{0}^{\theta}}< e^{2\rho \mc{M}_{0}^{\theta}}.
\end{equation*}
Moreover, $$\varphi_{\mc{M}_{0}} e^{-\rho |x-v|}< e^{2\rho \mc{M}_{0}^{\theta}} e^{-\rho \mc{M}_{0}}
< \frac{1}{10}.$$

In conclusion, we have
\begin{equation*}
  |\mc{G_{N}}(x,y)|< (A^{(\log \mc{K})^{C_{1}}}+ e^{\mc{M}_{0}^{b}})+\frac{1}{4} \max_{v\in [-\mc{N,N}]^{d}} |\mc{G_{N}}(v,y)|.
\end{equation*}
which further implies
\begin{equation*}
  \max_{x\in [-\mc{N,N}]^{d}} |\mc{G_{N}}(x,y)|< 2  (A^{(\log \mc{K})^{C_{1}}}+ e^{\mc{M}_{0}^{1-}})
\end{equation*}
for any $y\in [-\mc{N,N}]^{d}$.
By Schur's criterion, we finally get
\begin{equation}\label{7}
  \|\mc{G_{N}}\|< 2 \mc{N}^{d}  (A^{(\log \mc{K})^{C_{1}}}+ e^{\mc{M}_{0}^{b}}) < A^{(\log \mc{N})^{C_{1}}}
\end{equation}
by our assumptions on the constants.

\bigskip

\noindent\emph{2. Exponential decay estimate.}

For any $|x|, |y|\leq \mc{N}$, we apply \eqref{5} and \eqref{6}
to take iterations. At each step, there are four cases. See table
\eqref{table}. When $|x-y|>(\log \mc{N})^{C_{2}}$, the iteration would
start from \textbf{A2} or \textbf{A4}.
Note also $10 (\log \mc{N})^{C_{2}}< \mc{K}$ and $\log\mc{K}\sim \log\log \mc{N}$.

\begin{table}\caption{Four cases}\label{table}
  \centering
  \begin{threeparttable}
 \begin{tabular}{|c||c|c|c|c|}
  \hline
 Case & Condition & Estimate of $\mc{G_{N}}(x,y)$ & $|x-v|$\tnote{$\ast$}
& Action\tnote{$\dagger$ } \\
 \hline
  \textbf{A1} & $|x|\leq \mc{K}/2$, $|x-y|\leq (\log \mc{K})^{C_{2}}$
& apply \eqref{7} & N/A\tnote{$\ddagger$ } & off \\
  \hline
  \textbf{A2} & $|x|\leq \mc{K}/2$, $|x-y|> (\log \mc{K})^{C_{2}}$
& apply \eqref{5}  & $> \mc{K}/2$ & on \\
  \hline
  \textbf{A3} & $|x|> \mc{K}/2$, $|x-y|\leq  \mc{M}_{0}^{\theta}$
& apply \eqref{7} & N/A & off \\
  \hline
  \textbf{A4} & $|x|> \mc{K}/2, |x-y|> \mc{M}_{0}^{\theta}$ & apply \eqref{6} & $>\mc{M}_{0}$ & on \\
  \hline
\end{tabular}

 \begin{tablenotes}
        \footnotesize
        \item[$\ast$]  $v$ originates from the application of \eqref{5} or \eqref{6}.
        \item[$\dagger$]  This indicates the iteration is going on or called off.
        \item [$\ddagger$] N/A indicates not applied since the iteration terminates.
      \end{tablenotes}
\end{threeparttable}

\end{table}

 A sequence of iterations should obey the following rule
 \begin{align*}
&\cdots \rightarrow (A2 \rightarrow A4)  \rightarrow A4\rightarrow\cdots, \\
 \textrm{or}\quad &\cdots \rightarrow (A2\rightarrow A4)\rightarrow
 (A2\rightarrow A4)\rightarrow \cdots,\\
 \textrm{or}\quad &\cdots\rightarrow  A4\rightarrow  (A2\rightarrow A4)
 \rightarrow\cdots, \\
 \textrm{or}\quad &\cdots\rightarrow  A4\rightarrow A4 \rightarrow \cdots.
 \end{align*}
The iteration would stop in the following way
\begin{align*}
& \cdots A4\rightarrow A1~/~ A3,\\
\textrm{or}\quad&\cdots A2\rightarrow A3.
\end{align*}

Assume  we are able to iterate $\mathbf{(A2\rightarrow A4)}$ for $p$ times
and iterate $\mathbf{A4}$ alone for
$q$ times. Then we have
\begin{equation}\label{8}
  |\mc{G_{N}}(x,y)|< (p+q-1) \varphi_{\mc{K}}^{p}\varphi_{\mc{M}_{0}}^{p+q-1}
   e^{-(\alpha\wedge\rho) |x-y|}
  + \varphi_{\mc{K}}^{p}\varphi_{\mc{M}_{0}}^{p+q} e^{-(\alpha\wedge\rho) |x-v_{2p+q}|}
  ~|\mc{G_{N}}(v_{2p+q},y)|
\end{equation}
and
\begin{equation*}
  |x-v_{2p+q}|> p\frac{\mc{K}}{2}+q \mc{M}_{0}.
\end{equation*}

Let
\begin{equation*}
  p^{*}\frac{\mc{K}}{2}+q^{*} \mc{M}_{0}= 10 |x-y|
\end{equation*}
and thus $p^{*}\leq 20 |x-y|/\mc{K}$ and
$p^{*}+q^{*}\leq 10 |x-y|/\mc{M}_{0}$.
Moreover, we have
\begin{equation}\label{9}
  \log \varphi_{\mc{K}}^{p^{*}}\leq \log \mc{N}^{2 d p^{*}}\leq \frac{40 d \log \mc{N}}{\mc{K}}|x-y|<
   \frac{1}{10 (\log \mc{N})^{8}} |x-y|
\end{equation}
and also
\begin{equation*}\label{10}
  \log \varphi_{\mc{M}_{0}}^{p^{*}+q^{*}}\leq \frac{20\rho}{\mc{M}_{0}^{1-\theta}} |x-y|
< \frac{1}{10 (\log \mc{N})^{8}} |x-y|.
\end{equation*}

If $(p,q)=(p^{*},q^{*})$, it follows from \eqref{8} that
\begin{equation*}
\begin{aligned}
  |\mc{G_{N}}(x,y)|<&\frac{1}{2} \exp\left(-
  (\alpha\wedge\rho-\frac{1}{5(\log \mc{N})^{8}})|x-y|\right)\\
&+ \exp\left(
  \frac{1}{5(\log \mc{N})^{8}})|x-y|\right)
  A^{(\log \mc{N})^{C_{1}}} e^{-10 (\alpha\wedge\rho) |x-y|}\\
<& \exp\left(-(\alpha\wedge\rho
  -\frac{1}{(\log \mc{N})^{8}})|x-y|\right)
\end{aligned}
\end{equation*}
since $|x-y|>(\log \mc{N})^{C_{2}}$ and $C_{2}>2C_{1}+10$.

If we stop the iteration before $(p,q)$ arriving at $(p^{*},q^{*})$,
then we have
\begin{equation*}
\begin{aligned}
   \varphi_{\mc{K}}^{p}\varphi_{\mc{M}_{0}}^{p+q} e^{-(\alpha\wedge\rho) |x-v_{2p+q}|}
  ~|\mc{G_{N}}(v_{2p+q},y)|<& \exp\left(\frac{1}{5(\log \mc{N})^{8}})|x-y|\right)
  A^{(\log \mc{N})^{C_{1}}} e^{-(\alpha\wedge\rho) |x-y|} e^{\rho \mc{M}_{0}^{\theta}}\\
<&\frac{1}{2} \exp\left(-
  (\alpha\wedge\rho-\frac{1}{(\log \mc{N})^{8}})|x-y|\right),
  \end{aligned}
\end{equation*}
which together with \eqref{8} implies
\begin{equation*}
  |\mc{G_{N}}(x,y)|< \exp\left(-(\alpha\wedge\rho
  -\frac{1}{(\log \mc{N})^{8}})|x-y|\right).
\end{equation*}
This completes the proof. \qed
\bigskip

Now we turn to establish the Green's function estimates on $G_{N}$. Recall the two scales
 $0<M_{0}<K$ satisfying
\begin{equation*}
M_{0}= (\log N)^{C_{0}},\quad \log K=(\log M_{0})^{C_{7}}.
\end{equation*}
Moreover, we take $l_{0}$ and $l_{1}$ such that
\begin{equation*}
  l_{0}=C_{8} \log M_{0},\quad K=A^{l_{1}},
\end{equation*}
with $C_{8}> (1+\log 10)/(\log \frac{4}{3})$.
By the Iterative Lemma \ref{iter lemma}, we have
\begin{equation*}
\begin{aligned}
& \|G_{l_{1}; K}\|<  A^{l_{1}^{C_{1}}}< A^{(\log K)^{C_{1}}},\\
& |G_{l_{1};K}(k,k')|<  e^{-s_{l_{1}} |k-k'|}\quad \textrm{for}~
 \quad |k-k'|>l_{1}^{C_{2}}\sim (\log K)^{C_{2}},
  \end{aligned}
\end{equation*}
for any $\xi\in\Pi_{l_{1}}$.
Using \eqref{T regularity}, Lemma \ref{variation} and Lemma \ref{Neumann},
we have
\begin{equation}\label{T K decay}
  \begin{aligned}
& \|G_{ K}\|<   A^{(\log K)^{C_{1}}},\\
& |G_{K}(k,k')|<  e^{-s |k-k'|},\quad \textrm{for}~
 \quad |k-k'|> (\log K)^{C_{2}},
  \end{aligned}
\end{equation}
since
\begin{equation*}
\epsilon_{l_{1}}^{1/10}\cdot N+ \epsilon_{l_{1}}^{1/10} <\epsilon_{l_{1}}^{1/20}< A^{-l_{1}^{C_{1}}}.
\end{equation*}
(Indeed, $\log \log \epsilon_{l_{1}}^{-1}\sim l_{1}\sim (\log (m+1))^{C_{7}}\gg \log (m+1)\sim \log\log N$.)
Moreover, by verifying $A^{-(m+1)^{C_{3}}}<A^{-(\log K)^{C_{1}}}$, it follows
again from Lemma \ref{Neumann} that
\eqref{T K decay} remains valid on $\scr{O}(\Pi_{l_{1}}, A^{-(m+1)^{C_{3}}})$
 and hence on $\scr{O}(\Pi_{+}, A^{-(m+1)^{C_{3}}})$.

Then, using \eqref{T K decay} on $\scr{O}(\Pi_{+}, A^{-(m+1)^{C_{3}}})$ and
Lemma \ref{M_0 decay}, we obtain from Lemma \ref{2 scale couple} that
\begin{equation}\label{G N decay}
  \begin{aligned}
&\|G_{N}\|< A^{(\log N)^{C_{1}}},\\
& |G_{N}(k,k')|< e^{-(s- (\log N)^{-8})~ |k-k'| },\quad \textrm{for}~|k-k'|> (\log N)^{C_{2}}
  \end{aligned}
\end{equation}
holds on $\scr{O}(\Pi_{+}, A^{-(m+1)^{C_{3}}})$.
Note that $s-(\log N)^{-8}> s^{1}>s_{+}=s_{m+1}$.
This completes the proof of Lemma \ref{Green}.

%

\section {Appendix A: Large deviation  theorem}

The appendix is devoted to the proof of the large
deviation theorem (Lemma \ref{LDT}), which can be read independently. The proof follows exactly the same line
in  \cite{BGS02}
and we prove it here for completeness.
It is worthy noticing that the notations below are also independent of
the main body of the paper. For that reason, we write simply $\mc{T}$ by  $T$ and so on.

\subsection{Notations and phrases}
We consider matrix defined on $\ZZ^{d}$. For ${m}=(m_{1},\cdots, m_{d}), {n}=(n_{1},\cdots, n_{d})\in\ZZ^{d}$, we define
the distance by
\begin{equation*}
  |{m}-{n}|=\max_{1\leq j\leq d} |m_{j}-n_{j}|.
\end{equation*}
For $\Lambda\subset\ZZ^{d}$, we denote the diameter of $\Lambda$ by $|\Lambda|$. For a
matrix $A$ defined on $\ZZ^{d}$, we denote by $\|A\|$ the
the operator norm induced by
the $\ell^{2}$ norm  of a vector in $\ZZ^{d}$. The inverse of a matrix is always
denoted by $G$.

When applying in resolvent identity, we shall control the Green's
function $G_{\Lambda}$ with $\Lambda$ being the difference of two
boxes in $\ZZ^{d}$. For that reason, as in \cite{BGS02},
we introduce the elementary regions.
An elementary region is defined to be a set $\Lambda$ of the form
\begin{equation*}
  \Lambda=R\setminus (R+z)
\end{equation*}
where $z\in\ZZ^{d}$  is arbitrary and $R$ is a block in $\ZZ^{d}$, i.e.,
\begin{equation*}
  R= \{y=(y_{1},\cdots, y_{d})\in\ZZ^{d}: y_{i}\in [x_{i}-a_{i}, x_{i}+a_{i}], i=1\cdots, d\}.
\end{equation*}
The size of an elementary region $\Lambda$ is simply its diameter. For any integer $M>0$, the set of all
elementary regions of size $M>0$ is denoted by $\mathcal{ER}(M)$ and
are also referred as $M$-regions.
The class of elementary regions consists of $d$-dimensional rectangles, L-shaped regions and
$(d-1)$ -dimensional rectangles with normal vector parallel to the
axis.

Note that
these regions play only a role in the application of resolvent identity
in the presence of interior corners, but basically have no effect
on the other parts of the argument.

Given a elementary region $\Lambda$, we consider exhaustion
$\{S_{j}(m)\}_{j=0}^{l}$ of $\Lambda$ of width $2M$ centered at $m\in\Lambda$ defined
inductively by
\begin{equation}\label{exhaustion}
  \begin{aligned}
&S_{0}(m)= Q_{M}(m)\cap \Lambda, \quad Q_{M}(m)=\{n\in\ZZ^{d}: |n-m|\leq M\},\\
&S_{j}(m)=\bigcup_{n\in S_{j-1}(m)} (Q_{2M}(n)\cap\Lambda),\quad \textrm{for}\quad
    1\leq j\leq l,
  \end{aligned}
\end{equation}
where $l$ is maximal such that $S_{l}\neq \Lambda$.
Define the annulus between the exhaustion by
\begin{equation}\label{annulus}
   A_{j}(m)= S_{j}(m)\setminus S_{j-1}(m),\quad 0\leq j\leq l
\end{equation}
with $ S_{-1}(m) =\emptyset.$
We have the following two simple observations:

\begin{itemize}
\item Except the possible exception of a single
annulus, $Q_{M}(n)\cap A_{j}(m)$ is an elementary region for
all $n\in A_{j}(m)$. The exceptional annulus is the one that
contains the unique interior corner of $\Lambda$ (i.e., the
corner lying in the interior of the hull of $\Lambda$).

\item Any two cubes $Q_{M}(n_{1})$ and $Q_{M}(n_{2})$
with centers $n_{1}$ and $n_{2}$ lying in nonadjacent annuli
are disjoint.
\end{itemize}

\subsection{Coupling Lemma for long range operators}

We present and prove two kinds of coupling lemmas.

\begin{lemma}\label{CL1}

 Let $T$ be a  matrix defined on  a finite
set  $\Lambda\subset\ZZ^{d}$, $|\Lambda|=N$.
Let the various constants below satisfy
\begin{align*}
& 0<\theta<1,\quad 0<b<1,\quad 0<\tau<1, \quad b\tau<\theta,\quad \alpha>0,\quad \rho>0.
\end{align*}
Further let
\begin{equation*}
  (\log N)^{10/b}< M< N^{\tau}.
\end{equation*}

Assume the following properties hold.
\begin{enumerate}[(i)]

\item The matrix $T$ exhibits the off diagonal exponential decay
\begin{equation}\label{CL1-1}
  |T(x,y)|<  e^{-\rho|x-y|},\quad x\neq y,\quad x,y\in\Lambda.
\end{equation}

\item For every $m\in\Lambda$, there is a subinterval $U(m)\subset \Lambda$
containing $m$ with
\begin{equation}\label{CL1-2}
|U(m)|=M\quad  \textrm{and} \quad\emph{dist}~(m,\Lambda\setminus U(m))>\frac{M}{2}
\end{equation}
 such that
\begin{equation}\label{CL1-3}
  \|G_{U(m)}\|< e^{M^{b}}
\end{equation}
and
\begin{equation}\label{CL1-4}
  |G_{U(m)}(x,y)|< e^{-\alpha |x-y|}\quad  \emph{for }~|x-y|>N^{\theta},\quad x,y\in U(m).
\end{equation}

\end{enumerate}

Then, there is
\begin{equation*}
\begin{aligned}
& \|G_{\Lambda}\|< 2N^{d} e^{M^{b}}< e^{N^{b}},\\
& |G_{\Lambda}(x,y)|< e^{-\alpha' |x-y|}\quad \emph{for}~|x-y|> N^{\theta}.
  \end{aligned}
\end{equation*}
provided $N$ is large enough, i.e., $N\geq \underline{N}(\alpha, b, d,\rho, \theta, \tau)$. Moreover, the decay rate $\alpha'\geq (\alpha\wedge\rho)-(\log N)^{- 50}$.
\end{lemma}

The proof of Lemma \ref{CL1} is similar to that of Lemma \ref{2 scale couple} and is much simpler.
We omit it here.

\bigskip

\begin{lemma}\label{CL2}
Let $T$ be a matrix defined on $\Lambda_{0}\in\mathcal{ER}(N)\subset\ZZ^{d}$. Let the various constants below satisfy
\begin{equation}\label{103}
   0<\tau\leq b\leq \theta<1,\quad
   \theta\geq \frac{1-2\tau}{1-\tau},
\end{equation}
and let
\begin{equation*}
  N^{\tau}<M_{0}<2 N^{\tau}.
\end{equation*}
  Assume the following properties hold.
    \begin{enumerate}[(i)]
  \item  The matrix $T$ exhibits the off-diagonal exponential decay
  \begin{equation}\label{CL2-3}
    |T(m,n)|< e^{-\rho |m-n|},\quad m,n\in\Lambda_{0}, m\neq n.
  \end{equation}

 \item   For any $\Lambda\in\mathcal{ER}(L)$, $\Lambda\subset\Lambda_{0}$ with
  any $N^{\tau}< L< N$, there is a bounded inverse
  \begin{equation}\label{CL2-1}
    \|G_{\Lambda}\|< e^{L^{b}}.
  \end{equation}
  We say an elementary region
  $\Lambda\in \mathcal{ER}(L), \Lambda\subset\Lambda_{0}$ is
    good  if in addition to \eqref{CL2-1} the Green's function exhibits the
 off diagonal decay
\begin{equation}\label{CL2-4}
  |G_{\Lambda}(m,n)|< e^{-\alpha(L) |m-n|},\quad m,n\in\Lambda,|m-n|>L^{\theta}.
\end{equation}
Otherwise $\Lambda$ is called bad.

  \item  For any family $\mathcal{F}$ of pairwise disjoint
  bad $M'$-regions in $\Lambda_{0}$ with $M_{0}+1\leq M' \leq 2M_{0}+1$,
  \begin{equation}\label{CL2-5}
    \#\mathcal{F}< \frac{N^{b}}{M_{0}}.
  \end{equation}
  \end{enumerate}

Then, there is
  \begin{equation*}
    |G_{\Lambda_{0}}(m,n)|< e^{-\alpha' |m-n|},\quad\textrm{for all}\quad
    m,n\in\Lambda_{0}, |m-n|>  N^{\theta}.
  \end{equation*}
  provided $N$ is sufficiently large, i.e.,
  $N\geq \underline{N}(b,d, \tau,\theta)$.
Moreover, $\alpha'\geq (\alpha\wedge \rho)-N^{-\delta}$ for some $\delta=\delta(b,d, \tau,\theta)>0$ and $\alpha=\alpha(M_{0})$.
\end{lemma}

\noindent\textbf{Proof.}
The proof is based on an iteration procedure. In the first step,
we give a detailed analysis on the off diagonal decay of the
Green's function at small scale $M_{1}$. Then we list an induction
statement, whose proof is basically same to that in the first step
and hence is omitted. By the finitely many iterations, we
obtain the Green's function estimate at large scale $N$.

\bigskip

\noindent\underline{\textbf{The first step. }}
Let $M_{1}=[M_{0}^{\lambda}]$ with $\lambda>1$ and consider $\Lambda_{1}\in \mathcal{ER}(M_{1}), \Lambda_{1}\subset\Lambda_{0}$.
Fix any $m\in\Lambda_{1}$ and let $\{S_{j}(m)\}_{j=0}^{l}$ be the
exhaustion of $\Lambda_{1}$ of width $2M_{0}$ and centered at $m$
(see \eqref{exhaustion} and the associated annuli in \eqref{annulus}).

We say an annulus $A_{j}(m)$ is good if for any $n\in A_{j}(m)$ both
$Q_{M_{0}}(n)\cap A_{j}(m)$ and $Q_{M_{0}}(n)\cap \Lambda_{1}$ are good regions
in the sense of \eqref{CL2-1} and \eqref{CL2-4}.
Otherwise $A_{j}(m)$ is bad. Note that there is at most one annulus
$A_{j_{0}}$ (consisting the interior corner of $\Lambda_{1}$)
such that $Q_{M_{0}}(n)\cap A_{j_{0}}(n)$ possibly fails to be an elementary
region. In this case , $A_{j_{0}}$ is counted among the bad annuli.
Note also that for good annuli, the diameter of those $Q_{M_{0}}(n)\cap A_{j}(m)$
ranges from $M_{0}+1$ to $2M_{0}+1$.

With the above definition of good and bad annuli, we say that
an elementary region $\Lambda_{1}\in\mathcal{ER}(M_{1}), \Lambda_{1}\subset \Lambda_{0}$ is "GOOD"\footnote{We use "GOOD" here to
make a difference from the goodness of
an elementary region as in \eqref{CL2-1} and \eqref{CL2-4}} if, for any $m\in\Lambda_{1}$, there are at most
$B_{1}=\kappa\frac{M_{1}^{\theta}}{M_{0}}$ many bad annuli for the associated exhaustion centered at
$m$, where $\kappa$ will be determined below. Otherwise, the
$M_{1}$-region $\Lambda_{1}$ is called "BAD".

Let $\mathcal{F}_{1}$ be an arbitrary family of pairwise disjoint
"BAD" $M_{1}$-regions contained in $\Lambda_{0}$. If $\Lambda_{1}\in\mathcal{F}_{1}$,
we can find an exhaustion of $\Lambda_{1}$ centered at some $m\in\Lambda_{1}$ such that
there are at least $\frac{1}{2} B_{1}$ many nonadjacent annuli.
Each bad annuli $A_{j}(m)$ contains a bad $M'$-region
($Q_{M_{0}}\cap A_{j}(m)$ or $Q_{M_{0}}\cap \Lambda_{1}$) with
 $M_{0}\leq M'\leq 2M_{0}+1$, which does not intersect
with that in the nonadjacent bad annulus. As a result, we have
\begin{equation*}
  \#\mathcal{F}_{1}<\frac{\#\mathcal{F}}{B_{1}/2}< \frac{2 N^{b}}{\kappa M_{1}^{\theta}}.
\end{equation*}

Consider a "GOOD" region $\Lambda_{1}\in \mathcal{ER}(M_{1})$ and fix any pair
$m,n\in \Lambda_{1}, |m-n|> M_{1}^{\theta}$. Let $\{S_{j}(m)\}_{j=0}^{l}$ and $\{A_{j}(m)\}_{0}^{l}$
be the associated exhaustion and annuli of $\Lambda_{1}$ of width $2M_{0}$
centered at $m\in\Lambda_{1}$. Let $A_{j}, A_{j+1},\cdots, A_{j+s}$ be
adjacent good annuli and denote
\begin{equation*}
  U=\bigcup_{i=j}^{j+s} A_{i}.
\end{equation*}
Obviously, $|U|\geq 2M_{0}(s+1)$.
We claim the following Green's function estimate on $G_{U}$
\begin{equation}\label{CL2-9}
|G_{U}(x,y)|< e^{ \beta (2M_{0}-|x-y|)}, \quad\textrm{for all}\quad x,y\in U,
\end{equation}
where $$\beta=\alpha\wedge\rho=\alpha(M_{0})\wedge \rho.$$

Usually $U$ is no longer an elementary region and thus \eqref{CL2-1}
is not applicable to get a norm estimate on $G_{U}$,
Nevertheless, we can invoke Lemma \ref{CL1} to estimate $\|G_{U}\|$.
For any $n\in A_{i}\subset U$, by definition both $Q_{M_{0}}(n)\cap \Lambda_{1}$
and $Q_{M_{0}}(n)\cap A_{i}$ are good. Following the notations in Lemma \ref{CL1}, we take $U(n)=Q_{M_{0}}(n)\cap \Lambda_{1}$ when $Q_{M_{0}}(n)\subset U$ and
take $U(n)= Q_{M_{0}}(n)\cap A_{j}$ when $Q_{M_{0}}(n)\setminus U\neq \emptyset$.
Then we have
\begin{equation}\label{CL2-8}
  \|G_{U}\|< 2 M_{1}^{d} e^{(2M_{0}+1)^{b}}.
\end{equation}
Next we repeat the same analysis as Lemma \ref{CL1} and  obtain
\begin{equation}\label{CL2-6'}
  |G_{U}(x,y)|< e^{\beta M_{0}} e^{-\beta |x-y|}, \quad\textrm{for}~|x-y|>M_{0}
\end{equation}
as long as
\begin{equation*}
  1<\lambda<2-(b\vee \theta).
\end{equation*}
Then the claim \eqref{CL2-9} is an immediate result of \eqref{CL2-8} \eqref{CL2-6'}.

Now we are back to establish the off diagonal estimate for a good
$M_{1}$-region $\Lambda_{1}$, i.e., to establish
\begin{equation*}
  |G_{\Lambda_{1}}(m,n)|< e^{-(\beta-)|m-n|}\quad\textrm{for}\quad
  |m-n|> M_{1}^{\theta}, m,n\in \Lambda_{1}.
\end{equation*}

Recall the exhaustion of $\Lambda_{1}$ of width $2 M_{0}$ centered at
$m$. Suppose $S_{0}(m)$ is good and write an exhaustion
\begin{equation*}
   S_{0}(m)\subset J_{0}\subset J_{1}\subset\cdots\subset J_{g}=\Lambda_{1}
\end{equation*}
satisfying
\begin{itemize}
  \item $J_{s+1}\setminus J_{s}$ is the union of adjacent bad  annuli (resp. union of adjacent good annuli) if $s$ is even (resp. if $s$ is odd);
  \item The exhaustion is maximal in the sense that if
    \begin{equation*}
      J_{s+1}\setminus J_{s}=\bigcup_{j=j_{s}}^{j_{s+1}} A_{j} ,\quad s\in 2\NN,
    \end{equation*}
    then $A_{j_{s}-1}$, $A_{j_{s+1}+1}$ are  good and $A_{j}$ is bad for all
    $j_{s}\leq j\leq j_{s+1}$. The case of $s$ being odd is similar;
  \item $J_{s}$ is the elementary region in $\Lambda_{1}$ for all $1\leq s\leq g$.
\end{itemize}
By the "GOOD" property of $\Lambda_{1}$, there is
\begin{equation}\label{CL2-14}
\sum_{s~\textrm{even}} (j_{s+1}-j_{s})< \kappa \frac{M_{1}^{\theta}}{M_{0}},
\end{equation}
and hence
\begin{equation}\label{CL2-15}
  g< 2\kappa \frac {M_{1}^{\theta}}{M_{0}}.
\end{equation}
To begin with, we see from \eqref{CL2-9} that
\begin{equation*}
  |G_{J_{0}}(m,y)|< e^{\beta(2 M_{0}-|m-y|)}\quad\textrm{for}\quad y\in J_{0}.
\end{equation*}
Take
\begin{equation}\label{CL2-16}
\varphi_{0}= e^{2 \beta M_{0}}
\end{equation}
 and we assume by induction that
\begin{equation}\label{CL2-18}
  |G_{J_{s}}(m,y)|< \varphi_{s}~ e^{-\beta|m-y|}\quad \textrm{for}\quad y\in J_{s}.
\end{equation}

If $s+1$ is odd, $J_{s+1}\setminus J_{s}$ is made up of  bad annuli. For any $y\in J_{s+1}$,
we apply the resolvent identity
\begin{equation}\label{CL2-10}
  \begin{aligned}
    |G_{J_{s+1}}(m,y)|<& |G_{J_{s}}(m,y)|\chi_{J_{s}}(y)+\sum_{z\in J_{s}, z'\in J_{s+1}\setminus J_{s}} |G_{J_{s}}(m,z)| e^{-\rho |z-z'|} |G_{J_{s+1}}(z',y)|\\
<& \varphi_{s} e^{-\beta|m-y|} + \varphi_{s}
    \sum_{z\in J_{s}, z'\in J_{s+1}\setminus J_{s}}
    e^{-\beta |m-z|-\rho |z-z'|} |G_{J_{s+1}}(z',y)|\\
<& \varphi_{s} e^{-\beta |m-y|} + \varphi_{s} M_{1}^{2d} e^{M_{1}^{b}}
    \max_{z'\in J_{s+1}\setminus J_{s}} e^{-\beta|m-z'|}.
  \end{aligned}
\end{equation}
If $y\in J_{s}$, then $|m-z'|\geq |m-y|$ and if $y\in J_{s+1}\setminus J_{s}$
there is
\begin{equation*}
|m-z'|\geq \textrm{dist}(m, \partial S_{j_{s}-1})
\geq |m-y|-2M_{0} (j_{s+1}-j_{s}).
\end{equation*}
Consequently, we have
\begin{equation}\label{CL2-12}
|G_{J_{s+1}}(m,y)|< (1+M_{1}^{2d} e^{M_{1}^{b}}
e^{2 \beta M_{0} (j_{s+1}-j_{s})}) ~\varphi_{s} e^{-\beta |m-y|}
\end{equation}

If $s+1$ is even, $J_{s+1}\setminus J_{s}$ is made up of  good annuli. For any $y\in J_{s}$,
we repeat the  resolvent identity analysis in \eqref{CL2-10} and obtain from
$|m-z'|\geq |m-y|$ that
\begin{equation}\label{CL2-11}
|G_{J_{s+1}}(m,y)|< \varphi_{s} e^{-\beta |m-y|} (1+  M_{1}^{2d} e^{M_{1}^{b}}).
\end{equation}
If $y\in J_{s+1}\setminus J_{s}$, we have
\begin{equation*}
|G_{J_{s+1}}(m,y)|< \sum_{z'\in J_{s}, z\in J_{s+1}\setminus J_{s}}
|G_{J_{s+1}}(m,z')| ~e^{-\rho |z'-z|} ~|G_{J_{s+1}\setminus J_{s}} (z, y)|.
\end{equation*}
Applying \eqref{CL2-9} with $U=J_{s+1} \setminus J_{s}$ to $G_{J_{s+1}\setminus J_{s}}(z,y)$
and applying \eqref{CL2-11} to $G_{J_{s+1}}(m,z')$
we have
\begin{equation}\label{CL2-13}
|G_{J_{s+1}}(m,y)|< M_{1}^{2d} e^{M_{1}^{b}} \varphi_{s} (1+  M_{1}^{2d} e^{M_{1}^{b}}) e^{2 \beta M_{0}}
e^{-\beta|m-y|}.
\end{equation}

In conclusion, we can take
\begin{equation}
\varphi_{s+1} = \left\{
\begin{aligned}
& e^{3 \beta M_{0} (j_{s+1}-j_{s})}\varphi_{s}, \quad s~\textrm{is even};\\
& e^{3 \beta M_{0} } \varphi_{s}, \quad s~\textrm{is odd}.
\end{aligned}
\right.
\end{equation}
By \eqref{CL2-14},\eqref{CL2-15} and \eqref{CL2-16}, we get
\begin{equation}\label{CL2-17}
  \varphi_{g}< e^{6 \beta M_{0} g} < e^{15 \beta \kappa M_{1}^{\theta} }.
\end{equation}
Suppose $S_{0}(m)$ is bad, then $\varphi_{0}< e^{M_{1}^{b}} e^{\beta \kappa M_{1}^{\theta}}$
and \eqref{CL2-17} is also valid by the same analysis.
Therefore, we prove the induction statement \eqref{CL2-18} and get
\begin{equation}
|G_{\Lambda_{1}}(m,n)|< e^{15\kappa \beta M_{1}^{\theta} } e^{-\beta |m-n|}.
\end{equation}
Since $|m-n|>M_{1}^{\theta}$, it follow that
\begin{equation*}
  |G_{\Lambda_{1}}(m,n)|< e^{-\alpha_{1} |m-n|},\quad \alpha_{1}=\beta (1-15 \kappa),
\end{equation*}
which establishes the off diagonal decay of $G_{\Lambda_{1}}$
for a "GOOD" elementary region $\Lambda_{1}\in\mathcal{ER}(M_{1})$.

\bigskip

\noindent\underline{\textbf{Induction statement.}}~
Let $\kappa< 10^{-2}$ be specified later and let $\lambda$ satisfy
\begin{equation*}
1<\lambda <2-(b\vee \theta)=2-\theta,\quad b\lambda<1.
\end{equation*}
Indeed, due to our choice of $b\leq \theta$, we have $2-\theta< 1/b$.

Define inductively
$M_{t}= [M_{t-1}^{\lambda}], t\leq t_{*}$ and $t_{*}$ is specified later.
Consider $\Lambda_{t}\in \mathcal{ER}(M_{t}), \Lambda_{t}\subset\Lambda_{0}$.
Fix any $m\in\Lambda_{t}$ and let $\{S_{j}(m)\}_{j=0}^{l}$ be the
exhaustion of $\Lambda_{t}$ of width $2M_{t-1}$ and centered at $m$. Let
$\{A_{j}\}_{j=0}^{l}$ be the associated annuli.

We say an annulus $A_{j}(m)$ is good if for any $n\in A_{j}(m)$ both
$Q_{M_{t-1}}(n)\cap A_{j}(m)$ and $Q_{M_{t-1}}(n)\cap \Lambda_{t}$ are good regions
in the sense of \eqref{CL2-1} and \eqref{CL2-4} but with the  decay rate
$\alpha_{t-1}= \beta_{t-1} (1-15\kappa)=\beta (1-15\kappa)^{t-1}$ and $\beta_{t-1}=\alpha_{t-2}\wedge \rho$.
Otherwise $A_{j}(m)$ is bad.
We say that
an elementary region $\Lambda_{t}\in\mathcal{ER}(M_{t}), \Lambda_{t}\subset \Lambda_{0}$ is "GOOD" if, for any $m\in\Lambda_{t}$, there are at most
$B_{t}=\kappa\frac{M_{t}^{\theta}}{M_{t-1}}$ many bad annuli for the associated exhaustion centered at
$m$. Otherwise, the
$M_{t}$-region $\Lambda_{t}$ is called "BAD".
 Let $\mathcal{F}_{t-1}$ be the family of pairwise disjoint "BAD" $M_{t-1}$-regions
 contained in $\Lambda_{0}$.

 Assume
 \begin{equation*}
   \# \mathcal{F}_{t-1}< \left(\frac{2}{\kappa}\right)^{t-1} \frac{N^{b}}{M_{t-1}^{\theta}} (M_{t-2}M_{t-3}\cdots M_{1})^{1-\theta}
 \end{equation*}
 and \eqref{CL2-1} holds for all $N^{\tau}< L\leq N$.
 Then for any "GOOD" $M_{t}$-region $\Lambda_{t}\in\mathcal{ER}(M_{t})$,
 the Green's function $G_{\Lambda_{t}}$ exhibits off diagonal decay
 \begin{equation*}
   |G_{\Lambda_{t}}(m,n)|< e^{-\alpha_{t} |m-n|},\quad |m-n|>M_{t}^{\theta},
   m,n\in\Lambda_{t}
 \end{equation*}
 with $\alpha_{t}=\beta (1-15\kappa)^{t}$. Moreover, denoting by $\mathcal{F}_{t}$
 the family of pairwise disjoint "BAD" $M_{t}$-regions contained in $\Lambda_{0}$,
 there is
 \begin{equation*}
   \#\mathcal{F}_{t}< \left(\frac{2}{\kappa}\right)^{t} \frac{N^{b}}{M_{t}^{\theta}} (M_{t-1} M_{t-2} \cdots M_{1} )^{1-\theta}.
 \end{equation*}

 The proof of the above statement is the same to that in the first step and
 is omitted.
 \bigskip

 \noindent\underline{\textbf{Off diagonal estimate of $G_{N}$}.}
In order to reach size $N=M_{t_{*}}$, the number $t_{*}$ of steps
should satisfy
\begin{equation*}
  \lambda^{t_{*}} \log M = \log N
\end{equation*}
hence $\lambda^{t_{*}}\sim \frac{1}{\tau}$.
It then suffices to show that $[-N,N]^{d}$ is a "GOOD" $M_{t_{*}}$-region,
which is of course valid if
\begin{equation}\label{101}
  \left(\frac{2}{\kappa}\right)^{t_{*}-1} \frac{N^{b}}{M_{t_{*}-1}^{\theta}}
  (M_{t_{*}-2}\cdots M_{1})^{1-\theta}
< \kappa \frac{N^{\theta}}{M_{t_{*}-1}}.
\end{equation}
Obviously, \eqref{101} is equivalent to
\begin{equation*}
  \kappa> 2 \left(\frac{1}{2 N^{\gamma}}\right)^{1/t_{*}},
  \quad \gamma=\theta-b-
  \frac{\lambda\tau}{\lambda-1} (\lambda^{t_{*}-1}-1)(1-\theta).
\end{equation*}
To keep $\gamma>0$, it suffices to take
\begin{equation*}
  \lambda>\frac{1-b}{\theta(1-\tau)+\tau-b}
\end{equation*}
which is compatible with $\lambda< 2-\theta$
according to our choice of $\theta\geq (1-2\tau)/(1-\tau)$.

Take
\begin{equation*}
  \kappa=\kappa_{N}= 4 N^{-\frac{\gamma \log \lambda}{\log \tau^{-1}}}
\end{equation*}
and then
\begin{equation*}
  |G_{N}(m,n)|< e^{-\alpha_{t_{*}} |m-n|},\quad |m-n|> N^{\theta}.
\end{equation*}
The conclusion is valid with some choice of $\delta=\delta(b,d, \tau,\theta,\lambda(b,\tau,\theta))$ such
that
\begin{equation*}
 \alpha'=\alpha_{t_{*}}= \beta (1-15 \kappa)^{t_{*}} >(\alpha(M_{0})\wedge\rho)-N^{-\delta}.
\end{equation*}
This completes the proof.
\qed

\subsection{Matrix-valued Cartan's theorem}

The following matrix-valued Cartan's theorem as well as its proof is given in
\cite{Bou05_Green}.

\begin{lemma}\label{Cartan1}
  Let $A(\sigma)$ be a matrix valued function defined on $\sigma\in[-\delta,\delta]$
  with $A(\sigma)(m,n)\in\CC$ for
  $m,n \in\Lambda_{0}\subset\ZZ^{d}$, $|\Lambda_{0}|=N$. Assume
  \begin{enumerate}[(i)]
    \item $A(\sigma)$ is real analytic in $\sigma$, and there is a holomorphic
    extension to a strip
    \begin{equation}\label{Cartan1-1}
      |\textbf{Re}~z |< \delta, \quad |\textbf{Im} z|<\gamma
    \end{equation}
    satisfying
    \begin{equation}\label{Cartan1-2}
      \|A(z)\|<B_{1}.
    \end{equation}

    \item For each $\sigma\in[-\delta,\delta]$, there is a subset $\Lambda\subset \Lambda_{0}$ such that
        \begin{equation}\label{Cartan1-3}
          \sharp \Lambda^{c} < M
        \end{equation}
        and
        \begin{equation}\label{Cartan1-4}
          \|(R_{\Lambda} A(\sigma) R_{\Lambda})^{-1}\|< B_{2}.
        \end{equation}
    \item
    \begin{equation}\label{Cartan1-5}
    \emph{mes} \left\{\sigma\in[-\delta,\delta]: \|A(\sigma)^{-1}\|> B_{3}\right\}
< 10^{-3} \gamma (1+B_{1})^{-1} (1+B_{2})^{-1}.
    \end{equation}
  \end{enumerate}

  Then, letting
  \begin{equation}\label{Cartan1-6}
    \kappa< (1+B_{1}+B_{2})^{-10 M},
  \end{equation}
  we have
  \begin{equation}\label{Cartan1-7}
  \emph{mes}\left\{\sigma\in[-\frac{\delta}{2}, \frac{\delta}{2}]:
  \|A(\sigma)^{-1}\|>\frac{1}{\kappa}\right\}
< \emph{\textbf{exp}}\left\{-\frac{c \log \kappa^{-1}}{M \log (M+B_{1}+B_{2}+B_{3})}\right\}.
  \end{equation}
\end{lemma}

\bigskip

\begin{corollary}\label{Cartan2}
Let $T(\sigma)$ be a matrix valued function defined on $\sigma\in[-\delta,\delta]$
  with $T(\sigma)(m,n)\in\CC$ for
  $m,n \in\Lambda_{0}\subset [-N,N]^{d}\subset\ZZ^{d}$, $|\Lambda_{0}|=N$.
  Let the various constants below satisfy
  $$0<\alpha, b, \beta, \rho,  \theta <1,\quad 0<\tau<\frac{9}{10(1+d)}\beta, \quad C>1 $$
   and further let
  $$(\log N)^{2}< M< N^{\tau}.$$
  Assume
  \begin{enumerate}[(i)]
    \item $T(\sigma)$ is real analytic in $\sigma$, and exhibits the off diagonal decay
    \begin{equation}\label{Cartan2-1}
      |T(\sigma)(m,n)|< e^{-\rho |m-n|}, \quad m\neq n.
    \end{equation}
    Moreover, there is a holomorphic
    extension to a strip
    \begin{equation}\label{Cartan2-4}
      |\emph{\textbf{Re}}~z |< \delta, \quad |\emph{\textbf{Im}}~ z|<\gamma
    \end{equation}
    satisfying
    \begin{equation}\label{Cartan2-2}
      \|T(z)\|< N^{C}.
    \end{equation}

    \item For each $\sigma\in[-\delta,\delta]$, the set
    $\Omega(\sigma)$ of bad sites satisfies
    \begin{equation}\label{Cartan2-3}
      \# \Omega(\sigma)< N^{1-\beta}.
    \end{equation}
    Here we say $m\in \Lambda_{0}$ is a  good site if
    $Q_{M}(m)\cap \Lambda_{0}$ and the restriction of $T(\sigma)$
    on $Q=Q_{M}(m)$ is invertible. Also
    \begin{equation}\label{Cartan2-5}
      \|(R_{Q} T(\sigma) R_{Q})^{-1}\|< e^{M^{b}},
    \end{equation}
    and
    \begin{equation}\label{Cartan2-6}
      |(R_{Q} T(\sigma) R_{Q})^{-1}(x,y)|<e^{-\alpha |x-y|},\quad x,y\in Q_{m}(M), |x-y|
> M^{\theta}.
    \end{equation}
    Otherwise $m$ is called a bad site\footnote{Note that for those sites at
    the corner of $\Lambda_{0}$, the size of $Q_{M}\cap \Lambda_{0}$ might
    have very small diameter. For that reason, we think of all corners of $\Lambda_{0}$ as bad  sites.}.
    \item
    \begin{equation*}
      \textrm{mes}\left\{\sigma\in[-\delta,\delta]:
      \|T(\sigma)^{-1}\|> e^{N^{\frac{\beta}{4}}}
      \right\}< e^{-10 M}.
    \end{equation*}
  \end{enumerate}

  Then, we have
  \begin{equation}\label{Cartan2-7}
    \textrm{mes}\left\{\sigma\in[-\frac{\delta}{2},\frac{\delta}{2}]:
      \|T(\sigma)^{-1}\|> e^{N^{1-\frac{\beta}{10}}}
      \right\}< e^{-N^{\frac{\beta}{20}}}.
  \end{equation}
\end{corollary}

\noindent\textbf{Proof.}~Fix $\sigma\in[-\delta,\delta]$.
Consider a paving of $\Lambda_{0}$
by $Q_{M}(x)$ with $x\in 2 M \ZZ^{d}$.
Let
\begin{equation*}
  \Lambda=\bigcup_{Q_{M}(x)\cap \Omega(\sigma)=\emptyset} Q_{M}(x)\cap \Lambda_{0}.
\end{equation*}
From \eqref{Cartan2-3}, we have
\begin{equation*}
  \# \Lambda^{c} < 2^{d} M^{d} N^{1-\beta}.
\end{equation*}
By Lemma \ref{CL1}, we get a norm control on $G_{\Lambda}=(R_{\Lambda}T(\sigma)R_{\Lambda})^{-1}$
\begin{equation*}
  \|G_{\Lambda}\|< 2 N^{d} e^{M^{b}}.
\end{equation*}

Next we employ \ref{Cartan1} to prove the corollary. Obviously,
$B_{1}=N^{C}$.
Then Lemma \ref{Cartan1} (i) holds with $B_{2}= 2 N^{d} e^{M^{b}}$.
Letting $B_{3}= e^{N^{\frac{\beta}{4}}}$, we have
\begin{equation*}
  e^{-10 M}< 10^{-3} \gamma B_{1}^{-1} B_{2}^{-1}\sim \gamma N^{-C} e^{-M^{b}}
\end{equation*}
and thus Lemma \ref{Cartan1} (iii) holds.
Noticing that
\begin{equation*}
  \kappa= e^{-N^{1-\frac{\beta}{10}}}
< (B_{1}+B_{2})^{-10\cdot 2^{d} M^{d} N^{1-\beta}}\sim e^{-N^{(b+d)\tau+1-\beta}},
\end{equation*}
then the conclusion \eqref{Cartan2-7} is an immediate result of
\eqref{Cartan1-7} as long as $N$ is large enough.
\qed

\subsection{Multi-scale analysis}

\begin{lemma}\label{MA}
  Let $T(\sigma)$ be a matrix valued function defined on $\sigma\in[-\delta,\delta]$
  with $T(\sigma)(m,n)\in\CC$ for
  $m,n \in\Lambda_{0}, \Lambda_{0} \in\mathcal{ER}(N)$.
  Let the various constants below satisfy
  \begin{equation*}
    0<\beta\ll 1,\quad \alpha>0,\quad \rho>0,\quad 0<1-\frac{\beta}{10}<b\leq \theta<1,
  \end{equation*}
  and further let
  \begin{equation}\label{MA-1}
    M= [N^{\beta^{6}}],\quad L_{0}= [N^{\frac{\beta^{2}}{100}}].
  \end{equation}

  Assume the following properties hold.
  \begin{enumerate}[(i)]
  \item $T(\sigma)$ is real analytic in $\sigma$ and satisfies
   \eqref{Cartan2-1},\eqref{Cartan2-4} and \eqref{Cartan2-2}.

  \item  For any $I\in\mathcal{ER}(L_{0})$, except for $\sigma$ in a
  set $\mathscr{E}(I)$ of measure at most $e^{-L_{0}^{\beta^{3}}}$,
  \begin{equation}\label{MA-2}
    \|(R_{I} T(\sigma) R_{I})^{-1}\|< e^{L_{0}^{b}}
  \end{equation}
  and
  \begin{equation}\label{MA-3}
    |(R_{I} T(\sigma) R_{I})^{-1}(m,n)|< e^{-\alpha(L_{0})~ |m-n|}\quad\textrm{for}
    ~ m,n\in I, |m-n|>L_{0}^{\theta}.
  \end{equation}

  \item
  Define again $\Omega(\sigma)$ the set of bad sites in $\Lambda_{0}$ by
  condition \eqref{Cartan2-5} and \eqref{Cartan2-6}.
  Assume further that for any $J\in\mathcal{ER}(L)$ such that $$L> N^{\frac{\beta}{5}},$$
  we have
  \begin{equation}\label{MA-4}
    \# (J\cap \Omega(\sigma))< L^{1-\beta}.
  \end{equation}

  \end{enumerate}

  Then we have
  \begin{equation*}
    \|G_{\Lambda_{0}}\|< e^{N^{b}}
  \end{equation*}
  and
  \begin{equation*}
    |G_{\Lambda_{0}}(m,n)|< e^{-\alpha' |m-n|},\quad\textrm{for all}\quad
    m,n\in\Lambda_{0}, |m-n|>  N^{\theta}
  \end{equation*}
  except for $\sigma\in[-\frac{\delta}{2}, \frac{\delta}{2}]$ in a
  set of measure at most $e^{-N^{c \beta^{2}}}$, where $0<c=c(d)<1$ is an absolute constant.
  Moreover,  the decay rate $\alpha'>(\alpha\wedge\rho)-(\log N)^{-8}$.

\end{lemma}

\noindent\textbf{Proof.}~
Let
\begin{equation}
  \mathscr{E}_{0}=\bigcup_{I\in\mathcal{ER}(L_{0})} \mathscr{E}(I).
\end{equation}
It follows that
\begin{equation*}
  \textrm{mes}~\mathscr{E}_{0}< N^{d} e^{- L_{0}^{\beta^{3}}}.
\end{equation*}
For any $\sigma\in [-\delta,\delta]\setminus \mathscr{E}$, all
$L_{0}$-regions in $\Lambda_{0}$ are good in the sense
of \eqref{MA-2} and \eqref{MA-3}.
Using Lemma \ref{CL1} and  taking
\[
M\equiv L_{0},\quad N\equiv L,\quad \tau\equiv \frac{\beta}{20},
\]
we obtain that for any $J\in\mathcal{ER}(L), J\subset\Lambda_{0}$ with $L>N^{\frac{\beta}{5}}$, there is
\begin{equation*}
  \|G_{J}(\sigma)\|<2  L^{d} e^{L_{0}^{b}}<e^{L^{\frac{\beta}{4}}}.
\end{equation*}

In other words, for any such $J$, we have
\begin{equation*}
\textrm{mes}~\{\sigma\in[-\delta,\delta]: \|G_{J}(\sigma)\|>  e^{L^{\frac{\beta}{4}}}\}
< N^{d} e^{- L_{0}^{\beta^{3}}}< e^{-10 M}.
\end{equation*}
Combining \eqref{MA-4} and taking
\[
M\equiv M, \quad N\equiv L,\quad \tau\equiv \beta^{4},
\]
 it follows from the Cartan's estimate Corollary
\ref{Cartan2}
that
\begin{equation*}
  \textrm{mes}~\{\sigma\in[-\delta/2,\delta/2]: \|G_{J}(\sigma)\|>  e^{L^{1-\frac{\beta}{10}}}\}< e^{L^{-\frac{\beta}{20}}}.
\end{equation*}
Denoting
\begin{equation}\label{102}
  \mathscr{E}= \bigcup_{J\in\mathcal{ER}(L), L>N^{\alpha/5}} \{\sigma\in[-\delta/2,\delta/2]: \|G_{J}(\sigma)\|
> e^{L^{1-\frac{\beta}{10}}}\},
\end{equation}
then $\mathscr{E}$  is the desired exceptional set satisfying
\begin{equation*}
  \textrm{mes}~\mathscr{E}< N^{d+1} e^{L^{-\frac{\beta}{20}}}< e^{-N^{c \beta^{2}}}
\end{equation*}
with $0<c<\frac{1}{100}$ depending on $d$.

Fix any $\sigma\in\mathscr{E}^{c}=[-\delta,\delta]\setminus \mathscr{E}$ in what follows. We shall apply Lemma \ref{CL2} to prove the results.

Observe first that, for $\Lambda\in\mathcal{ER}(L)$, $L>N^{\frac{\beta}{5}}$ with $\Lambda\subset \Lambda_{0}$, it follows from
\eqref{102} that
\begin{equation*}
  \|G_{\Lambda}(\sigma)\|< e^{L^{1-\frac{\beta}{10}}}< e^{L^{b}}.
\end{equation*}

Next, for $I\in\mathcal{ER}(N^{\frac{\beta}{5}})$, $I\subset \Lambda_{0}$
with  $I\cap \Omega(\sigma)=\emptyset$, applying Lemma \ref{CL1} by taking
\begin{equation*}
  M\equiv M,\quad N\equiv N^{\frac{\beta}{5}},\quad \tau=\beta^{4},
\end{equation*}
we get
\begin{equation*}
  |G_{I}(\sigma)(x,y)|< e^{-\tilde{\alpha} |x-y|},\quad \textrm{for}~
  |x-y|> N^{\frac{\beta}{5}\cdot \theta}, x,y\in I,
\end{equation*}
where
\begin{equation*}
  \tilde{\alpha}>(\alpha\wedge\rho)-(\log N^{\beta/5})^{-50}>(\alpha\wedge\rho)
  -(\log N)^{-10}.
\end{equation*}
As a result, we call the above region $I$ is good and call
a $N^{\frac{\beta}{5}}$-region bad if it contains a bad site in $\Omega(\sigma)$.

Finally recalling \eqref{MA-4}, there are at most
$N^{1-\alpha}$ bad M-sites in $\Lambda_{0}$.
The set $\mathcal{F}$ of disjoint bad $N^{\frac{\beta}{5}}$-region satisfies
\begin{equation*}
  \# \mathcal{F}< N^{1-\beta}< \frac{N^{b}}{N^{\frac{\beta}{5}}}
\end{equation*}
since $b>1-\frac{\beta}{10}$. Then the conclusion follows from  Lemma \ref{CL2} by taking
$$ M_{0}\equiv N^{\frac{\beta}{5}},\quad N\equiv N,\quad \tau\equiv \frac{\beta}{5}.$$
The arithmetical condition \eqref{103} is valid since
$\theta\geq b>1-\frac{\beta}{10}$ and $\beta\ll 1$.

The decay rate $\alpha'$ satisfies
\begin{equation*}
  \alpha'> \tilde{\alpha}-N^{-\delta}>\alpha\wedge\rho-(\log N)^{-8}.
\end{equation*}
\qed

\subsection{Large deviation theorem}


Consider the matrix
 \begin{equation}
 T^{\sigma}= D^{\sigma}+ \varepsilon S
 \end{equation}
where $D^{\sigma}$ is a diagonal matrix with
\begin{equation}\label{D sigma}
  D^{\sigma}_{\pm,j,k}= \pm (\langle k,  \lambda'\rangle+\sigma)-\mu_{j}.
\end{equation}
 $S$ satisfies the Toeplitz property
 and $\|S\|<1$. The one dimensional parameter
$\sigma$ is defined on some open set $\mathscr{J}\subset\RR$.
Let $0<\beta\ll 1$ and $0<1-\frac{\beta}{10}<b<\theta<1$.

\underline{\textbf{Assume}} that
\begin{equation}\label{S decay}
  |S(x,y)|< e^{-\rho |x-y|},\quad\textrm{for some}~ \rho>0.
\end{equation}


Obviously, for any large $N$,  $T_{N}^{\sigma}$ has a holomorphic extension of $\sigma$ on
$\mathscr{J}$ to the complex domain
\begin{equation*}
  \{\sigma\in\CC: \textrm{dist}(\sigma, \mathscr{J})<1\}
\end{equation*}
such that
\begin{equation*}
  \|T^{\sigma}_{N}\|< N^{C}.
\end{equation*}
Note that $\sup_{\sigma\in\mathscr{J}}|\sigma|\sim N$. Otherwise, for $|\sigma|> 100 N$,
the matrix $T^{\sigma}_{N}$ is diagonal dominated and
a simple application of Neumann series yields
a desired Green's function estimate of $G_{N}^{\sigma}$.

\underline{\textbf{Assume}} $N_{0}$ is sufficiently large and
 the property
\begin{equation}
"N_{0}-\textrm{good}":\left\{
\begin{aligned}
&\|G^{\sigma}_{N_{0}}\|< e^{N_{0}^{b}},\\
& |G_{N_{0}}^{\sigma}(x,y)|< e^{-\alpha_{0} |x-y|},\quad
  \textrm{for}~|x-y|> N_{0}^{\theta}, |x|\leq N_{0}, |y|\leq N_{0}
  \end{aligned}
  \right.
\end{equation}
holds for all $\sigma$ except in a set $\mathscr{E}_{0}$ of measure
at most $ e^{-N_{0}^{\beta^{3}}}$.

Indeed, for low scale $N_{0}$, the "$N_{0}$-good" property can be derived
from a simple application of Neumann series. We show some details here.
Consider
\begin{equation*}
  |D^{\sigma}_{\pm,j,k}|=|\langle k,  \lambda'\rangle+\sigma\pm\mu_{j}|< \varepsilon_{1},
\end{equation*}
which is valid for $\sigma$ lying in an interval of size $2\varepsilon_{1}$.
Then, denoting
\begin{equation*}
  \mathscr{E}_{0}= \left\{\sigma\in\RR: \min_{1\leq j\leq d, |k|\leq N_{0}}
  |D^{\sigma}_{\pm,j,k}|< \varepsilon_{1}
  \right\},
\end{equation*}
there is
\begin{equation*}
  \textrm{mes}~\mathscr{E}_{0} < 8d  N_{0}^{d} \varepsilon_{1}.
\end{equation*}
For $\sigma\in\mathscr{J}\setminus \mathscr{E}_{0}$,
\begin{equation*}
  \|(D^{\sigma}_{N_{0}})^{-1}\|< \frac{1}{\varepsilon_{1}}.
\end{equation*}
Assume
\begin{equation}
  0<\varepsilon< e^{-4 \rho N_{0}^{\theta}},\quad \varepsilon_{1}\sim e^{-N_{0}^{b}}.
\end{equation}
By Lemma \ref{Neumann}, we obtain that
\begin{equation*}
  \|G_{N_{0}}^{\sigma}\|=\|(T^{\sigma}_{N_{0}})^{-1}\|< \frac{2}{\varepsilon_{1}}.
\end{equation*}
and
\begin{equation*}
  |G_{N_{0}}^{\sigma}(x,y)|<  e^{-\rho |x-y|}.
\end{equation*}
To ensure that
\begin{equation*}
  \textrm{mes}~ \mathscr{E}_{0}< e^{-N_{0}^{\gamma}}
\end{equation*}
for some $0<\gamma<1$, we take
\begin{equation*}
 1- \frac{\beta}{10}<\gamma< b.
\end{equation*}
Consequently, the "$N_{0}$-good" property
holds (with $\alpha_{0}=\rho$) for all $\sigma$ except in a set $\mathscr{E}_{0}$ of measure
at most $e^{-N_{0}^{\gamma}}< e^{-N_{0}^{\beta^{3}}}$.
Observe also that the matrix element of $T^{\sigma}$ is at most
linear in $\sigma$. Hence $\mathscr{E}_{0}(N_{0})$ is a semi-algebraic
set in $\sigma$ of degree at most $N_{0}^{C(d)}$.

Let $N_{0}\gg 1$ and let
\begin{equation*}
  \underline{N}_{0}=N_{0}^{100\beta^{4}},\quad \underline{N}_{0}<\overline{N}_{0}< \underline{N}_{0}^{C_{*}}
\end{equation*}
where  $C_{*}$ is to be determined later.

We apply Lemma \ref{MA} to get the Green's function estimate at larger scales.
Following the notations in Lemma \ref{MA}, we take
\begin{equation*}
  L_{0}\in [\underline{N}_{0}, \overline{N}_{0}]
\end{equation*}
and define
\begin{equation*}
N\equiv [L_{0}^{100/\beta^{2}}],\quad M\equiv L_{0}^{100\beta^{4}}(= N^{\beta^{6}}).
\end{equation*}
where $0<\beta\ll 1$ is a fixed constant.

For any $\Lambda_{0}\in\mathcal{ER}(N)$, we establish the
Green's function estimate on $G^{\sigma}_{\Lambda_{0}}$.
For any $I\in\mathcal{ER}(L_{0})$ and $I\subset\Lambda_{0}$, it follows
from the previous arguments
that, if $\varepsilon<e^{-4\rho \overline{N}_{0}^{\theta}}$, then
 $G_{I}^{\sigma}=G^{\sigma+ \langle k_{*},\lambda'\rangle}_{I'}$
satisfies the "$L_{0}$-good" property
except $\sigma+\langle k_{*},\lambda'\rangle\in \mathscr{E}_{0}$,
where $k_{*}+I'=I$, $I'\subset [-L_{0},L_{0}]^{d}$
and $I'\in\mathcal{ER}(L_{0})$.
The exceptional set $\mathscr{E}(I)$ is characterized by
\begin{equation*}
  \mathscr{E}(I)= \mathscr{E}_{0}(L_{0})- \langle k_{*},\lambda'\rangle
\end{equation*}
and thus
\begin{equation*}
  \textrm{mes}~\mathscr{E}(I)< e^{-L_{0}^{\gamma}}<e^{-L_{0}^{\beta^{3}}}.
\end{equation*}
This verifies conditions \eqref{MA-2} and \eqref{MA-3}.

For any $J\in \mathcal{ER}(L)$ and any $N>L>N^{\beta/5}$, we compute
$\# (J\cap \Omega(\sigma))$, where $\Omega(\sigma)$
is the set of the $M$-bad sites in $\Lambda_{0}$.
Roughly speaking,
\begin{equation*}
\begin{aligned}
  n\in\Omega(\sigma)\quad \Leftrightarrow\quad  & (n+[-M,M]^{d})\cap \Lambda_{0}~\textrm{is a}~M~\textrm{-bad region}\\
  \Leftrightarrow\quad  & \sigma+ \langle n,\lambda'\rangle\in\mathscr{E}_{0}(M),
  \end{aligned}
\end{equation*}
and a site $n\in\Lambda_{0}$ is taken as bad site whenever $(n+[-M,M]^{d})\cap \Lambda_{0}$ is not an elementary region.

Recall that $\mathscr{E}_{0}(M)$ is a semi-algebraic set of degree at most $M^{C(d)}$.
The number of connected components of $\mathscr{E}_{0}(M)$ does not
exceed $M^{C}$. The constant $C=C(d)$ might differ from line to line.
Moreover, the size of  each component of $\mathscr{E}_{0}(M)$ is less than
$\eta=e^{-M^{\gamma}}$.
Fix a component $[a-\frac{\eta}{2}, a+\frac{\eta}{2}]$ and consider the set
\begin{equation*}
  H= \{n\in J: |\sigma+\langle n,\lambda'\rangle-a|<\eta/2\}.
\end{equation*}
For two different $n, n'\in H$, we have
\begin{equation*}
  |\langle n-n',\lambda'\rangle|<\eta.
\end{equation*}

\underline{\textbf{Assume}} $\lambda'$ is diophantine
\begin{equation*}
  |\langle k,\lambda'\rangle|> \frac{\nu}{|k|^{\tau}},\quad 0\neq k\in\ZZ^{d},
  0<\nu<1, \tau>d+1.
\end{equation*}
Then
\begin{equation*}
  |n-n'|>\nu (\frac{1}{\eta})^{1/\tau}=
  \nu e^{M^{\gamma}/\tau}\gg N=M^{1/\beta^{6}}.
\end{equation*}
whenever $N_{0}$ is large.
As a result, we have
\begin{equation*}
  \# (J\cap \Omega(\sigma)) < M^{C}= N^{C \beta^{6}}<
  (N^{\frac{\beta}{5}})^{5 C \beta^{5}}<  L^{1-\beta}.
\end{equation*}
and this verifies \eqref{MA-4}.

Let
\begin{equation*}
  \underline{N}_{1}= \underline{N}_{0}^{100/\beta^{2}},\quad
  \overline{N}_{1}= \overline{N}_{0}^{100/\beta^{2}}.
\end{equation*}
By Lemma \ref{MA}, we have that
for any $\underline{N}_{1}<N_{1}<\overline{N}_{1}$ and any
$\Lambda_{0}\in \mathcal{ER}(N_{1})$, the property
\begin{equation}
"N_{1}-\textrm{good}":\left\{
\begin{aligned}
&\|G^{\sigma}_{\Lambda_{0}}\|< e^{N_{1}^{b}},\\
& |G_{\Lambda_{0}}^{\sigma}(x,y)|< e^{-\alpha_{1} |x-y|},\quad
  \textrm{for}~|x-y|> N_{1}^{\theta}, x, y\in\Lambda_{0}
  \end{aligned}
  \right.
\end{equation}
except for $\sigma$ is a set $\mathscr{E}_{1}=\mathscr{E}_{1}(N_{1})$ of measure
\begin{equation*}
  \mathscr{E}_{1}(N_{1})
< e^{-N_{1}^{c \beta^{2}}} N_{1}^{d}< e^{-N_{1}^{\beta^{3}}},
\end{equation*}
where $\alpha_{1}=(\alpha_{0}\wedge\rho) -(\log \underline{N}_{1})^{-8}$.

To iterate on, we impose the condition that
\begin{equation*}
  \underline{N}_{1}< \overline{N}_{0}
\end{equation*}
which results in
\begin{equation*}
  \frac{100}{\beta^{2}}< C_{*}.
\end{equation*}

We only write out the iteration statement, whose proof is essential the same
to that from the scale $N_{0}$ to  $N_{1}$. The following statement holds
for all $k\geq 0$.

For any $\underline{N}_{k-1}^{100/\beta^{2}}=\underline{N}_{k}
<N_{k}<\overline{N}_{k}=\overline{N}_{k-1}^{100/\beta^{2}}$ and any
$\Lambda_{0}\in \mathcal{ER}(N_{k})$, the property
\begin{equation}
"N_{k}-\textrm{good}":\left\{
\begin{aligned}
&\|G^{\sigma}_{\Lambda_{0}}\|< e^{N_{k}^{b}},\\
& |G_{\Lambda_{0}}^{\sigma}(x,y)|< e^{-\alpha_{k} |x-y|},\quad
  \textrm{for}~|x-y|> N_{k}^{\theta}, x, y\in\Lambda_{0}
  \end{aligned}
  \right.
\end{equation}
holds except for $\sigma$ is a set $\mathscr{E}_{k}=\mathscr{E}_{k}(N_{k})$ of measure
\begin{equation*}
  \mathscr{E}_{k}< e^{-N_{k}^{\beta^{3}}},
\end{equation*}
where $\alpha_{k}=(\alpha_{0}\wedge\rho)-(\log \underline{N}_{1})^{-8}-\cdots-(\log \underline{N}_{k})^{-8}$ .
One easily finds that
\begin{equation*}
  \varliminf_{k\rightarrow \infty} \alpha_{k}> \rho-(\log N_{0})^{-8}.
\end{equation*}

Since $\underline{N}_{k+1}< \overline{N}_{k}$, we are able to
iterate constantly and proves Lemma \ref{LDT}.

\bibliographystyle{plain}
\def\cprime{$'$} \def\cprime{$'$} \def\cprime{$'$}

\begin{thebibliography}{10}

\bibitem{Arn61}
V.~I. Arnol{\cprime}d.
\newblock Small denominators. {I}. {M}apping the circle onto itself.
\newblock {\em Izv. Akad. Nauk SSSR Ser. Mat.}, 25:21--86, 1961.

\bibitem{BMS76}
N.~N. Bogoljubov, Ju.~A. Mitropoliskii, and A.~M. Samo{\u\i}lenko.
\newblock {\em Methods of accelerated convergence in nonlinear mechanics}.
\newblock Hindustan Publishing Corp., Delhi; Springer-Verlag, Berlin-New York,
  1976.
\newblock Translated from the Russian edition (published in 1969) by V. Kumar and edited by I. N. Sneddon.

\bibitem{BW08}
J.~Bourgain and W.-M. Wang.
\newblock Quasi-periodic solutions of nonlinear random {S}chr\"{o}dinger
  equations.
\newblock {\em J. Eur. Math. Soc. (JEMS)}, 10(1):1--45, 2008.

\bibitem{Bou97_MRL}
Jean Bourgain.
\newblock On {M}elnikov's persistency problem.
\newblock {\em Math. Res. Lett.}, 4(4):445--458, 1997.

\bibitem{Bou98_Ann}
Jean Bourgain.
\newblock Quasi-periodic solutions of {H}amiltonian perturbations of 2{D}
  linear {S}chr\"odinger equations.
\newblock {\em Ann. of Math. (2)}, 148(2):363--439, 1998.

\bibitem{Bou05_Green}
Jean Bourgain.
\newblock {\em Green's function estimates for lattice {S}chr\"odinger operators
  and applications}, volume 158 of {\em Annals of Mathematics Studies}.
\newblock Princeton University Press, Princeton, NJ, 2005.

\bibitem{BGS02}
Jean Bourgain, Michael Goldstein, and Wilhelm Schlag.
\newblock Anderson localization for {S}chr\"{o}dinger operators on {$\bold
  Z^2$} with quasi-periodic potential.
\newblock {\em Acta Math.}, 188(1):41--86, 2002.

\bibitem{CW93}
Walter Craig and C.~Eugene Wayne.
\newblock Newton's method and periodic solutions of nonlinear wave equations.
\newblock {\em Comm. Pure Appl. Math.}, 46(11):1409--1498, 1993.

\bibitem{Eli88}
L.~H. Eliasson.
\newblock Perturbations of stable invariant tori for {H}amiltonian systems.
\newblock {\em Ann. Scuola Norm. Sup. Pisa Cl. Sci. (4)}, 15(1):115--147
  (1989), 1988.

\bibitem{Kol54}
A.~N. Kolmogorov.
\newblock On conservation of conditionally periodic motions for a small change
  in {H}amilton's function.
\newblock {\em Dokl. Akad. Nauk SSSR (N.S.)}, 98:527--530, 1954.

\bibitem{Kuk93}
Sergej~B. Kuksin.
\newblock {\em Nearly integrable infinite-dimensional {H}amiltonian systems},
  volume 1556 of {\em Lecture Notes in Mathematics}.
\newblock Springer-Verlag, Berlin, 1993.

\bibitem{Mel68}
V.~K. Mel\cprime~nikov.
\newblock A certain family of conditionally periodic solutions of a
  {H}amiltonian system.
\newblock {\em Dokl. Akad. Nauk SSSR}, 181:546--549, 1968.

\bibitem{Mos62}
J.~Moser.
\newblock On invariant curves of area-preserving mappings of an annulus.
\newblock {\em Nachr. Akad. Wiss. G\"{o}ttingen Math.-Phys. Kl. II},
  1962:1--20, 1962.

\bibitem{Pos89}
J\"{u}rgen P\"{o}schel.
\newblock On elliptic lower-dimensional tori in {H}amiltonian systems.
\newblock {\em Math. Z.}, 202(4):559--608, 1989.

\bibitem{Pos01}
J\"{u}rgen P\"{o}schel.
\newblock A lecture on the classical {KAM} theorem.
\newblock In {\em Smooth ergodic theory and its applications ({S}eattle, {WA},
  1999)}, volume~69 of {\em Proc. Sympos. Pure Math.}, pages 707--732. Amer.
  Math. Soc., Providence, RI, 2001.

\bibitem{You99}
Jiangong You.
\newblock Perturbations of lower-dimensional tori for {H}amiltonian systems.
\newblock {\em J. Differential Equations}, 152(1):1--29, 1999.

\end{thebibliography}

\def\cprime{$'$} \def\cprime{$'$} \def\cprime{$'$}

\end{document}